\def\updated{January 18, 2010}

\def\yr{2010}\def\vol{5}

\count100= 1 \count101= 91

\input colordvi



\font\titlefont=cmr17 \font\titlei=cmmi10 at 17pt
\font\titlesy=cmsy10 at 17pt \font\titleit=cmti10 at 17pt
\font\titlesl=cmsl10 at 17pt \font\titlebf=cmbx10 at 17pt
\font\Bbbt=msbm10 at 17pt \font\twelverm=cmr12 \font\twelvebf=cmbx10
at 12pt \font\twelvei=cmmi10 at 12pt \font\twelvesy=cmsy10 at 12pt
\font\ninerm=cmr9 \font\sevenrm=cmr7 \font\sixrm=cmr6
\font\fiverm=cmr5 \font\ninei=cmmi9 \font\seveni=cmmi7
\font\ninesy=cmsy9 \font\sevensy=cmsy9 
\font\sixbf=cmbx6 \font\fivebf=cmbx5 \font\ninebf=cmbx9
\font\nineit=cmti9 \font\ninesl=cmsl9 \font\nineex=cmex9
\font\nineBbb=msbm9  
  \font\Bbb=msbm10

\def\ninepoint{\def\rm{\fam0\ninerm}
    \textfont0 = \ninerm
    \textfont1 = \ninei
    \textfont2 = \ninesy
    \textfont3 = \nineex
    \scriptfont0 = \sevenrm
    \scriptfont1 = \seveni
    \scriptfont2 = \sevensy
    \scriptscriptfont0 = \fiverm
    \scriptscriptfont1 = \fivei
    \scriptscriptfont2 = \fivesy
    \textfont\itfam=\nineit \def\it{\fam\itfam\nineit}
    \textfont\bffam=\ninebf \scriptfont\bffam=\sixbf
    \scriptscriptfont\bffam=\fivebf \def\bf{\fam\bffam\ninebf}
    \textfont\slfam=\ninesl \def\sl{\fam\slfam\ninesl}
    \let\Bbb\nineBbb
    \baselineskip 10pt}

\def\titlepoint{\def\rm{\fam0\titlefont}
    \textfont0 = \titlefont
    \textfont1 = \titlei
    \textfont2 = \titlesy
    \textfont3 = \nineex
    \scriptfont0 = \twelverm
    \scriptfont1 = \twelvei
    \scriptfont2 = \twelvesy
    \textfont\itfam=\titleit \def\it{\fam\itfam\titleit}
    \textfont\bffam=\titlebf \scriptfont\bffam=\twelvebf
    \def\bf{\fam\bffam\titlebf}
    \textfont\slfam=\titlesl \def\sl{\fam\slfam\titlesl}
    \let\Bbb\Bbbt\titlefont}


\def\endproofsymbol{\makeblanksquare6{.4}}
\def\eop{\endproofsymbol\nopf}

\def\meop{~~\endproofsymbol}

\def\nopf{\medskip\goodbreak}

\def\pf{\noindent{\bf Proof: }}

\def\makeblanksquare#1#2{
\dimen0=#1pt\advance\dimen0 by -#2pt
      \vrule height#1pt width#2pt depth0pt\kern-#2pt
      \vrule height#1pt width#1pt depth-\dimen0 \kern-#1pt
      \vrule height#2pt width#1pt depth0pt \kern-#2pt
      \vrule height#1pt width#2pt depth0pt
}

\magnification\magstephalf

\hsize6.5truein\vsize8.6truein \voffset.5truein


\def\title#1{\toneormore#1||||:}
\def\titexp#1#2{\hbox{{\titlefont #1} \kern-.25em%
  \raise .90ex \hbox{\twelverm #2}}\/}
\def\titsub#1#2{\hbox{{\titlefont #1} \kern-.25em%
  \lower .60ex \hbox{\twelverm #2}}\/}

\def\author#1{\bigskip\bigskip\aoneormore#1||||:\smallskip\centerline{\updated}}

\def\abstract#1{\bigskip\bigskip\medskip%
 {\ninepoint
 \narrower{\bf Abstract.~}\rm#1\smallskip
  MSC: \mscnumbers\ifx\keywords\empty\else\smallskip
  Keywords: \keywords\fi\bigskip
 \printtochere}\starttoc\bigskip}

\def\toneormore#1||#2||#3:{\centerline{\titlepoint #1}%
    \def\next{#2}\ifx\next\empty\else\medskip\toneormore#2||#3:\fi}
\def\aoneormore#1||#2||#3:{\centerline{\twelverm #1}%
    \def\next{#2}\ifx\next\empty\else\smallskip\aoneormore#2||#3:\fi}

\newwrite\toc\def\tocone{0}\def\tochalf{.5}\def\toctwo{1}
\countdef\counter=255
\def\diamondleaders{\global\advance\counter by 1
  \ifodd\counter \kern-10pt \fi
  \leaders\hbox to 15pt{\ifodd\counter \kern13pt \else\kern3pt \fi
    .\hss}}
\newdimen\lextent
\def\printtochere{\immediate\closeout\toc\relax%
\begingroup
\def\\##1.  ##2.  {\setbox1=\hbox{##1}\ifnum\wd1>\lextent\lextent\wd1\fi}
\lextent0pt\inputifthere{\jobname.toc}\advance\lextent by 2em\relax
\def\\##1.  ##2.  {\centerline{\hbox to \lextent{\rm##1\def\next{##2}%
\iffoliointoc \ifx\next\empty\else\diamondleaders\fi\hfil\hbox to
2em{\hss##2}\fi}}}
\inputifthere{\jobname.toc}\endgroup}

\def\starttoc{\immediate\openout\toc=\jobname.toc}
\def\nexttoc#1{{\let\folio=0\edef\next{\write\toc{#1}}\next}}

\def\tocline#1#2#3{\nexttoc{\noexpand\noexpand\noexpand\\\hskip#2truecm #1.  #3.  }}

\def\footnoterule{\kern -3pt \hrule width 0truein \kern 2.6pt}
\def\leftheadline{\ifnum\pageno=\count100 \hfill%
  \else\hfil\it\shortauthor\hfil\llap{\rm\folio}\fi}
\def\rightheadline{\ifnum\pageno=\count100 \hfill%
  \else\hfil\it\shorttitle\hfil\llap{\rm\folio}\fi}

\nopagenumbers \headline{\ifodd\pageno\rightheadline
\else\leftheadline\fi} \footline{\hfil} \null\vskip 18pt
\centerline{} \pageno=\count100 \count102=\count100
\advance\count102 by -1 \advance\count102 by \count101


\def\copyright{\hbox{{\twelverm o}\kern-.61em\raise .46ex\hbox{\fiverm c}}}

\insert\footins{\sixrm
\medskip
\baselineskip 8pt \leftline{Surveys in Approximation Theory
  \hfill {\ninerm \the\pageno}}
\leftline{Volume \vol, \yr, pp.~\the\pageno--\the\count102.}
\leftline{Copyright \copyright\ \yr\ Surveys in Approximation
Theory.} \leftline{ISSN 1555-578X} \leftline{All rights of
reproduction in any form reserved.}
\smallskip
\par\allowbreak}


\def\sect#1{\startsect\edef\showsectno{\the\sectionno}%
   \let\tocindent\tocone\soneormore#1||||:\relax\medskip\noindent\ignorespaces}
\def\sectwopn#1{\foliointocfalse\sect{#1}\foliointoctrue}

\def\soneormore#1||#2||#3:{%
   \formsecthead{#1}
   \def\next{#2}\def\pfolio{\folio}\iffoliointoc\else\def\pfolio{}\fi%
   \ifx\next\empty\puttocline{\showsectno\ \ #1}{\pfolio}%
   \else\puttocline{\showsectno\ \ #1}{}\let\showsectno\skipsectno\let\tocindent\tochalf\soneormore#2||#3:\fi}

\def\formsecthead#1{\leftline{\bf\showsectno\hskip2em #1}}

\newif\iffoliointoc\foliointoctrue
\def\puttocline#1#2{\tocline{#1}{\tocindent}{#2}}
\def\skipsectno{\setbox0=\hbox{\the\sectionno}\hskip\wd0}

\def\subsect#1{\bigskip\formsubsecthead{#1}\medskip\let\tocindent\toctwo\puttocline{#1}{\folio}}

\def\formsubsecthead#1{{\bf #1}\hskip1em}

\newcount\sectionno\sectionno0
\def\presect{\the\sectionno.}
\newcount\subsectionno
\def\startsect{\ifx\empty\presect\else\restartnums\fi%
               \subsectionno0\global\advance\sectionno by 1\relax

               \goodbreak\bigskip\smallskip}
\def\startsubsect{\global\advance\subsectionno by 1\goodbreak\bigskip}


\newcount\blackmarks\blackmarks0
\newcount\eqnum
\newcount\labelnum
\def\restartnums{\eqnum0\labelnum0}
\def\singlecount{\let\labelnum\eqnum}

 \newread\testfl
 \def\inputifthere#1{\immediate\openin\testfl=#1
    \ifeof\testfl\message{(#1 does not yet exist)}
    \else\input#1\fi\closein\testfl}

 \inputifthere{\jobname.aux}
 \newwrite\aux
 \immediate\openout\aux=\jobname.aux

\def\plazieres{\expandafter\ifx\csname\griff\endcsname\relax%
  \xdef\esfehlt{\griff}\blackmark\else{\csname\griff\endcsname}\fi}

\def\definieres{\expandafter\xdef\csname\griff\endcsname{\inhalt}%
 \def\blankkk{ }\expandafter\immediate\write\aux{%
 \string\expandafter\def\string\csname%
 \blankkk\griff\string\endcsname{\inhalt}}}

\def\blackmark{\ifnum\blackmarks=0\global\blackmarks=1%
 \write16{============================================================}%
 \write16{Some forward reference is not yet defined. Re-TeX this file!}%
 \write16{============================================================}%
 \fi\immediate\write16{undefined forward reference: \esfehlt}%
 {\vrule height10pt width2pt depth2pt}\esfehlt%
 {\vrule height10pt width2pt depth2pt}}

\def\marginal#1{\strut\vadjust{\kern-\strutdepth%
\vtop to \strutdepth{\baselineskip\strutdepth\vss\llap{\fiverm#1\
}\null}}}
\def\strutdepth{\dp\strutbox}


\newif\ifdraft

\newcount\hour\newcount\minutes
\def\draft{\drafttrue
\headline={\sevenrm \hfill\ifx\filenamed\undefined\jobname\else\filenamed\fi%
(.tex) (as of \ifx\updated\undefined???\else\updated\fi)
 \TeX'ed at {\hour\time\divide\hour by 60{}%
\minutes\hour\multiply\minutes by 60{}%
\advance\time by -\minutes \the\hour:\ifnum\time<10{}0\fi\the\time\
on \today\hfill}} }

\def\today{\number\day\space%
\ifcase\month\or January\or February\or March\or April\or May\or
June\or
 July\or August\or September\or October\or November\or December\fi%
\space\number\year}




\gdef\formfirstauthor{\the\firstname\  \the\lastname}
\gdef\formnextauthor{, \the\firstname\the\lastname}
\gdef\formotherauthor{ and \the\firstname\the\lastname}
\gdef\formlastauthor{,\formotherauthor}

\gdef\formB{\the\au\ [\yr] ``\the\ti'', \the\pb, \the\pl.
\setcitelabel} \gdef\formD{\the\au\ [\yr] ``\the\ti'', dissertation,
\the\pl. \setcitelabel}
\gdef\formJ{\the\au\ [\yr] \the\ti, {\sl\the\jr}\ifx\vl\empty%
\else\ {\bf\vl}\fi, \pp. \setcitelabel} \gdef\formP{\the\au\ [\yr]
\the\ti, in {\sl\the\tit},
\getfirstchar\aut\ifx\firstchar\unknownx\else\the\aut, ed\edsop, \fi
\getfirstchar\pub\ifx\firstchar\unknownx\else\the\pub, \fi \the\pl,
\pp. \setcitelabel} \gdef\formR{\the\au\ [\yr]
\the\ti\ifx\is\empty\else, \is\fi. \setcitelabel}

\newtoks\lastname
\newtoks\firstname
\newtoks\au
\newtoks\aut
\newtoks\ti
\newtoks\tit
\newtoks\pb
\newtoks\pub
\newtoks\pl
\newtoks\jr

\newtoks\rhlau

\def\setcitelabel{\edef\griff{cit\rh}\edef\inhalt{\the\rhlau\ \yr}\definieres}
\def\setcitelabel{}

\def\getfirstchar#1{\edef\theword{\the#1}\expandafter\getit\theword:}
\def\getit#1#2:{\def\firstchar{#1}}
\def\unknownx{x}

\newif\ifonesofar
\def\concat#1{\edef\audef{{#1}}\au=\audef}
\def\decodeauthor#1, #2,#3;{\lastname={#1}\firstname={#2}%
\concat{\formfirstauthor}\onesofartrue%
\def\morerhlau{}%
\def\next{#3}\ifx\next\empty\else\def\morerhlau{ et al.}\decodemoreauthor#3;\fi
\edef\morerhlauu{{\the\lastname\morerhlau}}\rhlau=\morerhlauu}
\def\decodemoreauthor#1, #2,#3;{\lastname={#1}\firstname={#2}%
\def\next{#3}\ifx\next\empty\let\formaut=\formlastauthor%
\ifonesofar\ifx\formotherauthor\undefined\else\let\formaut=\formotherauthor%
\fi\fi\concat{\the\au\formaut}%
\else\onesofarfalse\concat{\the\au\formnextauthor}\decodemoreauthor#3;\fi}

\def\refproc #1(#2; #3; {\decodeproc#2; \xdef\yr{#3}}
\def\decodeproc#1), #2 (ed#3.),#4 (#5); {%
 \global\tit={#1}\global\aut={#2}\xdef\edsop{#3}\global
 \pub={#4}\global\pl={#5}}

\def\dd{\,{\rm d}} 
\def\ee{{\rm e}}   
\def\ii{{\rm i}}   
\font\Bbb=msbm10 

\def\CC{\hbox{{\Bbb{C}}}}

\def\NN{\hbox{{\Bbb{N}}}}

\def\RR{\hbox{{\Bbb{R}}}}

\title{Strong Uniqueness}
\author{Andr\'as Kro\'o and Allan Pinkus}

\def\shorttitle{Strong Uniqueness}
\def\shortauthor{A.~Kro\'o and A.~Pinkus}

\def\mscnumbers{41A52, 41A50, 41A65}

\def\keywords{strong uniqueness, best approximation, uniqueness}


\def\alp{\alpha}

\def\del{\delta}                
\def\eps{\varepsilon}
                
\def\lam{\lambda}               
\def\sig{\sigma}                
\def\ome{\omega}

\def\bfc{{\bf c}}               
\def\bfm{{\bf m}}               
\def\bfu{{\bf u}}               
\def\bfw{{\bf w}}               
\def\bfx{{\bf x}}

\def\incl{\subseteq}
\def\nek{,\ldots,}

\def\amp{\raise 6pt\hbox{$\scriptstyle\bf1$}}

\def\span{{\rm span}}

\def\sgn{\mathop{\rm sgn}}

\def\\{{\backslash}}

\def\om{{\overline m}}

\def\tilm{{\widetilde m}}
\def\tilx{{\widetilde x}}

\def\tilbfm{{\widetilde{\bf m}}}

\def\sqr#1#2{{\vcenter{\hrule height.#2pt\hbox{\vrule
width.#2pt height#1pt \kern#1pt \vrule width.#2pt}\hrule
height.#2pt}}}

\def\formsecthead#1{\centerline{\bf #1}}

\def\formsubsecthead#1{\centerline{\bf #1}}


\abstract{This is a survey paper on the subject of strong uniqueness
in approximation theory.}


\subsect{{0.} Foreword}

\medskip\noindent
This is a survey paper on the subject of {\sl Strong Uniqueness} in
approximation theory. The concept of strong uniqueness was
introduced by Newman, Shapiro in 1963. They proved, among other
things, that if $M$ is a finite-dimensional real Haar space in
$C(B)$, $B$ a compact Hausdorff space, and if $m^*$ is the best
approximant to $f$ from $M$, then there exists a $\gamma>0$ such
that
$$\|f-m\| - \|f-m^*\| \ge \gamma \|m-m^*\|\eqno{(0.1)}$$
for all $m\in M$, where $\gamma$ may depend upon $f$, $M$ and $m^*$,
but is independent of the specific $m$. An inequality of the above
form valid for all $m\in M$ is what we call {\sl classical strong
uniqueness}. This property is stronger than the uniqueness of the
best approximant (hence the name).

Strong uniqueness has been much studied, mainly during the 1970s,
1980s and 1990s, and there have been over 100 research papers
devoted to the subject. We thought this might be a good time to
reflect upon what was considered and accomplished in the study of
this subject. We hope that you will agree that the resulting theory
is not insignificant.

This paper is organized as follows. We have divided the survey into
three parts. The first part is concerned with what we termed above
classical strong uniqueness, i.e., inequality (0.1).  In Section 1,
we exactly characterize the optimal strong uniqueness constant
(largest possible $\gamma$ in (0.1)) via the one-sided Gateaux
derivative.  In Section 2, we consider some general results
regarding classical strong uniqueness in the uniform norm. We look
in detail at the case where the approximating set is a
finite-dimensional subspace, give a characterization of when we have
strong uniqueness, and also prove that in this setting the set of
functions with a strongly unique best approximant is dense in the
set of functions with a unique best approximant. Finally, we provide
a general upper bound on the optimal strong uniqueness constant
based on projection constants. In Section 3, we consider the
relationship between local Lipschitz continuity of the best
approximation operator at a point and classical strong uniqueness at
the same point. A general result is that the latter always implies
the former. We prove that the converse holds when approximating from
a finite-dimensional subspace in the uniform norm. In Section 4, we
restrict our attention to approximation from finite-dimensional Haar
spaces in the uniform norm. We first prove the result of Newman,
Shapiro [1963] mentioned above, i.e., that in this case, we always
have strong uniqueness. We then consider various properties of the
optimal strong uniqueness constant with respect to the function
being approximated. We prove, for example, that the optimal strong
uniqueness constant is upper semi-continuous but not necessarily
continuous, and it is not uniformly bounded below by a positive
constant (i.e., bounded away from $0$) if the subspace is of
dimension at least $2$ and the underlying domain is not discrete. We
find lower bounds (that sometimes provide equality) for the optimal
strong uniqueness constant via specific elements of the
approximating subspace. Assuming the subspace is of dimension $n$
and the error function only attains its norm at $n+1$ points (the
generic case), we obtain that the optimal strong uniqueness constant
is bounded above by $1/n$. We also give another characterization of
the optimal strong uniqueness constant and consider the question of
when strong uniqueness and uniqueness are equivalent concepts. In
Section 5, we look at classical strong uniqueness in the $L^1$ norm.
Among other results, we prove that if $\nu$ is a non-atomic positive
measure then the set of functions in $L^1(K, \nu)$ that have a
strongly unique best approximant from any finite-dimensional
subspace is dense in $L^1(K, \nu)$. In addition, under the above
assumptions, we show, as in the uniform norm, that a function in
$L^1(K, \nu)$ has a strongly unique best approximant if and only if
the best approximation operator from the same finite-dimensional
subspace is locally Lipschitz continuous. Furthermore, in this case,
we get explicit upper and lower bounds for the optimal strong
uniqueness constant based on the Lipschitz continuity of the best
approximation operator. We also present some results on classical
strong uniqueness in the problem of one-sided $L^1$ approximation.
In Section 6, we consider approximation by rational functions of the
form
$$R_{m,n}:= \{r=p/q: p\in \Pi_m, q\in \Pi_n, q(x)>0, x\in [a,b]\},$$
where $\Pi_n = \span\{1,x\nek x^n\}$. The main result reported on is
that we have the equivalence of strong uniqueness to a function $f$
from $R_{m,n}$, the operator of best rational approximation from
$R_{m,n}$ being continuous at $f$, and the fact that the unique best
rational approximant to $f$ from $R_{m,n}$ is not contained in
$R_{m-1,n-1}$.

In the second part of this survey, we discuss what we call {\sl
non-classical strong uniqueness}. By this, we mean the existence of
a nonnegative strictly increasing function $\phi$  defined on
$\RR_+$, and a constant $\gamma>0$ that may depend upon $f$ and $M$
(and thus on $m^*$), for which
$$\|f-m\| - \|f-m^*\| \ge \gamma \phi(\|m-m^*\|)$$
for all $m\in M$ (a global estimate), or for all $m\in M$ such that
$$\|f-m\| -\|f-m^*\| \le \sig$$
for some $\sig>0$ (a local estimate). We start in Section 7 with the
case of a uniformly convex norm and prove the basic inequality
$$\|f-m\| - \|f-m^*\| \ge  \|f-m\|\, \del\left({{\|m-m^*\|}\over {\|f-m\|}}\right)$$
valid for all $m\in M$ where $m^*\in M$ is a best approximant to
$f$. Here $\del$ is the usual modulus of convexity of the uniformly
convex space. We  apply this result and consider also variants
thereof. In Section 8, we return to a consideration of the uniform
norm. The main result we report on therein is that non-classical
strong uniqueness holds for $\phi(t)=t^2$ if $M\subset C^{2}[a,b]$
is a finite-dimensional unicity space (and not necessarily a Haar
space) with respect to $C^{2}[a,b]$. In Section 9, we return to a
consideration of the $L^1$-norm and obtain two local non-classical
strong uniqueness estimates. In one case, we prove that if $M$ is a
finite-dimensional unicity space for all $f\in C[a,b]$ in the
$L^1[a,b]$ norm then we have non-classical strong uniqueness with
$\phi(t) = t \ome^{-1}(f; D t)$ for some constant $D>0$ that depends
only on $f$ and $M$. Here $\ome$ is the standard modulus of
continuity. See Theorem 9.2 for a more detailed explanation. A
second result concerns the one-sided $L^1$ approximation problem. We
prove that if $M\subset C^1[a,b]$ is a finite-dimensional unicity
subspace for $C^1[a,b]$ in the one-sided $L^1$-norm, then  we have
non-classical strong uniqueness with $\phi(t) = t H_f^{-1}(D t)$ for
some constant $D>0$ that depends only on $f$ and $M$, where $H$ is
based on the moduli of continuity of $f'$ and $m'$ for all $m\in M$.
In Section 10, we consider strong uniqueness when we approximate
complex-valued functions in the uniform norm. The main result
therein is non-classical strong uniqueness with $\phi(t) = t^2$. We
also discuss conditions under which classical strong uniqueness
holds.

The third part of this survey concerns various applications of
strong uniqueness results. The main idea behind these applications
is that instead of solving a best approximation problem in a given
norm, we replace it by considering another norm, close to the
original norm, one that leads to a simpler approximation problem.
Strong uniqueness is then applied in order to show that the best
approximant in this new norm is sufficiently close to the original
best approximant. An estimation of the ``closeness'' is given in
Section 11. Typically, the original norm is modified by replacing it
by a similar discrete norm (as considered in Section 12), or by
introducing a weight function into the norm. In Section 13, we use
strong uniqueness results in order to solve approximation problems
in the case when the norm is altered by a weight function. This
approach is used to derive asymptotic representations for weighted
Chebyshev polynomials.

We have tried to make this survey as self-contained as seemed
reasonable. As such there are many additional peripheral results
included in this paper that, we hope, put the strong uniqueness
results into a reasonable context.

On the other hand, while this survey is lengthy, it is far from
complete. Certain specific problems and general topics within the
area of strong uniqueness have not been surveyed. For example, a lot
of effort and many papers have been concerned with the study of
$\gamma_n(f)$, the optimal strong uniqueness constant (in (0.1))
when approximating $f\in C[a,b]$ in the uniform norm from
$\Pi_n=\span\{1,x\nek x^n\}$. Many different questions were asked
concerning this important case. Nevertheless, we have not considered
this class of problems. Thus, as is noted in Section 4, Poreda
[1976] raised the question of describing the asymptotic behaviour of
the above sequence $(\gamma_n(f))$ as $n$ tends to $\infty$ for a
given function $f$. Henry and Roulier [1978] conjectured that
$$\liminf_{n\to\infty} \gamma_n(f) =0$$
if $f$ is not a polynomial. This conjecture, sometimes called the
Poreda conjecture, was proved for various classes of functions over
a string of papers. It was finally solved in the positive by Gehlen
[1999]. It should be noted that the above $\liminf$ cannot be
replaced by the simple limit. There are functions $f\in C[a,b]$ for
which
$$\liminf_{n\to\infty}\gamma_n(f) =0,\qquad \limsup_{n\to\infty}\gamma_n(f)=1;$$
see Schmidt [1978]. These same $\gamma_n(f)$ were also studied for
specific functions $f$, or specific classes of function; see for
example Henry, Huff [1979], Henry, Swetits [1981], Henry, Swetits,
Weinstein [1981], Henry, Swetits [1982] and  Henry, Swetits,
Weinstein [1983]. There are also papers devoted to looking for sets
in $C[a,b]$ over which the strong uniqueness constants are uniformly
bounded below by a positive constant; see for example Bartelt,
Swetits [1983], Marinov [1983], and Bartelt, Swetits [1988], as well
as papers devoted to looking at the optimal strong uniqueness
constant as a function of the domain; see for example Henry, Roulier
[1977], Bartelt, Henry [1980], and Paur, Roulier [1981]. For a study
of strong uniqueness when approximating with constraints, leading to
non-classical strong uniqueness, see for example Fletcher, Roulier
[1979],  Schmidt [1979], Chalmers, Metcalf, Taylor [1983] and
Kro\'o, Schmidt [1991]. All the above is with regards to
approximation from $\Pi_n$. In addition, strong uniqueness when
approximating by splines, with either fixed or variable knots, was
considered in N\"urnberger, Singer [1982], N\"urnberger [1982/83],
N\"urnberger [1994], Sommer, Strauss [1993], Zeilfelder [1999] and
Zwick [1987]. Again this is not discussed in what follows.

We have tried to give exact references, and have also endeavored to
make the list of references complete. This list contains all
references we have found on the subject of strong uniqueness. Note
that not all papers on the list of references are referred to in the
body of this survey. We apologize for any omissions and would
appreciate information with regard to any additional references.

\goodbreak

\bigskip
\sectwopn{{Part I.} Classical Strong Uniqueness}

\subsect{{1.} Introduction}

\medskip\noindent
Let $X$ be a normed linear space and $M$ a subset of $X$. For any
given $f\in X$, we let $P_M(f)$ denote the set of best approximants
to $f$ from $M$. That is, $m^*\in P_M(f)$ if $m^*\in M$ and
$$\|f-m^*\| \le \|f-m\|$$
for all $m\in M$. This set $P_M(f)$ may be empty (non-existence), a
singleton (unicity) or larger. In this section, we consider
conditions on a linear subspace or convex subset $M$ of $X$ under
which we obtain {\sl Classical Strong Uniqueness}. By this, we mean
the existence of a strictly positive constant $\gamma>0$ for which
$$\|f-m\| - \|f-m^*\| \ge \gamma \|m-m^*\|$$
for all $m\in M$ where $m^*\in P_M(f)$. The constant $\gamma>0$ can
and will depend upon $f$, $m^*$ and $M$, but must be independent of
$m$. Note that the form of this inequality is the best one might
expect since, by the triangle inequality, we always have
$$\|f-m\| -\|f-m^*\| \le \|m-m^*\|.$$

We first exactly characterize the optimal (largest) $\gamma$ for
which such an inequality holds. This will be done via the one-sided
Gateaux derivative that is defined as follows. Given $f,g \in X$,
set
$$\tau_+(f,g) =\lim_{t\to 0^+}{{\|f+tg\|-\|f\|}\over{t}}.\eqno{(1.1)}$$
The value $\tau_+(f,g)$ is termed the one-sided Gateaux derivative
of the norm at $f$ in the direction $g$. The first thing we prove is
that this one-sided Gateaux derivative always exists.

\proclaim Proposition 1.1. The one-sided Gateaux derivative of the
norm at $f$ in the direction $g$, namely $\tau_+(f,g)$, exists for
every $f,g\in X$.

\pf Set
$$r(t) := {{\|f+tg\|-\|f\|}\over{t}}.$$
We claim that $r(t)$ is a non-decreasing function of $t$ on
$(0,\infty)$ and is bounded below thereon. As this is valid then
$\tau_+(f,g)$ necessarily exists.

To see that $r(t)$ is bounded below on $(0, \infty)$, note that
$$\|f+tg\| \ge \|f\| - \|tg\| = \|f\| -t \|g\|.$$
Thus for $t>0$, we have $r(t)\ge -\|g\|$. The non-decreasing
property of $r$ can be shown as follows. Let $0<s<t$. Then
$$t\|f+sg\| = \|tf+tsg\| = \|s(f+tg)+(t-s)f\| \le s\|f+tg\| + (t-s)\|f\|.$$
Thus
$$t(\|f+sg\|-\|f\|) \le s(\|f+tg\| - \|f\|),$$
and that implies $r(s)\le r(t)$. \eop

Using this functional $\tau_+$, it is now a simple matter to
characterize best approximants from linear subspaces. Namely,

\proclaim Theorem 1.2. Let $M$ be a linear subspace of $X$. Then
$m^*\in P_M(f)$ if and only if $\tau_+(f-m^*, m)\ge 0$ for all $m\in
M$.

\pf ($\Rightarrow$) Assume $m^*$ is a best approximant to $f$ from
$M$. Therefore
$$ \|f-m^* + tm\|\ge \|f-m^*\|$$
for every $m\in M$ and $t>0$,  immediately implying that
$\tau_+(f-m^*, m)\ge 0$.

\smallskip\noindent
($\Leftarrow$) Assume $\tau_+(f-m^*, m)\ge 0$ for all $m\in M$. As
$M$ is a linear subspace, this implies that $\tau_+(f-m^*, m^*-m)
\ge 0$. From the proof of Proposition 1.1, $r(t)$ is a
non-decreasing function of $t$ on $(0,\infty)$. Setting $t=1$
therein with respect to $\tau_+(f-m^*, m^*-m)$, we obtain
$$\|f-m^*+m^*-m\| -  \|f-m^*\| \ge  \tau_+(f-m^*, m^*-m)\ge 0.$$
Thus $\|f-m\| \ge \|f-m^*\|$ for all $m\in M$. \eop

If $M$ is a convex subset of $X$, then the same arguments prove the
following.

\proclaim Corollary 1.3. Let $M$ be a convex subset of $X$. Then
$m^*\in P_M(f)$ if and only if $\tau_+(f-m^*, m^*-m)\ge 0$ for all
$m\in M$.

\medskip
If $\tau_+(f-m^*, m^*-m)>0$ for all $m\in M$, $m\ne m^*$, then we
may have a $\gamma >0$ for which
$$ \|f-m \| - \|f-m^*\| \ge \gamma \|m-m^*\| $$
for all $m\in M$. When there exists a $\gamma>0$ for which
$$ \|f-m \| - \|f-m^*\| \ge \gamma \|m-m^*\|$$
for all $m\in M$, then we say that $m^*$ is a {\sl strongly unique}
best approximant to $f$ from $M$. The reason for this terminology is
simply that strong uniqueness is stronger than uniqueness. That is,
if $m^*$ is a strongly unique best approximant to $f$ from $M$, then
it is certainly a unique best approximant. While the converse does
not hold in general, it does and can hold in various settings. This
concept was introduced in Newman, Shapiro [1963] with respect to
certain specific spaces. We will state and prove their results in
later sections.

This next result characterizes the optimal (largest) $\gamma$ for
which such an inequality holds. We prove this result under the
assumption that $M$ is a linear subspace. The case where $M$ is only
a convex subset follows analogously, and we will subsequently
formally state it as a corollary.

\proclaim Theorem 1.4. Let $M$ be a subspace of $X$ and $f\in X\\
{M}$. Assume $m^*\in P_M(f)$. Set
$$\gamma (f) := \inf\{ \tau_+(f-m^*, m): m\in M, \|m\|=1\}.$$
Then $\gamma (f) \ge 0$ and, for all $m\in M$, we have
$$\|f-m \| - \|f-m^*\| \ge \gamma (f)\, \|m-m^*\|.$$
Furthermore, if $\gamma > \gamma (f)$ then there exists an $\tilm\in
M$ for which
$$ \|f-\tilm \| - \|f-m^*\| < \gamma\, \|\tilm-m^*\|.$$

\pf As $\tau_+(f-m^*, m^*-m)\ge 0$ for all $m\in M$ ($M$ is a
subspace), we have $\gamma (f) \ge 0$ from Theorem 1.2. Assume
$\gamma (f)>0$. From the fact that $\tau_+(f,cg) = c\, \tau_+(f,g)$
for all $c>0$, it follows that
$$\tau_+(f-m^*, m^*-m ) \ge \gamma (f) \|m-m^*\|$$
for all $m\in M$. From the proof of Proposition 1.1,
$${{\|f-m^* +t(m^*-m)\| - \|f-m^*\|}\over t}  \ge \tau_+(f-m^*, m^*-m )$$
for all $t>0$. Setting $t=1$, we obtain
$$\|f-m\| - \|f-m^*\|\ge \tau_+(f-m^*, m^*-m )  \ge \gamma (f) \|m-m^*\|$$
for all $m \in M$.

Now assume $\gamma > \gamma (f)$. By definition, there exists an
$\om\in M$, $\|\om \|=1$, for which
$$\tau_+(f-m^*, \om ) < \gamma \|\om\|.$$
Thus for $t>0$ sufficiently small, we have
$${{\|f-m^* +t\om\| - \|f-m^*\|}\over t}  < \gamma \|\om\|,$$
implying
$$\|f-m^* +t\om\| - \|f-m^*\|  < \gamma \|t\om\|.$$
Setting $\tilm := m^* - t\om$, we obtain
$$ \|f-\tilm \| - \|f-m^*\|< \gamma \|\tilm-m^*\| . \meop $$

For a convex set $M$, this same reasoning gives:

\proclaim Theorem 1.5.  Let $M$ be a convex subset of $X$ and $f\in X\\
{M}$. Assume $m^*\in P_M(f)$. Set
$$\gamma (f) := \inf\left\{ {{\tau_+(f-m^*, m^*-m)}\over  {\|m^*-m\|}}: m\in M, m\ne m^*\right\}.$$
Then $\gamma (f) \ge 0$ and for all $m\in M$, we have
$$\|f-m \| - \|f-m^*\| \ge \gamma (f)\, \|m-m^*\|.$$
Furthermore, if $\gamma > \gamma (f)$ then there exists an $\tilm\in
M$ for which
$$ \|f-\tilm \| - \|f-m^*\| < \gamma\, \|\tilm-m^*\|.$$

In much of this paper, we consider approximation from linear
subspaces and not from convex subsets. However, most of the results
obtained can be easily generalized to convex subsets, as above. We
leave their exact statements to the interested reader.

Note that if $M$ is a finite-dimensional subspace (or closed convex
subset thereof) then a compactness argument implies that the {\sl
inf} in the definition of $\gamma(f)$ in both Theorems 1.4 and 1.5
is in fact a {\sl min}. In the former case, we rewrite an essential
consequence of Theorem 1.4 as:

\proclaim Proposition 1.6.  Let $M$ be a finite-dimensional subspace of $X$ and $f\in X\\
{M}$. Assume $m^*\in P_M(f)$. Then $m^*$ is a strongly unique best
approximant to $f$ from $M$, i.e., $\gamma(f)$, as defined in
Theorem 1.4, is strictly positive, if and only if there does not
exist an $m\in M$, $m\ne 0$, for which $\tau_+(f-m^*, m)=0$.

For many normed linear spaces, we have that
$$\tau(f, g) =\lim_{t\to 0}{{\|f+tg\|-\|f\|}\over{t}}$$
exists for all $f, g\in X$, $f\ne 0$. (Note that this is the
two-sided limit that need not in general exist.) If this is the case
then, for any linear subspace $M$, we have that $m^*\in P_M(f)$ if
and only if
$$\tau(f-m^*, m) =0$$
for all $m\in M$. Thus, we never have classical strong uniqueness in
such spaces. If the functional $\tau(f,g)$ exists for all $f, g\in
X$, $f\ne 0$, then the space $X$ is said to be {\sl smooth}. Strong
uniqueness cannot hold in smooth spaces as strong uniqueness implies
that the unit ball has corners. The space $X$ being smooth is
equivalent to the existence, for each $f\in X$, $f\ne 0$, of a
unique continuous linear functional $h\in X^*$ of norm $1$
satisfying $h(f) = \|f\|$. This is true, for example, in any Hilbert
space and in $L^p$ for every $p\in (1,\infty)$. In all these cases,
we must look for other non-classical type strong uniqueness
formulae.

The approach that was taken in this section to characterizing best
approximants and classical strong uniqueness may be found in Papini
[1978] and Wulbert [1971].

\goodbreak

\medskip
\subsect{{2.} Classical Strong Uniqueness in the Uniform Norm}

\medskip\noindent
In this section, we assume that $B$ is a compact Hausdorff space and
consider $C(B)$, the normed linear space of continuous real-valued
functions defined on $B$, with norm
$$\|f\| = \max_{x\in B} |f(x)|.$$
For each $f\in C(B)$, set $A_f = \{ x : |f(x)|=\|f\|\}$ and
$$\sgn f(x) =\cases{ 1, & $f(x)>0$,\cr 0, & $f(x)=0$,\cr -1, &
$f(x)<0$.\cr}$$ The formula for the one-sided Gateaux derivative
$\tau_+(f,g)$ in $C(B)$ is the following.

\proclaim Theorem 2.1. For $f,g\in C(B)$, $f\ne 0$, we have
$$\tau_+(f,g) = \max_{x\in A_{f}}\,[\sgn f(x)]\, g(x).$$

\pf We first prove that the right-hand-side is a lower bound for
$\tau_+(f,g)$. To this end, let $x\in A_f$. Thus $|f(x)|=\|f\|$ and
$|f(x) + tg(x)| \le \|f+tg\|$. Therefore
$$\eqalign{\tau_+(f,g) = & \lim_{t\to 0^+} {{\|f+tg\|-\|f\|}\over
t}\cr \ge & \lim_{t\to 0^+} {{|f(x)+tg(x)|-|f(x)|}\over t}\cr = &
\lim_{t\to 0^+} {{[\sgn f(x)]\,(f(x)+tg(x))-[\sgn f(x)]\,f(x)}\over
t}\cr = & \lim_{t\to 0^+} {{t[\sgn f(x)]\,g(x)}\over t}\cr = &
\,[\sgn f(x)]\,g(x).\cr}$$ This implies that
$$\tau_+(f,g) \ge \max_{x\in A_{f}}\,[\sgn f(x)]\, g(x).$$

We prove the converse direction as follows. For each $t>0$, let
$x_t$ satisfy
$$\|f+tg\| = |f(x_t)+ tg(x_t)|.$$
As $B$ is compact, there exists an $x^*\in B$ that is a limit point
of the $x_t$ as $t\to 0$. Let $t_n$ denote a sequence, decreasing to
zero, along which $x_{t_n} := x_n$ converges to $x^*$. We first
claim that we must have $x^*\in A_f$. If not, then
$$\|f+t_n g\| = |f(x_n)+ t_n g(x_n)| \rightarrow |f(x^*)|<\|f\|,$$
contradicting the continuity of the norm. As $x^*\in A_f$, it
therefore follows that for $n$ sufficiently large,
$$\sgn(f(x_n)+t_n g(x_n)) = \sgn f(x^*).$$
Thus
$$\eqalign{{{\|f+t_n g\|-\|f\|}\over {t_n}} = & {{ [\sgn f(x^*)] [f(x_n)+t_n
g(x_n)] - [\sgn f(x^*)]f(x^*)}\over {t_n}}\cr = & [\sgn f(x^*)]
\left[ {{f(x_n)-f(x^*)}\over {t_n}}\right] + [\sgn f(x^*)]
g(x_n).\cr}$$ Note that
$$ [\sgn f(x^*)] \left[ {{f(x_n)-f(x^*)}\over {t_n}}\right] \le
0$$ since $t_n>0$ and  $[\sgn f(x^*)] f(x_n) \le \|f\| =|f(x^*)| =
[\sgn f(x^*)] f(x^*)$. Take limits on both sides of the above
equality. The left-hand-side has the finite limit $\tau_+(f,g)$, the
first term on the right-hand-side is nonpositive while
$$\lim_{n\to\infty} [\sgn f(x^*)]\, g(x_n) = [\sgn f(x^*)]\, g(x^*)
\le \max_{x\in A_{f}}\,[\sgn f(x)]\, g(x).$$ Thus
$$\tau_+(f,g) \le \max_{x\in A_{f}}\,[\sgn f(x)]\, g(x). \meop$$

As a consequence of Theorem 2.1, Theorem 1.2 and Theorem 1.4, we
have:

\proclaim Theorem 2.2. Let $M$ be a linear subspace of $C(B)$. Then
$m^*\in P_M(f)$ if and only if
$$\tau_+(f-m^*, m) = \max_{x\in A_{f-m^*}}\,[\sgn (f-m^*)(x)]\, m(x)\ge 0$$
for all $m\in M$. Furthermore, $m^*$ is a strongly unique best
approximant to $f$ from $M$ if and only if
$$\gamma(f) = \inf_{{{m\in M}\atop {\|m\|=1}}} \max_{x\in A_{f-m^*}}\,[\sgn (f-m^*)(x)]\, m(x) >0.$$

This characterization of a best approximant is called the Kolmogorov
criterion. This characterization of the optimal strong uniqueness
constant $\gamma(f)$ first appeared in Bartelt, McLaughlin [1973].
If $M$ is finite-dimensional, then the above {\sl inf} can be
replaced by a {\sl min}. This is not true in infinite dimensions. An
example illustrating that fact may be found in Bartelt, McLaughlin
[1973]; see Example 4.

In the finite-dimensional setting, the above characterization
results can be further refined. We firstly state the following
well-known result.

\proclaim Theorem 2.3. Let $M$ be an $n$-dimensional subspace of
$C(B)$. Given $f\in C(B)$, we have  $m^*\in P_M(f)$ if and only if
there exist $k$ distinct points $x_1\nek x_k$ in $A_{f-m^*}$ , and
strictly positive values $\lam_1\nek \lam _k$, $1\le k\le n+1$, such
that
$$\sum_{i=1}^k \lam_i [\sgn (f-m^*)(x_i)] m(x_i)= 0$$
for all $m\in M$.

(A variation on a proof of Theorem 2.3 is given by the proof of
Theorem 10.3.)

We also have the following result that is proved in Brosowski
[1983], and was also later proved in Smarzewski [1990], generalizing
results from Bartelt [1974] and N\"urnberger [1980].

\proclaim Theorem 2.4. Let $M$ be an $n$-dimensional subspace of
$C(B)$. Given $f\in C(B)\\ M$, we have that $m^*\in M$ is the
strongly unique best approximant to $f$ from $M$ if and only if
there exist $k$ distinct points $x_1\nek x_k$ in $A_{f-m^*}$, and
strictly positive values $\lam_1\nek \lam _k$, with $n+1\le k\le
2n$, such that
$$\sum_{i=1}^{k} \lam_i [\sgn (f-m^*)(x_i)] m(x_i)= 0\eqno{(2.1)}$$
for all $m\in M$, and
$$\dim M|_{\{x_1\nek x_k\} } = n.$$

\pf We first prove the ``easy'' direction in this theorem.

\noindent ($\Leftarrow$). Assume there exist points
$\{x_i\}_{i=1}^k$ satisfying the conditions of the theorem. As the
$x_1\nek x_k$ are in $A_{f-m^*}$, we have
$$\gamma(f)  \ge \gamma  = \min_{{{m\in M}\atop {\|m\|=1}}} \max_{i=1\nek k}\sgn (f-m^*)(x_i) m(x_i).$$
We claim that $\gamma>0$. If this holds then $\gamma(f)>0$ and $m^*$
is the strongly unique best approximant to $f$ from $M$. To see that
$\gamma>0$, recall that since $M$ is a finite-dimensional subspace
then from a compactness argument the above minimum is always
attained. Thus if $\gamma \le 0$, there exists an $\tilm\in M$, $\|
\tilm\|=1$, for which
$$\sgn (f-m^*)(x_i)\tilm (x_i) \le 0$$
for $i=1\nek k$. Moreover as $\tilm \in M$, $\tilm \ne 0$, then from
the fact that $\dim M|_{\{x_1\nek x_k\} } = n$, it follows that
$\tilm(x_i)\ne 0$ at at least one of the $x_i$. As $\lam_i>0$,
$i=1\nek k$, this then implies
$$ \sum_{i=1}^{k} \lam_i \sgn (f-m^*)(x_i) \tilm(x_i) <0$$
which contradicts (2.1). Thus $\gamma>0$.

\smallskip\noindent
($\Rightarrow$). Assume $m^*\in M$ is a strongly unique best
approximant to $f$ from $M$. Without loss of generality, we will
assume $m^*=0$.

We start with the case where $\dim M =1$. Let $M=\span \{m_1\}$.
Since
$$\gamma(f) =\max_{x\in A_f} [\sgn f(x)] (\pm m_1(x)) >0,$$
there must exist points $x_1, x_2\in A_f$ for which we have that
both $[\sgn f(x_1)] m_1(x_1)>0$ and $[\sgn f(x_2)] m_1(x_2) < 0$.
Thus, there exist $\lam_1, \lam_2 >0$ satisfying
$$\lam_1 [\sgn f(x_1)] m(x_1) + \lam_2 [\sgn f(x_2)] m(x_2) =0$$
for all $m\in M$, and $\dim M|_{\{x_1, x_2\}} =1$. This is the
desired result.

We now consider the general case where $\dim M=n$. Since
$$\max_{x\in A_f} [\sgn f(x)] m(x)\ge c >0$$
for all $m\in M$, $\|m\|=1$, it follows that ${\bf 0}\in \RR^n$ is
in the strict interior of the convex hull of
$$E = \{([\sgn f(x)] m_1(x) \nek [\sgn f(x)] m_n(x)): \,x\in A_f\} \subset \RR^n$$
where $\{m_1\nek m_n\}$ is any basis for $M$. From a generalization
of Carath\'eodory's Theorem, essentially due to Steinitz, see
Danzer, Gr\"unbaum, Klee [1963], the vector ${\bf 0}\in \RR^n$ is in
the strict interior of the convex hull of some set of at most $2n$
points of $E$. That is, there exist $x_1\nek x_k\in A_f$, $k\le 2n$,
such that ${\bf 0}\in \RR^n$ is in the strict interior of the convex
hull of
$$E^* = \{([\sgn f(x_i)] m_1(x_i) \nek [\sgn f(x_i)] m_n(x_i)): \, i=1\nek k\}.$$

From the fact that ${\bf 0}\in \RR^n$ is in the strict interior of
the convex hull of $E^*$, it easily follows that there exist
$\lam_1\nek \lam _k > 0$, $\sum_{i=1}^k \lam_i=1$, for which
$$0=\sum_{i=1}^k \lam_i [\sgn f(x_i)] m_j(x_i), \qquad j=1\nek n,$$
implying
$$0=\sum_{i=1}^k \lam_i [\sgn f(x_i)] m(x_i),$$
for all $m\in M$, and also that the vectors
$$([\sgn f(x_i)] m_1(x_i) \nek [\sgn f(x_i)] m_n(x_i)), \qquad i=1\nek
k,$$ span $\RR^n$. Since $\sgn f(x_i)\ne 0$ for each $i$, this is
equivalent to the fact that the vectors
$$(m_1(x_i) \nek  m_n(x_i)), \qquad i=1\nek k,$$
span $\RR^n$, i.e.,
$$\dim M|_{\{x_1\nek x_k\} } = n.$$
From this fact and
$$0=\sum_{i=1}^k \lam_i [\sgn f(x_i)] m(x_i),$$
for all $m\in M$, it follows that $k\ge n+1$. \eop

Is the bound on $k$, namely $n+1\le k\le 2n$, the correct bound? The
answer is yes. The value $k=n+1$ is minimal. If $k\le n$ there is
always a nontrivial $m\in M$ that vanishes at $x_1\nek x_{k-1}$ and
then $\min\{[\sgn (f-m^*)(x_k)] (\pm m(x_k))\} \le 0$ implying
$\gamma(f)=0$. In the next section, we see that in the Haar space
setting, we can and do have $k=n+1$. The case $k=2n$ is necessary
if, for example, $M$ is spanned by $n$ functions with disjoint
support. To see this, assume $m_1\nek m_n$ is a basis for $M$ where
these basis functions have disjoint support. Assume $f\in C(B) \\ M$
and consider (2.1). In order that $\dim M|_{\{x_1\nek x_k\} } = n$,
it is necessary that among the $x_1\nek x_k$ there is at least one
point in the support of each $m_j$. But we must have at least two
points in the support of each $m_j$ among the $x_1\nek x_k$ if (2.1)
is to hold. Thus $k=2n$ is necessary for strong uniqueness in this
case and, by the above, there do exist $f\in C(B)\\ M$ with strongly
unique best approximants from $M$.

It might be conjectured that uniqueness and strong uniqueness are
equivalent properties in $C(B)$. This, however, is not true, as can
be seen from this example taken from Cheney [1966, p.~82].

\medskip \noindent{\bf Example.}  Let $M=\span\{x\}$ in $C[-1,1]$, and set
$f(x) = 1-x^2$. As is easily verified, the unique best approximant
to $f$ from $M$ is $m^*(x)=0$. On the other hand, $m^*$ is not a
strongly unique best approximant to $f$ from $M$ since $A_f =\{0\}$,
and $m(0)=0$, i.e., $\tau_+(1-x^2, x) =0$.

\medskip
Nevertheless, while uniqueness and strong uniqueness are not
equivalent properties in $C(B)$, the set of functions with a
strongly unique best approximant is dense in the set of functions
with a unique best approximant when approximating from a
finite-dimensional subspace. This next result and the subsequent
Corollary 2.7 are from N\"urnberger, Singer [1982]. The proof given
here is from Smarzewski [1988].

\proclaim Theorem 2.5. Let $M$ be a finite-dimensional subspace of
$C(B)$. Then the set of functions with a strongly unique best
approximant is dense in the set of functions with a unique best
approximant.

\medskip
Note that we have an algebraic characterization for a best
approximant and an algebraic characterization for a strongly unique
best approximant from any finite-dimensional subspace (Theorems 2.2,
2.3 and 2.4). However there is no known algebraic characterization
for when we have a unique best approximant. We get around this by
proving the next proposition. We start with some notation.

Let us assume, for ease of presentation, that the zero function is a
best approximant to $f$ from $M$. Thus from Theorem 2.2, we have
$$\max_{x\in A_f} [\sgn f(x)] m(x)\ge 0$$
for all $m\in M$. Let $(k_n)$ be any strictly decreasing sequence of
positive numbers that converges to zero, and assume $k_n < 1/2$ for
all $n$. Set
$$B_n := \{ x : |f(x)| > ( 1-k_n) \|f\|\}, \qquad
n=2,3,\ldots\ .$$ Then we have:

\proclaim Proposition 2.6. Let $M$ be a finite-dimensional subspace
of $C(B)$. If the zero function is the unique best approximant to
$f$ from $M$ then for each $n$, there exists an $a_n>0$ such that
$$\max_{x\in B_n} [\sgn f(x)] m(x) \ge a_n \|m\|$$
for all $m\in M$, where $B_n$ is as above.

\smallskip\noindent
{\bf Remark.} The above is a necessary, but not a sufficient
condition, for the uniqueness of the best approximant. As an
example, consider $f(x) = 1 - |x|$ on $[-1,1]$, and $M=\span \{m_1
\}$ where $m_1(x)=x$ on $[-1,1]$. Then, as is readily verified, $B_n
= [-k_n, k_n]$ and
$$\max_{x\in B_n} [\sgn f(x)] m(x) = k_n \|m\|$$
for all $m\in M$. Moreover, while the zero function is a best
approximant to $f$ from $M$, it is not the unique best approximant
as we also have $\pm m_1\in P_M(f)$.

\medskip
\pf Assume to the contrary that the desired inequality does not hold
for some $n\ge 2$. This inequality trivially holds for $m=0$. As
such, we can consider it over the boundary of the unit ball of $M$,
and since $M$ is a finite-dimensional subspace, we can use
compactness to affirm that there exists an $\tilm\in M$,
$\|\tilm\|=1$, such that
$$\max_{x\in B_n} [\sgn f(x)] \tilm (x)\le 0.$$

Choose any $\alp>0$ satisfying $\alp \le k_n \|f\|$. Thus
$$|\alp \tilm (x) | < |f(x)|$$
for all $x\in B_n$, and therefore, we also have thereon
$$|f(x) +\alp \tilm (x)| \le |f(x)| +\alp [\sgn f(x)] \tilm (x) \le
|f(x)| \le \|f\|.$$ On the other hand, for $x\in B\\ B_n$, we have
$$ |f(x) +\alp \tilm (x)| \le |f(x)| +\alp |\tilm (x)| \le ( 1-k_n) \|f\| +
k_n \|f\| = \|f\|.$$ Thus $-\alp \tilm \in P_M(f)$, contradicting
the fact that the zero function is the unique best approximant to
$f$ from $M$. \eop

One immediate consequence of the above proposition is the following.

\proclaim Corollary 2.7. If $M$ is  a subspace of $\ell^m_\infty$,
then every unique best approximant to $f\in \ell^m_\infty$ from $M$
is also a strongly unique best approximant to $f$ from $M$.

\pf For $n$ sufficiently large, we always have, in this case, $B_n =
A_f$. Apply Proposition 2.6 and the characterization of strong
uniqueness found in Theorem 2.2. \eop

\medskip\noindent
{\bf Proof of Theorem 2.5:} If $f\in M$, there is nothing to prove.
As such, we assume that $f\in C(B) \\ M$.  Without loss of
generality, we assume that the zero function is the unique best
approximant to $f$ from $M$. Let $B_n$ and $a_n$ be as in
Proposition 2.6.

By the Tietze-Urysohn Theorem, see e.g., Kuratowski [1966], there
exists a $f_n\in C(B)$ such that
$$f_n(x) = \cases { \|f\|\sgn f(x), & $x\in B_{n+2}$\cr f(x), &
$x\in B_n \cap \overline{(B \\ B_{n+2})}$\cr}$$ and
$$( 1 - k_n) \|f\| \le f_n(x) \le \|f\|$$
for $x\in B_n$. If $B_n=B_{n+2}$, then we simply set $f_n(x) =
\|f\|\sgn f(x)$ thereon. We extend $f_n$ to all of $C(B)$ by setting
$$f_n(x) = f(x)$$
for $x\in B\\ B_n$. Note that we do not lose continuity in the case
where $B_n = B_{n+2}$ since in that case
$$\{ x: ( 1-k_{n+2}) \|f\| \ge |f(x)|\ge ( 1- k_n)
\|f\|\} =\emptyset.$$ From this construction, we have
$$|f(x) - f_n(x) | \le k_n \|f\|$$
for all $x\in B$. That is, we have
$$\lim_{n\to\infty} \|f-f_n\|=0.$$

Now, since $A_{f_n} \supseteq B_{n+2} \supseteq A_f$, it follows
from Theorem 2.2 that $0\in P_M(f_n)$ for all $n$. Furthermore from
Proposition 2.6, we have
$$\min_{{{m\in M}\atop {\|m\|=1}}} \max_{x\in A_{f_n}} [\sgn f_n(x)]
m(x) \ge  \min_{{{m\in M}\atop {\|m\|=1}}} \max_{x\in B_{n+2}} [\sgn
f(x)] m(x) \ge a_{n+2} >0.$$ Thus, again applying Theorem 2.2, we
see that the zero function is the strongly unique best approximant
to $f_n$ from $M$. \eop

The above begs the question of when the set of functions with a
unique best approximant from a given finite-dimensional subspace $M$
is dense in $C(B)$. A similar question was considered by Garkavi
[1964], [1965] who proved that the $n$-dimensional subspace $M$ of
$C(B)$ has the property that every function in $C(B)$ has a unique
best approximant from $M$ except for a set of first category in
$C(B)$ if and only if on each open subset $D$ of $B$ there
identically vanish at most $(n- |D|)_+ =\max\{n-|D|, 0\}$ linearly
independent functions of $M$, where $|D|$ is the number of points in
$D$. Thus, for example, the set of functions with a strongly unique
best approximant from $M$ is dense in $C[a,b]$ if no $m\in M$, $m\ne
0$, vanishes identically on an open interval.

\medskip Bounds for $\gamma(f)$ are difficult to obtain. The following
result, from Grothmann [1989], gives an upper bound for $\gamma(f)$.
In general, using it to find explicit bounds is very difficult.

Before stating the result, we recall the definition of a {\sl
projection constant} of a subspace $L$ in a normed linear space $X$.
It is given by
$$\lam (L;X) := \inf \{ \|P\|:\, P: X\to L {\rm \ is\ a\
projection}\}.$$

\proclaim Proposition 2.8. Let $M$ be a linear subspace of $C(B)$.
Let $f\in C(B)$ and assume $m^*\in M$ is a strongly unique best
approximant to $f$. Then
$$\gamma(f)\le {{\lam(M|_{A_{f-m^*}}; C(A_{f-m^*}))}\over {\lam
(M; C(B))}}.$$

\pf Since $m^*$ is a strongly unique best approximant to $f$ from
$M$, we have that
$$\gamma(f) = \min_{{{m\in M}\atop {\|m\|=1}}} \max_{x\in A_{f-m^*}}[\sgn (f-m^*)(x)] m(x)  > 0.$$
This also implies that
$$\|m\|_{A_{f-m^*}} = \max_{x\in A_{f-m^*}} |m(x)|$$
is a norm on $M$ since no $m\in M$, $m\ne 0$, can vanish identically
on $A_{f-m^*}$. From the above formula for $\gamma(f)$, we have
$$\eqalign{\gamma(f) \le &  \min_{{{m\in M}\atop {\|m\|=1}}} \max_{x\in
A_{f-m^*}} |m(x)|\cr = & \min_{m\in M} {{\,\,\|m\|_{A_{f-m^*}}}\over
{\|m\|\,\,}}\cr = & \left(\max_{m\in M} {{\|m\|\,\,}\over
{\,\,\|m\|_{A_{f-m^*}}}}\right)^{-1}.\cr}$$ Let $P:C(A_{f-m^*})\to
M|_{A_{f-m^*}}$ be a projection. Define $\sig : M|_{A_{f-m^*}}\to M$
and $\phi: C(B) \to C(A_{f-m^*})$ by
$$\sig\left( m|_{A_{f-m^*}}\right) = m$$
for $m\in M$ and
$$\phi(f) = f|_{A_{f-m^*}}$$
for $f\in C(B)$. The linear operator $\sig$ is well-defined since no
$m\in M$, $m\ne 0$, can vanish identically on $A_{f-m^*}$. Set
$$Q(f) := \sig(P(\phi(f))).$$
$Q$ is a projection from $C(B)$ onto $M$. Thus
$$\lam (M; C(B)) \le \|Q\| \le \|\sig\|\|P\| \|\phi\|.$$
Now $\|\phi\| = 1$ and, by the above, $\|\sig\| \le 1/\gamma(f)$.
Thus
$$\lam (M; C(B)) \le {1\over {\gamma(f)}} \|P\|$$
for all projections $P :C(A_{f-m^*})\to M|_{A_{f-m^*}}$. Therefore
$$\gamma(f)\le {{\lam(M|_{A_{f-m^*}}; C(A_{f-m^*}))}\over {\lam
(M; C(B))}}. \meop$$

In Section 4, where we assume $M$ is a Haar space, we discuss other
somewhat more computable bounds for $\gamma(f)$.

\goodbreak

\medskip
\subsect{{3.} Local Lipschitz Continuity and Classical Strong
Uniqueness}

\medskip\noindent
We start once again with an arbitrary normed linear space $X$.
Assume $M$ is a subset of $X$, and to $f\in X$ there exists a unique
best approximant from $M$. If, for a given $f\in X$, we have an
inequality of the form
$$\| P_M (f) - P_M (g)\| \le \sigma \|f-g\|$$
valid for all $g\in X$ and any element of $P_M(g)$, then we say that
the best approximation operator from $M$ is {\sl locally Lipschitz
continuous} at $f$, and call $\sigma$ a {\sl local Lipschitz
constant}.

One of the ``uses'' of classical strong uniqueness is that it
implies local Lipschitz continuity. Before proving this result, we
recall that we always have, for every $f, g\in X$,
$$|\, \|f-P_M(f)\| - \|g-P_M(g)\| \,| \le \|f-g\|.$$
To verify this, note that, assuming $ \|f-P_M(f)\| \ge
\|g-P_M(g)\|$, then
$$\|f-P_M(f)\| \le \|f -P_M(g)\| \le \|f-g\| + \|g-P_M(g)\|,$$
from which the result follows.

\proclaim Theorem 3.1. Assume $M$ is a subset of a normed linear
space $X$, $f\in X$, and for some $\gamma>0$, we have
$$ \|f-m\| - \|f-P_M(f)\| \ge \gamma \|m-P_M(f)\|$$
for all $m\in M$. Then for each $g\in X$ and any element of $P_M(g)$
$$\| P_M(f) - P_M(g)\| \le {2\over {\gamma}} \|f-g\|.$$

\pf By assumption, and an application of the previous inequality,
$$ \gamma \|P_M (f)- P_M (g)\| \le   \|f-P_M(g)\| - \|f-P_M(f)\|$$
$$\le \|f-g\| + \|g -P_M(g)\| - \|f-P_M(f)\| \le 2 \|f-g\|.$$ Thus
$$\| P_M(f) - P_M (g)\| \le {2\over {\gamma}} \|f-g\|. \meop$$

In the specific case where $M$ is a finite-dimensional Haar space,
see Section 4, local Lipschitz continuity of the best approximation
operator was proved in Freud [1958] by a different method of proof.
For polynomial approximation on an interval, this result can be
found in Kirchberger [1902, p.~18--21] (see also Borel [1905,
p.~89--92]).

We have proven that strong uniqueness at $f$ implies local Lipschitz
continuity at $f$. The converse direction, i.e., local Lipschitz
continuity implying strong uniqueness, does not necessarily hold. It
certainly does not hold in an inner product space where we always
have
$$\| P_M(f) - P_M(g)\| \le  \|f-g\|$$
for all $f, g$ if $M$ is a subspace. That is, we have local
Lipschitz continuity and do not have strong uniqueness, see the
discussion at the end of Section 1. However in $C(B)$ the converse
does hold, i.e., local Lipschitz continuity at $f\in C(B)$ does
imply strong uniqueness at this same $f$, assuming $M$ is a
finite-dimensional subspace. This next theorem is to be found in
Bartelt, Schmidt [1984].

\proclaim Theorem 3.2. Let $B$ be a compact Hausdorff space and $M$
a finite-dimensional subspace of $C(B)$. For given $f\in C(B)$, the
following are equivalent:
\item{(I)} There exists a $\gamma >0$ such that
$$ \|f-m\| - \|f-P_M (f)\| \ge \gamma \|m-P_M (f)\|$$
for all $m\in M$.
\item{(II)} There exists a $\sig>0$ such that
$$\| P_M(f) - P_M(g)\| \le \sig \|f-g\|$$
for all $g\in C(B)$.

\pf From Theorem 3.1, we have that (I) implies (II). It remains to
prove the converse direction. We prove the result by contradiction.
That is, we assume that (I) does not hold for a given $f$ and will
prove that (II) does not hold for this same $f$. Note that if  $f\in
M$ then (I) and (II) hold, while if the best approximant to $f$ from
$M$ is not unique, then neither (I) nor (II) hold. As such, we
assume $f\notin M$ and the best approximant to $f$ from $M$ is
unique.

Without loss of generality, we can and will assume that the zero
function is the unique best approximant to $f$, and $\|f\|=1$. Now,
since (I) does not hold, there exists, for each $\eps>0$, an
$m_\eps\in M$ satisfying
$$\|f-m_\eps\| < \|f\| + \eps \|m_\eps\| = 1 + \eps \|m_\eps\|.$$
Note that this implies that $m_\eps \ne 0$. We divide most of the
proof into a series of four lemmas.

\proclaim Lemma 3.3. There exist sequences $(\del_n)$ and $(\alp_n)$
of positive numbers tending to zero and an $\tilm\in M$, $\|\tilm
\|=1$, such that
\item{(i)} $\| f- \alp_n \tilm \| \le 1 + \del_n \alp_n$,
\item{(ii)} $\tilm(x)\sgn f(x) \ge 0$ for all $x\in A_f$.

\pf Let $m_\eps$ be as above. We first claim that for $\eps \le
1/2$, we necessarily have $\| m_\eps \| \le 4$. To see this, assume
$\|m_\eps \| =C$. Then for $0<\eps \le 1/2$, we have
$$ {C \over 2} \ge \eps \|m_\eps\| \ge \|f-m_\eps\| -1 \ge
\|m_\eps\| - \|f\| -1 = C-2,$$ whence $C\le 4$.

As $\| m_\eps\| \le 4$, we have
$$\lim_{\eps \to 0^+} \eps \|m_\eps \| =0.$$
Thus
$$ \lim_{\eps \to 0^+} \|f-m_\eps\| = \|f\| = \min_{m\in M}
\|f-m\|.$$ From a compactness argument, it therefore follows that
$$\lim_{\eps \to 0^+} \|m_\eps\| = 0,$$
i.e., $m_\eps$ tends to the  zero function. Set
$$m_\eps := \tilm_\eps \| m_\eps \|$$
and let $\alp_\eps =\|m_\eps\|$. Thus $\|\tilm_\eps \| = 1$ for all
$\eps$, and $\alp_\eps$ tends to $0$ as $\eps$ tends to $0$. As the
$\tilm_\eps$ are all functions of norm $1$ in a finite-dimensional
subspace, on a subsequence $\eps_n \to 0^+$, we have
$$\lim_{n\to \infty} \tilm_{\eps_n} = \tilm$$
where $\tilm \in M$ and $\|\tilm \| =1$. For convenience, let
$\tilm_n :=\tilm_{\eps_n}$, $m_n :=m_{\eps_n}$ and $\alp_n
:=\alp_{\eps_n}$.

We claim that $\tilm(x) \sgn f(x) \ge 0$ for all $x\in A_f$. Assume
not. There then exists an $x^*\in A_f$ for which
$$\tilm (x^*) \sgn f(x^*) = -c <0.$$
Thus for $n$ sufficiently large,
$$\tilm_{n} (x^*) \sgn f(x^*) < -{c\over 2},$$
and therefore
$$\eqalign{ 1 + \eps_n \alp_{n} &\;=\; \|f\| + \eps_n \| m_{n}\|  >
\|f- m_{n} \|\cr &\;\ge\; | f(x^*) - m_{n} (x^*)| \;=\; |f(x^*)| +
\alp_{n} |\tilm_{n}(x^*)| \;>\; 1 + \alp_{n} {c\over 2}.\cr}$$ But
this cannot possibly hold for $n$ sufficiently large as $\eps_n \to
0$. Thus $\tilm(x) \sgn f(x) \ge 0$ for all $x\in A_f$.

Let $\beta_n = \|\tilm_{n} -\tilm \|$. Therefore $\beta_n$ tends to
zero as $n$ tends to infinity. Now
$$\eqalign{\left[ \|f-\alp_{n} \tilm \| -\|f\|\right]
- \left[ \|f-\alp_{n} \tilm_{n} \| -\|f\|\right] &\;=\; \|f-\alp_{n}
\tilm \| - \|f-\alp_{n} \tilm_{n} \|\cr &\;\le\; \| \alp_{n} \tilm -
\alp_{n} \tilm_{n} \| \;=\; \alp_{n} \beta_n.\cr}$$ Thus
$$ 0 \le \|f- \alp_{n} \tilm \| -\|f\| \le
\|f-\alp_{n} \tilm_{n} \| -\|f\| + \alp_{n} \beta_n \le \alp_{n}
\eps_n + \alp_{n} \beta_n = \alp_{n}(\eps_n + \beta_n).$$ Set
$\del_n := \eps_n + \beta_n$. Then
$$\lim_{n\to\infty}\del_n=0$$
and this proves the lemma. \eop

Recall that we assumed that $0\in P_M(f)$. This implies, by Theorem
2.3, the existence of $k$ distinct points $x_1\nek x_k \in A_f$,
$1\le k\le n+1$, and strictly positive values $\lam_1\nek \lam _k$
such that
$$\sum_{i=1}^k \lam_i [\sgn f(x_i)] m(x_i)= 0$$
for all $m\in M$.

\proclaim Lemma 3.4. $\tilm(x_i)=0$, for $i=1\nek k$.

\pf From Lemma 3.3 (ii), we have
$$0 \le \tilm(x_i) \sgn f(x_i)$$
for all $i=1\nek k$. Thus if we have strict inequality for any $i$,
then we contradict the fact that
$$ \sum_{i=1}^k \lam_i [\sgn f(x_i)] \tilm(x_i)= 0.\ \meop$$

For each given $n$, let
$$G_n= \{x : |\tilm(x)|< \del_n/2\}.$$
Note that $G_n$ is an open neighborhood of the $\{x_1\nek x_k\}$. We
define a function $\phi_n \in C(B)$ as follows. Firstly set
$\phi_n(x_i) =\alp_n\del_n$, $i=1\nek k$, and $\phi_n(x) =0$ for all $x\in B\\
G_n$. Note that on $\{x_1\nek x_k\} \cup B\\ G_n$, we have
$$0 \le \phi_n(x) \le | \alp_n \del_n -\alp_n |\tilm(x)|\,|.$$
This follows from the fact that $\tilm(x_i)=0$, $i=1\nek k$. Extend
$\phi_n$ continuously to all of $B$ so that it continues to satisfy
$$0 \le \phi_n(x) \le | \alp_n \del_n -\alp_n |\tilm(x)|\,|$$
for all $x\in B$. Now set $g_n(x) := f(x) \left[ 1+
\phi_n(x)\right]$.

\proclaim Lemma 3.5. We have
$$\|f-g_n\| = \alp_n \del_n.$$

\pf From the definition of $g_n$, $f(x) - g_n(x) = - f(x)
\phi_n(x)$. Since $\|f\|=1$, we therefore have
$$\|f-g_n\| \le \|\phi_n\|.$$
On $G_n$, where $\phi_n$ need not vanish, we have
$$|\tilm(x)| < {{\del_n}\over 2}.$$
Thus for $x\in G_n$
$$0 \le \phi_n(x) \le | \alp_n \del_n -\alp_n |\tilm(x)|\,| =  \alp_n \del_n -\alp_n
|\tilm(x)| \le \alp_n \del_n.$$ Therefore $\|f-g_n\| \le \alp_n
\del_n$. Equality holds since
$$|f(x_i) - g_n(x_i) | =  |f(x_i) \phi_n(x_i)| = | \phi_n(x_i) | =
\alp_n \del_n. \meop$$

\proclaim Lemma 3.6. We have
$$\min_{m\in M} \| g_n -m\| = \|g_n - \alp_n \tilm\| = 1 +
\alp_n\del_n.$$

\pf For each $i=1\nek k$, since $\tilm(x_i)=0$, we have
$$ |g_n(x_i) -\alp_n \tilm(x_i)| = |g_n(x_i)| = |f(x_i)| \left| 1+
\phi_n(x_i)\right| =  1+ \phi_n(x_i) = 1 + \alp_n\del_n.$$

Now for $x\in B\\ G_n$, since $\phi_n(x) =0$, we have
$$ |g_n(x) -\alp_n \tilm(x)| = |f(x) -\alp_n \tilm(x)| \le
\|f-\alp_n \tilm\| \le 1 + \alp_n \del_n.$$ The latter inequality is
from Lemma 3.3 (i). For $x\in G_n$, since $|f(x)|\le 1$ and
$|\tilm(x)| < {\del_n}/2$, we have
$$\eqalign{ |g_n(x) -\alp_n \tilm(x)|
&\;=\;| f(x) + f(x)\phi_n(x)  -\alp_n \tilm(x)|\cr &\; \le\; 1 + |
\alp_n \del_n -\alp_n |\tilm(x)|\,| + \alp_n |\tilm(x)|\cr &\; \le\;
1 +  \alp_n \del_n -\alp_n |\tilm(x)| + \alp_n |\tilm(x)| \;=\; 1 +
\alp_n \del_n.\cr}$$ This proves that
$$ \|g_n - \alp_n \tilm\| = 1 +
\alp_n\del_n.$$

As
$$ g_n(x_i) -\alp_n \tilm(x_i)= g_n(x_i) = f(x_i) \left[1 + \alp_n\del_n \right]$$
for $i=1\nek k$, it follows from Theorem 2.3 characterizing best
approximants that $\alp_n\tilm$ is a best approximant to $g_n$ from
$M$. \eop

\medskip\noindent
{\bf Proof of Theorem 3.2 (cont'd):} If (II) holds, then there
exists a $\sig >0$ such that
$$ \| P_M f - P_M g\| \le \sig \|f-g\|$$
for all $g\in C(B)$. Recall that we assumed, without loss of
generality, that the zero function is the best approximant to $f$
from $M$ and $\|f\|=1$. Taking $g=g_n$ and applying the above
lemmas, we obtain
$$ \alp_n = \|\alp_n \tilm\| \le \sig \|f- g_n\| = \sig \alp_n \del_n.$$
Thus $1 \le \sig \del_n$ for all $n$. But $\lim_{n\to\infty} \del_n
=0$. This is a contradiction. \eop

Based on Theorems 2.4 and 3.2, we now have a characterization for
local Lipschitz continuity of the best approximation operator in the
uniform norm from finite-dimensional subspaces. A weaker sufficient
condition can be found in Kovtunec [1984].

\goodbreak

\medskip
\subsect{{4.} Strong Uniqueness in Haar Spaces in the Uniform Norm}

\medskip\noindent
We start with the definition of a Haar space.

\proclaim Definition 4.1. An $n$-dimensional subspace $M$ of $C(B)$
is said to be a {\sl Haar space} if no nontrivial $m\in M$ vanishes
at more than $n-1$ distinct points of $B$.

A {\sl unicity space} is any subspace $M$ of a normed linear space
$X$ with the property that each $f\in X$ has a unique best
approximant to $f$ from $M$. The following result was proved by Haar
[1918] and is built upon earlier results of Young [1908].

\proclaim Theorem 4.1.  An $n$-dimensional subspace $M$ of $C(B)$ is
a unicity space if and only if it is a Haar space.

As a result of this theorem, the term Haar space is often applied to
any unicity space in any normed linear space. But here we use the
term Haar space as that given in the above definition. There are
many equivalent definitions of the Haar space property. Here are two
that will prove useful.

\smallskip\noindent
1) An $n$-dimensional subspace $M$ of $C(B)$ is a Haar space if and
only if $\dim M|_{\{x_1\nek x_n\}} =n$ for every choice of $n$
distinct points $x_1\nek x_n$ in $B$.

\smallskip\noindent
2) Let $m_1\nek m_n$ be a basis for $M$. Then $M$ is a Haar space if
and only if
$$\det(m_i(x_j))_{i,j=1}^n\ne 0$$
for every choice of $n$ distinct points $x_1\nek x_n$ in $B$.

\smallskip
For $n\ge 2$, there are rather restrictive conditions on $B$ needed
to ensure that $C(B)$ can contain a Haar space of dimension $n$.
Essentially, $B$ must be homeomorphic to a subset of $S^1$. Exact
conditions are the content of Mairhuber's Theorem; see Mairhuber
[1956], Sieklucki [1958], Curtis [1959], Schoenberg, Yang [1961],
and McCullough, Wulbert [1985].

When we have a Haar space, the characterization result Theorem 2.3
can be further strengthened.

\proclaim Theorem 4.2. Let $M$ be an n-dimensional Haar space on
$C(B)$. Given $f\in C(B)$, we have that $m^*$ is the best
approximant to $f$ from $M$ if and only if there exist $n+1$
distinct points $x_1\nek x_{n+1}$ in $A_{f-m^*}$ and strictly
positive values $\lam_1\nek \lam_{n+1}$ such that
$$\sum_{i=1}^{n+1} \lam_i [\sgn (f-m^*)(x_i)] m(x_i)= 0$$
for all $m\in M$.

\pf This does not differ much from Theorem 2.3. The difference is in
the claim that the $k$ therein must be $n+1$. To see this, let
$$\sigma_i = \lam_i  [\sgn (f-m^*)(x_i)],\qquad i=1\nek k.$$
From Theorem 2.3, we have
$$\sum_{i=1}^k \sigma_i m(x_i) = 0$$
for all $m\in M$. Let $m_1\nek m_n$ be any basis for $M$. Then
$$\sum_{i=1}^k \sigma_i m_j(x_i) = 0,\qquad j=1\nek n.$$
If $k\le n$, then this implies that
$${\rm rank}\, (m_j(x_i))_{i=1}^k{}_{j=1}^n <k$$
since $(\sigma_1\nek \sigma_k) \ne {\bf 0}$. But this contradicts
the fact that
$${\rm rank}\, (m_j(y_i))_{i,j=1}^{n} =n$$
for all distinct $y_1\nek y_n$ in $B$. Thus $k=n+1$. \eop

A Haar space on an interval is generally called a {\sl Chebyshev
space} (or Tchebycheff space) and often abbreviated a {\sl
$T$-space}. A basis for a $T$-space is sometimes called a {\sl
$T$-system}. For a $T$-space, because of the connectedness of the
interval, we have:

\proclaim Proposition 4.3. Let $m_1\nek m_n$ be any basis for $M$.
Then $M$ is a $T$-space if and only if
$$\eps\, \det (m_i(x_j))_{i,j=1}^n  > 0$$
for some fixed $\eps\in \{-1,1\}$ and all $x_1<\cdots <x_n$.

For $T$-spaces, we can further specialize Theorem 4.2 into a final
and more geometric form. We have

\proclaim Theorem 4.4. Let $M$ be an $n$-dimensional $T$-space on
$C[a,b]$. Given $f\in C[a,b]$, we have that $m^*$ is the best
approximant to $f$ from $M$ if and only if there exist points $a\le
x_1< \cdots < x_{n+1}\le b$ and a $\delta\in \{-1,1\}$ such that
$$(-1)^i \delta (f-m^*)(x_i) = \|f-m^*\|,\qquad i=1\nek n+1.$$

\pf ($\Leftarrow$). Assume $m^*$ exists satisfying the above. If
$m^*$ is not a best approximant to $f$ from $M$, then there exists
an $\tilm\in M$ for which
$$\|f-\tilm\| < \|f-m^*\|.$$
Then for each $i\in \{1\nek n+1\}$, we have
$$\eqalign{(-1)^i \delta (f-m^*)(x_i)
&\;=\;  \|f-m^*\| \;>\;  \|f-\tilm\|\cr &\;\ge\;  |(f-\tilm)(x_i)|
\;\ge\;  (-1)^i \delta (f-\tilm)(x_i).\cr}$$ Thus
$$(-1)^i \delta (\tilm-m^*)(x_i) >0,\qquad i=1\nek n+1.$$
It therefore follows that $\tilm-m^*$ has at least one zero in each
of the intervals $(x_i, x_{i+1})$, $i=1\nek n$. That is, $\tilm -
m^*$ has at least $n$ distinct zeros on $[a,b]$. But $\tilm-m^*\in
M$ where $M$ is an $n$-dimensional $T$-space. Thus $\tilm=m^*$,
which contradicts our hypothesis.

\smallskip\noindent
($\Rightarrow$). Assume $m^*$ is the best approximant to $f$ from
$M$. Then from Theorem 4.2, we have $n+1$ distinct points $a\le
x_1<\cdots < x_{n+1}\le b$ in $A_{f-m^*}$, and strictly positive
values $\lam_1\nek \lam _{n+1}$ such that
$$\sum_{i=1}^{n+1} \lam_i [\sgn (f-m^*)(x_i)] m(x_i)= 0$$
all $m\in M$.

\smallskip
We wish to prove that
$$ [\sgn (f-m^*)(x_i)] [\sgn (f-m^*)(x_{i+1})]<0,\qquad i=1\nek n.$$
To this end, set
$$\sigma_i = \lam_i [\sgn (f-m^*)(x_i)],\qquad i=1\nek n+1.$$
Thus, we wish to prove that
$$\sigma_i \sigma_{i+1} < 0,\qquad i=1\nek n.$$

Let $m_1\nek m_n$ be any basis for $M$. From the above, it follows
that
$$\sum_{i=1}^{n+1} \sigma_i m_j(x_i) = 0,\qquad j=1\nek n.$$
Since
$${\rm rank}\, (m_j(x_i))_{i=1}^{n+1}{}_{j=1}^n =n,$$
we have
$$\sigma_r = \alp (-1)^r \det (m_j(x_i))_{{r=1}\atop {r\ne
i}}^{n+1}{}_{j=1}^n$$ for all $r=1\nek n+1$, where $\alp\ne 0$. Thus
from Proposition 4.3,
$$\sigma_r \sigma_{r+1} < 0,\qquad r=1\nek n. \meop$$

The following was first proven in Newman, Shapiro [1963] where they
introduced the concept of (classical) strong uniqueness.

\proclaim Theorem 4.5. Let $M$ be a finite-dimensional Haar subspace
of $C(B)$. Then the unique best approximant to $f$ from $M$ is a
strongly unique best approximant to $f$ from $M$.

\pf The proof follows immediately from Theorems 2.4 and 4.2, and the
definition of a Haar space that implies that we always have $\dim
M|_{\{x_1\nek x_{n+1}\}} =n$. \eop

Thus, when $M$ is a finite-dimensional Haar space, the functional
$\gamma(f)$ is strictly positive on $C(B)$. It also has certain
simple properties. For example, it is easy to verify that $\gamma(
af-m) =\gamma(f)$ for all $a\in \RR$, $a\ne 0$, and $m\in M$. Other
properties are less obvious and for good reason.

For example, the functional $\gamma(f)$ is in general not a
continuous function of $f$. This is because $A_{f-m^*}$ is a highly
noncontinuous function of $f$.

\medskip \noindent {\bf Example.} Assume $B=[-1,1]$ and $M= \span \{1, x\}$.
$M$ is a $T$-space on $C[-1,1]$. If $A_{f-m^*} = \{-1, 0, 1\}$ with
sign of the error $+1, -1, +1$, respectively, then $\gamma(f) =
1/3$. However, if $A_{f-m^*} = \{-1, [-1/2,1/2], 1\}$ with signs
$+1, -1, +1$, respectively, i.e., $f-m^*$ attains its norm
positively at $\pm 1$ and negatively on the full segment $[-1/2,
1/2]$, then $\gamma(f)= 3/5$. It is not at all difficult to
construct an $f$ and a sequence $f_n$ of functions in $C[-1,1]$, all
with the zero function as their best approximants from $M$, and such
that $A_{f_n} = \{-1, 0, 1\}$ for all $n$, while $A_f = \{-1,
[-1/2,1/2], 1\}$.

\smallskip
While we do not have continuity of $\gamma(f)$, we do have upper
semi-continuity when $M$ is a finite-dimensional Haar space. This
result is contained in Bartelt [1975] who writes that it is due to
Phelps.

\proclaim Proposition 4.6. Assume $M$ is a finite-dimensional Haar
subspace of $C(B)$. Then the optimal strong unicity constant
$\gamma$ is an upper semi-continuous function on $C(B)$.

\pf Assume to the contrary that there exists a sequence of functions
$f_n$ in $C(B)$ that tend to $f\in C(B)$ such that
$$\gamma(f_n) \ge \gamma(f)+\eps$$
for some $\eps >0$ and all $n$. Let $m^*\in P_M(f)$ and $m_n\in
P_M(f_n)$ for each $n$. For given $m\in M$, it follows from the
definition of $\gamma(f_n)$ that
$$\|f_n -m\| \ge \|f_n -m_n\| + \gamma(f_n) \|m-m_n\|
\ge \|f_n -m_n\| + (\gamma(f)+\eps) \|m-m_n\|.$$ Since $M$ is a Haar
space, it is well-known that the best approximation operator is a
continuous operator. That is, from the fact that $f_n$ converges to
$f$, it follows that $m_n$ converges to $m^*$. Thus it also follows
that as $n\to\infty$, we have that $\|f_n-m\|$ tends to $\|f-m\|$,
$\|f_n -m_n\|$ tends to $\|f -m^*\|$ and $\|m - m_n\|$ tends to $\|m
- m^*\|$. Taking limits in the above inequality, we obtain
$$\|f -m\| \ge \|f -m^*\| + (\gamma(f)+\eps) \|m-m^*\|.$$
Since this holds for every $m\in M$, this contradicts the definition
of $\gamma(f)$. \eop

As noted previously, the reason for the lack of continuity of
$\gamma(f)$ is in the fact that $A_{f-m^*}$ is highly noncontinuous.
When we have continuity of this set, then we have continuity of the
$\gamma(f)$. To explain what we mean, let us first define a distance
between sets. In what follows, we assume that $B$ is a compact
metric space. For sets $C, D\subseteq B$, we let
$$d(C, D) = \sup_{y\in D} \inf_{x\in C} \rho (x,y)$$
where $\rho$ is a metric on $B$. There are many different such
`distances'. We will make do with this one. This result is contained
in Bartelt [1975]

\proclaim Proposition 4.7. Assume $M$ is a finite-dimensional Haar
subspace of $C(B)$, where $B$ is a compact metric space. Let $f_n$
be a sequence in $C(B)$ converging uniformly to $f\in C(B)\\ M$. Let
$m_n\in P_M(f_n)$ and $m^*\in P_M(f)$. If
$$\lim_{n\to\infty}  d(A_{f_n-m_n}, A_{f-m^*})=0,$$
then
$$\lim_{n\to\infty} \gamma(f_n) =\gamma(f).$$

\pf From Proposition 4.6, it follows that we need only prove that
for any given $\eps >0$ there exists an $N$ such that for all $n>N$,
we have
$$\gamma(f_n) +\eps >\gamma(f).$$

We start by simplifying things somewhat. Set $g_n := f_n-m_n$ and $g
:=f-m^*$. Then $A_{g_n}=A_{f_n-m_n}$ and $A_g = A_{f-m^*}$,
$\gamma(g_n) = \gamma(f_n)$ and $\gamma(g) = \gamma(f)$, and the
zero function is a best approximant from $M$ to each of the $g_n$
and $g$. Since $M$ is a Haar space, we again have that $m_n$
converges to $m^*$. Thus the $g_n$ converge to $g$, and $g\notin M$.
We will, without loss of generality, assume $\|g_n\|=\|g\|=1$ for
all $n$. We also use the fact that the unit sphere in $M$ is
uniformly equicontinuous on $B$. That is, given $\eps>0$, there
exists a $\del>0$ such that if $\rho(x, y)<\del$ then
$$|m(x)-m(y)| <\eps$$
for all $m\in M$, $\|m\|=1$. With these preliminaries, we can now
prove the desired result.

By assumption, given $\eps>0$, there exists a $\del_1 > 0$ such that
if $\rho(x,y) <\del_1$ then
$$|m(x)-m(y)| <\eps$$
for all $m\in M$, $\|m\|=1$. Similarly from the uniform continuity
of $g$, there exists a $\del_2>0$ such that if $\rho(x, y)<\del_2$
then
$$|g(x)-g(y)| < 1/2.$$
By assumption, there exists an $N_1$ such that for all $n> N_1$, we
have
$$\|g_n -g\| <1/2.$$
And finally, by assumption, given $\del>0$, there exists an $N_2$
such that for all $n> N_2$, we have $d(A_{g_n}, A_{g})<\del$, i.e.,
for each $x\in A_g$, there exists a $y\in A_{g_n}$ such that
$\rho(x, y)<\del$. For our given $\eps>0$, we therefore set $\del =
\min \{\del_1, \del_2\}$ and $N=\max\{N_1, N_2\}$.

Given any $m\in M$, $\|m\|=1$, let $x^*\in A_g$ satisfy
$$\max_{x\in A_g} [\sgn g(x)] m(x) =  [\sgn g(x^*)] m(x^*).$$
By assumption for all $n>N$, we have
$$ d(A_{g_n}, A_{g})<\del.$$
Thus there exists a $y^*\in A_{g_n}$ such that $\rho(x^*, y^*)<\del$
implying that
$$|m(x^*)-m(y^*)|<\eps.$$
In addition, as $|g(x^*)|=1$, $n>N$ and $\rho(x^*, y^*)<\del$, we
have
$$\sgn g(x^*) = \sgn g(y^*)=\sgn g_n(y^*).$$
Thus
$$\max_{x\in A_g} [\sgn g(x)] m(x) =  [\sgn g(x^*)] m(x^*) <  [\sgn g_n(y^*)] m(y^*) + \eps.$$
As this is valid for all $m\in M$, $\|m\|=1$, it implies that for
all $n>N$
$$\gamma(f) =\gamma(g) < \gamma(g_n) +\eps =\gamma(f_n) +\eps.\ \meop$$

\noindent {\bf Remark.} In Propositions 4.6 and 4.7, we did not
really use the full Haar space property. We used the fact that the
$f$ and $f_n$ all had strongly unique best approximants from $M$.

\smallskip
It would be desirable if we could bound $\gamma(f)$ from below away
from zero for all $f\in C(B)$ and thus dispense with the dependence
upon $f$. This is certainly possible if $M$ is a Haar space of
dimension 1. In this case, $M=\span \{ m\}$ where $m\in C(B)$ has no
zero, and from the definition of $\gamma(f)$, we obtain
$$\gamma(f) \ge \min_{{{m\in M}\atop {\|m\|=1}}}  |m(x)| = \min_{m\in M}  {{|m(x)|}\over {\|m\|}} >0,$$
a lower bound independent of $f$.

Moreover, if $B$ is a finite set, then we can always uniformly bound
$\gamma(f)$ from below away from zero if $M$ is a Haar space. This
follows simply from a compactness argument. However if $B$ is not a
discrete set and $M$ is a Haar space of dimension $n$, $n\ge 2$,
then $\gamma(f)$ cannot be bounded away from zero. The following is
a variation on a result proved in Cline [1973].

\proclaim Proposition 4.8. Let $B$ be a compact Hausdorff space that
is not discrete. Assume $M \subset C(B)$ is a Haar space of
dimension $n$, $n\ge 2$. Then
$$\inf_{{{f\in C(B)}\atop {\|f\|=1}}} \gamma(f)=0.$$

\pf Let $\tilx$ be any accumulation point of $B$, i.e., every
neighbourhood of $\tilx$ contains at least one other point (and thus
an infinite number of points) of $B$. As $M$ is of dimension $n$,
$n\ge 2$, there exists an $\tilm\in M$ satisfying $\tilm(\tilx)=0$
and $\|\tilm\|=1$. Given $\eps >0$, choose $n+1$ arbitrary distinct
points $\{x_i\}_{i=1}^{n+1}$ in $B$ where $|\tilm(x_i)|<\eps$. Such
points can all be chosen from a suitable neighbourhood of $\tilx$.
Since $M$ is a Haar space of dimension $n$, there exist $a_1\nek
a_{n+1}$, all nonzero, such that
$$\sum_{i=1}^{n+1} a_i m(x_i)=0\eqno{(4.1)}$$
for all $m\in M$. Construct an $f\in C(B)$ satisfying  $f(x_i) =
\sgn a_i$, $i=1\nek n+1$, $\|f\| =1$ and $\|f-\tilm\| \le 1+\eps.$
Such a construction is clearly possible since $|(f-\tilm)(x_i)| \le
|f(x_i)| + |\tilm(x_i)| < 1+ \eps $ for all $i=1\nek n+1$.

From (4.1) and the fact that $f(x_i) = \sgn a_i$, $i=1\nek n+1$, we
have that $0 \in P_M(f)$. Thus from the strong uniqueness property
of Haar spaces, we necessarily have
$$\|f-\tilm \| - \|f\| \ge \gamma(f) \|\tilm\|$$
implying
$$\eps = (1+\eps) - 1 \ge \|f-\tilm \| - \|f\| \ge \gamma(f) \|\tilm\| = \gamma(f),$$
i.e., $\gamma(f) \le \eps$. Thus
$$\inf_{{{f\in C(B)}\atop {\|f\|=1}}} \gamma(f)=0. \meop$$

\medskip
While there is no uniform strictly positive lower bound on
$\gamma(f)$, there are results concerning sets in $C(B)$ over which
the strong uniqueness constants are uniformly bounded from below by
a positive constant; see for example Bartelt, Swetits [1983],
Marinov [1983], and Bartelt, Swetits [1988]. In addition, there are
specific bounds on $\gamma(f)$ related to the structure of the
$A_{f-m^*}$. To explain how these latter bounds arise, we start with
a construction.

As usual, we assume that $M$ is an $n$-dimensional Haar space in
$C(B)$. Given any $n+1$ distinct point $x_1\nek x_{n+1}\in B$, there
exist $a_i$, all nonzero, such that
$$\sum_{i=1}^{n+1} a_i m(x_i)= 0$$
for all $m\in M$. Normalize the $a_i$ so that
$$\sum_{i=1}^{n+1} |a_i| = 1.$$
As $M$ is a Haar space of dimension $n$, there also exist $m_j\in M$
satisfying
$$ m_j(x_i) = \sgn a_i,\qquad i=1\nek n+1,\, i\ne j.$$
Thus for each $j\in \{1\nek n+1\}$,
$$0 = \sum_{i=1}^{n+1} a_i m_j(x_i) = \sum_{{i=1}\atop {i\ne j}}^{n+1}
|a_i| + a_j m_j(x_j) = (1- |a_j|) + a_j m_j(x_j),$$ and therefore
$$[\sgn a_j] m_j(x_j) = 1 - {1\over {|a_j|}} <0.$$
These $m_j$, $j=1\nek n+1$, will play a special role in determining
and bounding $\gamma(f)$. The $m_j$ were introduced in Cline [1973].
The important property of the $m_j$ is the following.

\proclaim Proposition 4.9. Assume $M$ is an $n$-dimensional Haar
subspace of $C(B)$. Given any $n+1$ distinct points
$\{x_i\}_{i=1}^{n+1}$ in $B$, let the $\{a_i\}_{i=1}^{n+1}$ and
$\{m_j\}_{j=1}^{n+1}$ be as constructed above. Then
$$\min_{{{m\in M}\atop {\|m\|=1}}} \max_{i=1\nek n+1} [\sgn a_i] m(x_i)
= \min_{j=1\nek n+1} {1 \over {\|m_j\|}}.$$

\pf We first prove that any $n$ of the above $\{m_j\}_{j=1}^{n+1}$
are linearly independent, i.e., span $M$. Assume not, and without
loss of generality let us assume that $m_2\nek m_{n+1}$ are linearly
dependent. Thus
$$\sum_{k=2}^{n+1} b_k m_k=0$$
for some nontrivial choice of $b_2\nek b_{n+1}$. At $x_1$, we have
$m_k(x_1)= \sgn a_1$ for all $k=2\nek n+1$. Thus
$$ \sum_{k=2}^{n+1} b_k =0.$$
At $x_j$, $j\ne 1$, we have
$$\eqalign{ 0 =  [\sgn a_j] \sum_{k=2}^{n+1} b_k m_k(x_j)
&\;=\;  [\sgn a_j] b_j m_j(x_j) +  \sum_{{k=2}\atop {k\ne j}}^{n+1}
b_k\cr&\;=\;  b_j \left(1 - {1\over {|a_j|}}\right) +
\sum_{{k=2}\atop {k\ne j}}^{n+1} b_k =  - {{b_j}\over
{|a_j|}}.\cr}$$ As $|a_j|> 0$, we have $b_j=0$ for all $j\ne 1$, a
contradiction.

Let $m\in M$ be normalized so that
$$\max_{i=1\nek n+1} [\sgn a_i] m(x_i) =1.$$
Note that the above maximum is always strictly positive for any
$m\in M$, $m\ne 0$, since we have $\sum_{i=1}^{n+1} a_i m(x_i)=0$,
and $m\ne 0$ can only vanish at at most $n-1$ distinct points.
Assume
$$[\sgn a_r] m(x_r) =1.$$
As any $n$ of the above $\{m_j\}_{j=1}^{n+1}$ are linearly
independent, we have
$$m = \sum_{{i=1}\atop {i\ne r}}^{n+1} c_i m_i.$$
At $x_r$, we have
$$1 = [\sgn a_r] m(x_r) = \sum_{{i=1}\atop {i\ne r}}^{n+1} c_i
[\sgn a_r] m_i(x_r) = \sum_{{i=1}\atop {i\ne r}}^{n+1} c_i.$$ We now
prove that $c_i \ge 0$ for all $i$ as above. To this end, note that
for $j\in \{1\nek n+1\}\\ \{r\}$,
$$\eqalign{ 1 \ge [\sgn a_j] m(x_j) &\;=\; [\sgn a_j] \sum_{{i=1}\atop {i\ne r}}^{n+1} c_i
m_i(x_j)  = \sum_{{i=1}\atop {i\ne r, j}}^{n+1} c_i + [\sgn a_j] c_j
m_j(x_j)\cr &\;=\; (1-c_j) + c_j (1 - {1\over {|a_j|}}) = 1 -
{{c_j}\over {|a_j|}}.\cr}$$ Thus
$$1 \ge 1 - {{c_j}\over {|a_j|}}$$
and as $|a_j|>0$, this implies that $c_j\ge 0$.

Returning to our definition of $m$, we have, using the fact that
$c_i \ge 0$ and $\sum_{{i=1}\atop {i\ne r}}^{n+1} c_i = 1$,
$$ \|m\| =\| \sum_{{i=1}\atop {i\ne r}}^{n+1} c_i m_i \|\le  \sum_{{i=1}\atop {i\ne r}}^{n+1} c_i \|m_i\|
\le \max_{i=1\nek n+1} \|m_i\|.$$ That is, for every $m\in M$
satisfying
$$\max_{i=1\nek n+1} [\sgn a_i] m(x_i) =1,$$
we have
$$ {1 \over {\|m\|}}\ge \min_{j=1\nek n+1} {1 \over {\|m_j\|}}.$$
In addition, we have
$$\max_{i=1\nek n+1} [\sgn a_i] m_j(x_i) = 1$$
for $j=1\nek n+1$. \eop

\smallskip
Let $f\in C(B)\\ M$ and $m^*$ be the best approximant to $f$ from
$M$. Recall that as $M$ is an $n$-dimensional Haar subspace of
$C(B)$, there exist $n+1$ distinct point $x_1\nek x_{n+1}\in
A_{f-m^*}$ and strictly positive $\lam_1\nek \lam_{n+1}$, with
$\sum_{i=1}^{n+1} \lam_i=1$, satisfying
$$\sum_{i=1}^{n+1} \lam_i [\sgn (f-m^*)(x_i)] m(x_i)= 0$$
for all $m\in M$. In what follows, we let the $m_j\in M$, $j=1\nek
n+1$, be as constructed above, with respect to these $x_1\nek
x_{n+1}$ and $a_i = \lam_i [\sgn (f-m^*)(x_i)]$, $i=1\nek n+1$. From
Proposition 4.9, we immediately obtain the following results.

\proclaim Theorem 4.10. Assume $M$ is an $n$-dimensional Haar
subspace of $C(B)$. Let $f\in C(B)\\ M$ and $m^*$ be the best
approximant to $f$ from $M$. Let the $m_j$, $j=1\nek n+1$, be as
constructed above. Then
$$\gamma(f) \ge \min_{j=1\nek n+1} {1 \over {\|m_j\|}}.$$
Furthermore, if $A_{f-m^*} = \{x_1\nek x_{n+1}\}$ then
$$\gamma(f) = \min_{j=1\nek n+1} {1 \over {\|m_j\|}}.$$

\pf We have that
$$\eqalign{\gamma(f) &\;=\; \min_{{{m\in M}\atop {\|m\|=1}}} \max_{x\in A_{f-m^*}}\sgn (f-m^*)(x) m(x)\cr
&\;\ge\; \min_{{{m\in M}\atop {\|m\|=1}}} \max_{i=1\nek n+1}\sgn
(f-m^*)(x_i) m(x_i).\cr}$$ From Proposition 4.9, we see that this
latter quantity equals
$$\min_{j=1\nek n+1} {1 \over {\|m_j\|}}.$$
This proves the first statement of the theorem. The second statement
follows analogously since the inequality is then an equality. \eop

The first statement in Theorem 4.10 is in Cline [1973]. The second
statement follows easily from results in Cline [1973] together with
the characterization of Theorem 2.2, as noted in Henry, Roulier
[1978]. According to Blatt [1986], the first inequality was
partially implicitly arrived at in Freud [1958].

If $A_{f-m^*}$ has exactly $n+1$ points, then from the above
construction each of the $m_j$ is ``admissible'', i.e., satisfies
$$\max_{x\in A_{f-m^*}} [\sgn (f-m^*)(x)] m_j(x) =1,\qquad j=1\nek n+1.$$
However, if $A_{f-m^*}$ has more than $n+1$ points then it may be
that some of the above constructed $m_j$ satisfy
$$\max_{x\in A_{f-m^*}} [\sgn (f-m^*)(x)] m_j(x) >1.$$
This is why equality need not hold in the first part of Theorem
4.10.

Given $A_{f-m^*}$, it is also not clear, a priori, which choices of
$n+1$ points $x_1\nek x_{n+1}$ in $A_{f-m^*}$ can be nade so as to
satisfy
$$\sum_{i=1}^{n+1} \lam_i [\sgn (f-m^*)(x_i)] m(x_i)= 0$$
for all $m\in M$ with $\lam_i>0$, $i=1\nek n+1$ (except of course if
$A_{f-m^*}$ contains only $n+1$ points). However if $B$ is a
connected interval then the ``eligible'' $x_1<\cdots <x_{n+1}$ are
exactly those for which $[\sgn (f-m^*)(x_i) ] [\sgn
(f-m^*)(x_{i+1})]<0$, $i=1\nek n$, (see Theorem 4.4).

In Bartelt, Henry [1980], there is an example of a function $f$ for
which $\gamma(f)$ is strictly greater than $\min_{j=1\nek n+1} {1
/{\|m_j\|}}$ when one varies over all possible ``eligible''
$\{x_i\}_{i=1}^{n+1}$ (see the previous paragraph). It is the
following.

\medskip\noindent {\bf Example.}
Let $M= \Pi_1$ be the polynomials of degree $\le1$ in $C[0,1]$. Let
$f\in C[0,1]$ be piecewise linear with nodes $f(0)=f(1/2)=f(1) = 1$
and $f(1/4) = f(3/4) =-1$. It can be readily verified that
$\gamma(f) = 3/5$ while the maximum of $\min_{j=1\nek n+1}
{1/{\|m_j\|}}$ where we choose 3 arbitrary alternating points is
$1/5$ (attained choosing the points $1/4, 1/2, 3/4$).

\smallskip
A simple consequence of Theorem 4.10 is the following, as was noted
in Blatt [1984b].

\proclaim Corollary 4.11. Assume $M$ is an $n$-dimensional Haar
subspace of $C(B)$. Let $f\in C(B)\\ M$ and $m^*$ be the best
approximant to $f$ from $M$. Assume $A_{f-m^*} =\{x_1\nek
x_{n+1}\}$. Then
$$\gamma(f) \le {1\over n}.$$

\pf From Theorem 4.10, we have
$$\gamma(f) = \min_{j=1\nek n+1} {1 \over {\|m_j\|}}.$$
From the construction of the $m_j$
$$|m_j(x_j)| = {1\over {\lam_j}} - 1.$$
Since $\lam_i>0$ and $\sum_{i=1}^{n+1}\lam_i=1$, there exists a
$\lam_j$ satisfying $\lam_j \le 1/(n+1)$. Thus, for this $j$,
$$\|m_j\| \ge |m_j(x_j)| = {1\over {\lam_j}} - 1 \ge n,$$
implying $\gamma(f) \le 1/n$. \eop

This bound of $1/n$ is independent of both $B$ and $M$. There is no
reason to assume that it is attainable. However Cline [1973] proved
that if $B= [-1,1]$, $M=\Pi_{n-1} = \span\{1, x\nek x^{n-1}\}$, and
the $\{x_i\}_{i=1}^{n+1}$ are the extreme points of the Chebyshev
polynomial of degree $n$, then $\gamma(f) = 1/(2n-1)$. So at least
in these cases the order $O(n)$ is attainable.

If $M=\Pi_{n-1}$ and $B=[a,b]$, then it is well known that if $f\in
C^{n}[a,b]$ and $f^{(n)}(x)$ is strictly of one sign for all $x\in
[a,b]$, then $A_{f-m^*}$ has exactly $n+1$ points. As such, it
follows from Corollary 4.11 that the strong unicity constant
associated with, for example, $\ee^x$ necessarily tends to $0$ as
$n\uparrow \infty$. It has been shown by Gehlen [1999] that if $f$
is not a polynomial, then necessarily $\liminf \gamma_n(f) =0$,
where $\gamma_n(f)$ is the strong unicity constant for $f$ on
$[a,b]$ with respect to the polynomial approximating space
$\Pi_{n-1}$. This was formally conjectured in Henry, Roulier [1978]
and is sometimes called the Poreda conjecture; see Poreda [1976]. A
lot of work went into this conjecture until it was finally proven by
Gehlen. It should be noted that the above $\liminf$ cannot be
replaced by the simple limit. There are functions $f\in C[a,b]$ for
which
$$\liminf_{n\to\infty}\gamma_n(f) =0,\qquad \limsup_{n\to \infty}\gamma_n(f)=1;$$
see Schmidt [1978]. From this example, we also see that
$\gamma_n(f)$ is not, in general, a monotone function of $n$.
Nothing is known concerning the behavior of the associated
$\gamma_n(f)$ on other dense subspaces.

\smallskip If $A_{f-m^*}$ contains more than $n+1$ points, then via
the $\|m_j\|$, we cannot necessarily calculate the exact value of
$\gamma(f)$. It is therefore natural to ask how we can calculate
$\gamma(f)$. The following, from Schmidt [1980], characterizes the
elements of $M$ that attain the minimum in the formula for the
determination of $\gamma(f)$.

\proclaim Theorem 4.12.  Assume $M$ is an $n$-dimensional Haar
subspace of $C(B)$. Let $f\in C(B)\\ M$ and $m^*$ be the best
approximant to $f$ from $M$. Assume $\tilm\in M$, $\|\tilm\|=1$,
satisfies
$$\gamma(f) = \max_{x\in A_{f-m^*}}[\sgn (f-m^*)(x)] \tilm(x).$$
Then given any $x^*\in B$ satisfying
$$|\tilm(x^*)| = \|\tilm\|,$$
there exist points $x_1\nek x_n\in A_{f-m^*}$ and strictly positive
values $\lam_1\nek \lam_{n+1}$ such that
$$[\sgn (f-m^*)(x_i) ] \tilm (x_i)= \gamma(f),\qquad i=1\nek n,$$
and
$$\sum_{i=1}^n \lam_i [\sgn (f-m^*)(x_i) ] m(x_i) - \lam_{n+1}
[\sgn \tilm(x^*) ] m(x^*)=0$$ for all $m\in M$.

\pf There necessarily exists an $\tilm\in M$, $\|\tilm\|=1$,
satisfying
$$\gamma(f) = \max_{x\in A_{f-m^*}}[\sgn (f-m^*)(x)] \tilm(x).$$
Set
$$D := \{ x: x\in A_{f-m^*}, [\sgn (f-m^*)(x)] \tilm(x)=\gamma(f) \}$$
and let $x^*$ be as above, i.e., $|\tilm(x^*)| = \|\tilm\|$. We
first claim that there does not exist an $m\in M$ satisfying
$$ [\sgn (f-m^*)(x)] m(x)<0,\qquad x\in D,$$
and
$$\tilm(x^*) m(x^*) >0.$$
For if such an $m$ existed, then for $\eps>0$ sufficiently small
$m_\eps = \tilm + \eps m$ would satisfy
$$[\sgn (f-m^*)(x)] m_\eps(x)< \gamma(f),\qquad x\in A_{f-m^*}$$
and
$$\|m_\eps\| > \|\tilm\| =1,$$
contradicting the fact that
$$\gamma(f) = \min_{{{m\in M}\atop {\|m\|=1}}} \max_{x\in
A_{f-m^*}}\sgn (f-m^*)(x) m(x) = \max_{x\in A_{f-m^*}}\sgn
(f-m^*)(x) \tilm(x).$$

Let $m_1\nek m_n$ be any basis for $M$. For each $x\in D$, set
$$M(x) = ( [\sgn (f-m^*)(x)]m_1(x)\nek [\sgn (f-m^*)(x)]m_n(x))\in
\RR^n,$$ and let
$$N(x^*) = (-[\sgn \tilm(x^*)]m_1(x) \nek -[\sgn
\tilm(x^*)]m_n(x)).$$ As $M$ is a linear subspace and there exists
no $m\in M$ such that
$$ [\sgn (f-m^*)(x)] m(x)<0,\qquad x\in D$$
and
$$-[\sgn \tilm(x^*)] m(x^*) <0,$$
it follows that ${\bf 0}\in \RR^n$ is in the convex hull of
$\{M(x)\}_{x\in D} \cup N(x^*)$. From Carath\'eodory's Theorem, this
implies that ${\bf 0}$ is a convex combination of at most $n+1$
points in this set. If these $k$ points, $1\le k \le n+1$, do not
include $N(x^*)$, then there exist $\lam_i \ge 0$, $\sum_{i=1}^{k}
\lam_i=1$, such that
$$\sum_{i=1}^{k} \lam_i [\sgn (f-m^*)(x_i)] m(x_i) = 0$$
for all $m\in M$. Substituting $m=\tilm$, and since $[\sgn
(f-m^*)(x_i)] \tilm(x_i)=\gamma(f)>0$ for $i=1\nek k$, we get a
contradiction. Thus there exist $\lam_i \ge 0$, $\sum_{i=1}^{k}
\lam_i=1$, $\lam_{k}>0$, such that
$$\sum_{i=1}^{k-1} \lam_i [\sgn (f-m^*)(x_i)] m(x_i) -\lam_{k}[\sgn \tilm(x^*)] m(x^*) = 0$$
for all $m\in M$. In fact, since $M$ is a Haar space of dimension
$n$, we must have $k=n+1$ and $\lam_i>0$ for all $i=1\nek n+1$. \eop

\medskip\noindent {\bf Example.} Let us return to the example
presented just before Corollary 4.11, namely let $M= \Pi_1$ be the
polynomials of degree $\le1$ in $C[0,1]$. Let $f\in C[0,1]$ be
piecewise linear with nodes $f(0)=f(1/2)=f(1) = 1$ and $f(1/4) =
f(3/4) =-1$. It can be readily verified that $\gamma(f) = 3/5$ while
the maximum of $\min_{j=1\nek n+1} 1/\|m_j\|$ where we choose 3
alternating points is $1/5$ (attained choosing the points $1/4, 1/2,
3/4$). The $\tilm$ that gives this $\gamma(f) =3/5$ is obtained by
considering the linear function $\tilm$ satisfying $\tilm(1/4) =
-3/5$ and $\tilm(1) =3/5$, i.e., $\tilm(x) = (8x-5)/5$.

\medskip\noindent {\bf Remark.} In the case where $B$ is
well-ordered and $M$ is a $T$-system thereon, then the $x_1<\cdots
<x_n$ as in the statement of Theorem 4.12 are such that either
$f-m^*$ strictly alternates thereon, in which case the $\tilm$ takes
its norm at $x^*<x_1$ with $\sgn \tilm(x^*) = \sgn (f-m^*)(x_1)$ or
$\tilm$ takes its norm at $x^* > x_n$ with $\sgn \tilm(x^*) = \sgn
(f-m^*)(x_n)$, or $f-m^*$ strictly alternates thereon except for
having the same sign at $x_i$ and $x_{i+1}$ for some $i\in \{1\nek
n\}$ and $\tilm$ takes its norm at $x_i < x^*<x_{i+1}$ with $\sgn
\tilm(x^*) = \sgn (f-m^*)(x_i) = \sgn (f-m^*)(x_{i+1})$.

\medskip
For which subspaces of $C(B)$ are uniqueness and strong uniqueness
equivalent? That is, which subspaces $M$ have the property that if
$f\in C(B)$ has a unique best approximant from $M$, then it is also
a strongly unique best approximant? We know that they are equivalent
if $M$ is a finite-dimensional Haar space since in this case, we
always have both uniqueness and strong uniqueness of the best
approximant. In McLaughlin, Somers [1975], it is proved that if $B$
is an interval $[a,b]$, then uniqueness and strong uniqueness for
all functions in the space are equivalent if and only if $M$ is a
Haar space, i.e., if $M$ is not a Haar space on $C[a,b]$ then there
exist functions with a unique best approximant, but not a strongly
unique best approximant. The proof is  difficult and complicated,
and is restricted to the case $B=[a,b]$.

If $B$ is a nonconnected set, there may exist subspaces $M$ of
$C(B)$ such that uniqueness and strong uniqueness are equivalent,
and yet $M$ is not a Haar space. As an example, consider $M=
\span\{(1, 0)\}$ in $\RR^2$. Then for $f=(x,y)$, we can have $f\in
M$, i.e., $y=0$, in which case we trivially have both uniqueness and
strong uniqueness of the best approximant. On the other hand, if
$f\notin M$, i.e., $y\ne 0$, then  $f$ never has a unique best
approximant. The vectors $c(1,0)$ are best approximants to $(x,y)$
from $M$ for all $c\in [x-|y|, x+|y|\,]$ as for all such $c$, we
have
$$\|(x,y) - c(1,0)\| =\max\{ |x-c|, |y|\} = |y|.$$
Thus uniqueness and strong uniqueness hold simultaneously or do not
hold at all. A more general result is the following

\proclaim Proposition 4.13. Assume $B$ is a compact Hausdorff space
and $M$ is an $n$-dimensional subspace of $C(B)$. Assume $B= B_1\cup
\cdots \cup B_k$ where each $B_j$ is both open and closed, $M|_{B_j}
= M_j$ is a Haar space of dimension $n_j$ (assuming $n_j\ge 1$) on
$B_j$, and
$$\sum_{j=1}^k n_j = n.$$
That is, we assume $M_j$ (if $n_j\ge 1$) has a basis of functions
that vanish off $B_j$. Then uniqueness and strong uniqueness are
equivalent on $C(B)$.

\pf From the above assumptions, we can write
$$M= M_1 \oplus M_2 \oplus \cdots \oplus M_k$$
and each $m\in M$ has a unique decomposition of the form
$$m=\sum_{j=1}^k m_j$$
where $m_j\in M_j$ (and $m_j$ vanishes off $B_j$). Thus
$$\|f-m\|_B =  \|f -\sum_{j=1}^k m_j\|_B
= \max_{j=1\nek k} \|f-m\|_{B_j} = \max_{j=1\nek k}
\|f-m_j\|_{B_j}.$$

We first claim that $f$ has a unique best approximant from $M$ if
and only if
$$\min_{m\in M} \|f-m\|_B = \min_{m_j\in M_j} \|f-m_j\|_{B_j}$$
for each $j\in \{1\nek k\}$ for which $n_j\ge 1$. Assume $m^*\in M$
is a best approximant to $f$ from $M$, and there exists an
$\ell\in\{1\nek k\}$ for which $n_\ell\ge 1$ and
$$\min_{m\in M} \|f-m\|_B > \min_{m_\ell\in M_\ell} \|f-m_\ell\|_{B_\ell}.$$
Let $m^*_\ell\in M_\ell$ be the best approximant to $f$ from
$M_\ell$ on $B_\ell$. As $n_\ell\ge 1$, there exists an
$\tilm_\ell\in M_\ell$, $\tilm_\ell \ne 0$, and $\tilm_\ell$
vanishes off $B_\ell$. Thus
$$\min_{m\in M} \|f-m\|_B > \|f-(m^*_\ell +\eps \tilm_\ell)\|_{B_\ell}$$
for all $\eps$ sufficiently small, implying that $m^* +\eps
\tilm_\ell$ are best approximants to $f$ from $M$ for all $\eps$
sufficiently small. Thus, if  $f$ has a unique best approximant from
$M$ then
$$\min_{m\in M} \|f-m\|_B = \min_{m_j\in M_j} \|f-m_j\|_{B_j}$$
for each $j\in \{1\nek k\}$ with $n_j\ge 1$.

Let us now assume that
$$\min_{m\in M} \|f-m\|_B = \min_{m_j\in M_j} \|f-m_j\|_{B_j}$$
for each $j\in \{1\nek k\}$ for which $n_j\ge 1$. We claim that in
this case, we always have strong uniqueness that also, of course,
implies uniqueness. Let $m^*_j\in M_j$ be the unique best
approximant to $f$ from $M_j$ on $B_j$ (that vanishes off $B_j$).
Set $m^*:=\sum_{j=1}^k m^*_j$. Thus $m^*$ is a best approximant to
$f$ from $M$ and
$$\|f-m^*\|_B = \|f-m^*_j\|_{B_j}$$
for all $j$ with $n_j\ge 1$. As $M_j$ is a Haar space (assuming
$n_j\ge 1$) on $B_j$, there exists a $\gamma_j(f)>0$ for which
$$\|f-m_j\|_{B_j} - \|f-m^*_j\|_{B_j} \ge \gamma_j(f) \|m_j- m^*_j\|_{B_j}$$
for all $m_j\in M_j$. Set
$$\gamma (f) :=\min_{{j=1\nek k}\atop {n_j\ge 1}}\gamma_j (f).$$
Let $m\in M$ and $m=\sum_{j=1}^k m_j$ where $m_j\in M_j$. For each
$j\in \{1\nek k\}$ with $n_j\ge 1$, we have, since $\|f-m^*\|_B =
\|f-m_j^*\|_{B_j}$,
$$\eqalign{\|f-m\|_B - \|f-m^*\|_B &\;\ge\; \|f-m_j\|_{B_j} - \|f-m_j^*\|_{B_j}\cr
&\;\ge\; \gamma_j(f) \|m_j- m^*_j\|_{B_j}\cr  &\;\ge\; \gamma (f)
\|m_j- m^*_j\|_{B_j}.\cr}$$ If $n_j=0$, then $m_j=m^*_j=0$ and thus
$$\|f-m\|_B - \|f-m^*\|_B \ge 0 = \gamma (f) \|m_j- m^*_j\|_{B_j}.$$
As we now have
$$\|f-m\|_B - \|f-m^*\|_B \ge \gamma (f) \|m_j- m^*_j\|_{B_j}$$
for all $j=1\nek k$, this implies that
$$\|f-m\|_B - \|f-m^*\|_B \ge \gamma (f) \|m- m^*\|_{B}. \meop$$

\smallskip
It would be interesting to determine exact necessary and sufficient
conditions for when uniqueness and strong uniqueness are equivalent
on $C(B)$.

Unfortunately, the conditions of Proposition 4.13 are not necessary.
For example, consider the space $M= \span \{m_1, m_2\}$ on $B= [0,1]
\cup [2,3]$ where
$$m_1(x) = \cases{ x, & $x\in [0,1]$\cr 1, & $x\in [2,3]$\cr}$$
and
$$m_2(x) = \cases{ 1, & $x\in [0,1]$\cr 0, & $x\in [2,3]$\cr}.$$
In this example, uniqueness and strong uniqueness are equivalent,
but this example does not satisfy the conditions of Proposition
4.13. Uniqueness and strong uniqueness are equivalent because for
uniqueness of the best approximant from $M$ to hold it is necessary
that the error function attains its norm alternately at least three
times in $[0,1]\cup [2,3]$, and two of these alternants must be in
$[0,1]$. Details are left to the reader. But in these cases, we also
get strong uniqueness. This latter fact follows from Theorem 2.4.

For $\dim M =1$, we have the following result. The argument is from
McLaughlin, Somers [1975].

\proclaim Proposition 4.14. Assume $B$ is a compact Hausdorff space
and $M=\span\{m\}$ is a $1$-dimensional subspace of $C(B)$. Assume
uniqueness and strong uniqueness from $M$ are equivalent on $C(B)$.
Then either $M$ is a Haar space on $B$, i.e., $m$ vanishes nowhere
on $B$, or $B=B_1\cup B_2$ where $B_1, B_2$ are both open and
closed, and $m$ vanishes nowhere on $B_1$ but vanishes identically
on $B_2$.

\pf To prove the proposition, it suffices to prove the following. If
$x_0\in B$ is such that $m(x_0)=0$, then $m$ vanishes identically in
a neighbourhood of $x_0$. Assume this is not the case. That is,
there exists an $x_0\in B$ such that $m(x_0)=0$ and $m$ does not
vanish identically in any neighbourhood of $x_0$. Then, we can
choose some $x_0$ such that $m(x_0)=0$ and there exists a small,
closed, nondegenerate set $N$ such that $x_0\in \partial N$ and
$\eps m(x)>0$ for all $x\in N\\ \{x_0\}$ for some $\eps \in
\{-1,1\}$. We assume, taking $cm$ if necessary for some $c\in\RR$,
that $\eps =1$ and $\|m\|=1$. Let $(x_n)$ be a sequence of points in
$N$ converging to $x_0$. Thus $m(x_n)>0$ for all $n$, and $m(x_n)$
converges to zero. Let $y\in B\\ N$ with $m(y)\ne 0$. We define
$f\in C(B)$ as follows. First set $f(y)= -\sgn m(y)$ and $f(x_n) =
1- m^2(x_n)$. Note that $f(x_0)=1$. Extend $f$ to
be in $C(B)$ satisfying $\|f\|=1$ and $|f(x)|<1$ for all $x \in B\\
\{x_0, y\}$.

As $\|f\|=1$ and $f(x_0)=1$, it follows that the $0$ function is a
best approximant to $f$ from $M$. From the construction of $f$, we
have $A_f=\{x_0, y\}$. A simple calculation gives
$$\gamma (f) =\min_{\pm} \max \{\pm m(x_0), \mp(\sgn m(y)) m(y)\} =0.$$
Thus, we do not have strong uniqueness of the best approximant to
$f$ from $M$. It remains to prove we have uniqueness.

For $\alp >0$,
$$|f(y) -\alp m(y)| = |-\sgn m(y) -\alp (\sgn m(y)) |m(y)|| = 1 + \alp |m(y)|>1,$$
implying $\alp m$ is not a best approximant for $\alp >0$. Assume
$\alp <0$ is such that $\alp m$ is a best approximant. Thus for each
$x_n$, we have
$$|f(x_n) -\alp m(x_n)| = |1 -m^2(x_n) -\alp m(x_n)|\le 1.$$
This implies that we must have
$$ m^2(x_n)+ \alp m(x_n) \ge 0$$
for all $n$. But $\alp<0$ and $m(x_n)>0$ converges to zero, which is
impossible. Thus the $0$ function is the unique best approximant to
$f$ from $M$. \eop

\goodbreak

\medskip
\subsect{{5.} Classical Strong Uniqueness in the $L^1$ Norm}

\medskip\noindent
In this section, as the title indicates, we will study strong
uniqueness when our approximating norm is $L^1$. To this end, let us
assume that $K$ is a set, $\Sigma$ a $\sigma$-field of subsets of
$K$ and $\nu$ a positive measure on $\Sigma$. By $L^1(K,\nu)$, we
mean the set of real-valued functions $f$ that are $\nu$-measurable,
and such that $|f|$ is integrable. The norm on $L^1(K,\nu)$ is given
by
$$\|f\|_1 = \int_K |f(x)|\dd\nu(x).$$

Two important sets when approximating in this $L^1(K, \nu)$ norm are
the following. For each $f\in L^1(K,\nu)$, we define its zero set as
$$Z(f) = \{ x: f(x)=0\}$$
and also
$$N(f) = K\\ Z(f).$$
Note that $Z(f)$ is $\nu$-measurable.

To understand when classical strong uniqueness might hold, we first
calculate the $\tau_+(f,g)$ of (1.1).

\proclaim Theorem 5.1. For given $f,g\in L^1(K,\nu)$, we have
$$\tau_+(f,g) = \int_{Z(f)}|g(x)|\dd\nu(x) + \int_K [\sgn f(x)] g(x) \dd\nu(x).$$

\pf Let us assume that $f,g \ne 0$. By definition,
$$\tau_+(f,g) = \lim_{t\to 0^+}  {{\|f+tg\|_1 - \|f\|_1}\over t}.$$
For $t>0$,
$$\eqalign{ {{\|f+tg\|_1 - \|f\|_1}\over t} &\;=\; {1\over t} \left[ \int_K (|f+tg|-|f|)\dd\nu\right] \cr
&\;=\; {1\over t} \left[ \int_{Z(f)} t|g|\dd\nu +\int_{N(f)}
(|f+tg|-|f|)\dd\nu  \right] \cr &\;=\; \int_{Z(f)} |g|\dd\nu +
{1\over t}\left[ \int_{N(f)} (|f+tg|-|f|)\dd\nu  \right]. \cr}$$ On
$N(f)$,
$$\left| {{|f+tg| -|f|}\over t}\right| \le |g|$$
and
$$\eqalign{ {{|f+tg| -|f|}\over t} &\;=\; {{|f+tg|^2 -|f|^2}\over {t\left[|f+tg| + |f|\right]} }\cr
&\;=\; {{2fg + t |g|^2}\over { |f+tg| + |f|}}.\cr}$$ Thus on $N(f)$,
$$\lim_{t\to 0^+} {{|f+tg| -|f|}\over t} = {{2fg}\over {2|f|}}= [\sgn f] g.$$
Applying Lebesgue's Dominated Convergence Theorem, we obtain
$$\tau_+(f,g) = \int_{Z(f)}|g|\dd\nu + \int_K [\sgn f] g \dd\nu. \meop$$

As a consequence of Theorem 5.1, Theorem 1.2 and Theorem 1.4, we
have:

\proclaim Theorem 5.2. Let $M$ be a linear subspace of $L^1(K,\nu)$
and $f\in L^1(K,\nu) \\ M$. Then $m^*\in P_M(f)$ if and only if
$$\left|\int_K [\sgn (f-m^*)(x)] m(x)\dd\nu(x)\right| \le  \int_{Z(f-m^*)} |m(x)|\dd\nu(x)$$
for all $m\in M$. Furthermore, $m^*\in M$ is a strongly unique best
approximant to $f$ from $M$ if and only if
$$\gamma(f) = \inf_{{{m\in M}\atop {\|m\|_1 =1}}} \left\{ \int_{Z(f-m^*)} |m(x)|\dd\nu(x)
+ \int_K [\sgn (f-m^*)(x)] m(x) \dd\nu(x)\right\} >0.$$

\smallskip
The first statement in this theorem can be found in James [1947,
p.~291].

Recall (see Proposition 1.6) that if $M$ is a finite-dimensional
subspace, then the infimum in the definition of $\gamma(f)$ is a
minimum, i.e., it is attained. To prove strong uniqueness, it
therefore suffices to verify that
$$ \int_{Z(f-m^*)} |m(x)|\dd\nu(x) + \int_K [\sgn (f-m^*)(x)] m(x) \dd\nu(x) >0$$
for each $m\in M$, $\|m\|=1$, and that is equivalent to showing that
$$\left|\int_K [\sgn (f-m^*)(x)] m(x)\dd\nu(x)\right| <  \int_{Z(f-m^*)} |m(x)|\dd\nu(x)$$
for all $m\in M$, $m\ne 0$.

We consider two scenarios. In the first case, $\nu$ is a discrete
positive measure with a  finite number of points of support, and in
the second case $\nu$ is a non-atomic positive measure. The two
cases radically differ.

Assuming $\nu$ is a discrete positive measure with a finite number
of points of support, we are effectively considering  approximation
in the $\ell_1^m(\bfw)$ norm where
$$\ell_1^m(\bfw):= \{\bfx=(x_1\nek x_m): \bfx\in \RR^m, \|\bfx\|_\bfw = \sum_{i=1}^m |x_i|w_i\}$$
and $\bfw = (w_1\nek w_m)$ is a fixed strictly positive vector
(weight). In this case, it transpires that uniqueness and strong
uniqueness are equivalent. It is also true that for a fixed weight
$\bfw$ most subspaces are unicity spaces, i.e., for most subspaces
$M$ there is a unique best approximant to each vector $\bfx$ from
$M$. However for some subspaces, we do not have uniqueness of the
best approximant. For example, in $\ell_1^2$ with weight $(1, 1)$,
i.e., norm $\|\bfx\|_\bfw = |x_1| + |x_2|$, and $M= \span \{
(1,1)\}$ then to each and every $\bfx\notin M$, there is no unique
best approximant. On the other hand, all the 1-dimensional subspaces
in $\ell_1^2$ other than $M=\span\{(1, \pm 1)\}$, are unicity
spaces.

The equivalence of uniqueness and strong uniqueness in this setting
is a consequence of the following, that is essentially contained in
Rivlin [1969, Theorem 3.6]; see also Watson [1980, p.~122] and
Angelos, Schmidt [1983].

\proclaim Theorem 5.3. Let $M$ be a finite-dimensional subspace of
$\ell_1^m(\bfw)$ and $\bfx\in \RR^m$. Then $\bfm^*\in M$ is the
unique best approximant to $\bfx$ from $M$ if and only if
$$\left| \sum_{i=1}^m [\sgn (x_i-m_i^*)] m_i w_i \right| < \sum_{i\in Z(\bfx-\bfm^*)} |m_i|w_i$$
for all $\bfm\in M$, $\bfm\ne {\bf 0}$. Thus every unique best
approximant to $\bfx$ from $M$ is also a strongly unique best
approximant to $\bfx$ from $M$.

\pf Assume $\bfm^*\in M$ is the unique best approximant to $\bfx$
from $M$. Thus, by the characterization of the best approximant, we
have
$$\left| \sum_{i=1}^m [\sgn (x_i-m_i^*)] m_i w_i \right| \le \sum_{i\in Z(\bfx-\bfm^*)} |m_i|w_i$$
for all $\bfm\in M$. Assume $\tilbfm\in M$, $\tilbfm\ne {\bf 0}$, is
such that
$$\left| \sum_{i=1}^m [\sgn (x_i-m_i^*)] \tilm_i w_i \right| = \sum_{i\in Z(\bfx-\bfm^*)} |\tilm_i|w_i.$$
We shall prove that the uniqueness of the best approximant implies
that such a $\tilbfm$ cannot exist. From Theorem 5.2, this implies
the equivalence of uniqueness and strong uniqueness.

Since we are in the finite-dimensional subspace $\RR^m$,
$$\min \{ |x_i-m^*_i|: i\notin Z(\bfx-\bfm^*)\} = c > 0.$$
Thus there exists an $\eps>0$ such that for all $t$ satisfying $|t|<
\eps$, we have
$$\sgn (x_i - m^*_i) = \sgn (x_i - m^*_i -t\tilm_i)$$
for each $i\notin Z(\bfx-\bfm^*)$. For each such $t$,
$$\eqalign{ \|\bfx -\bfm^* -t\tilbfm\|_{\bfw} - \|\bfx -\bfm^*\|_{\bfw} &\;=\;
\sum_{i=1}^m \left[ |x_i - m^*_i -t\tilm_i| - |x_i -
m^*_i|\right]w_i\cr &\;=\; \sum_{i=1}^m [\sgn (x_i -m^*_i)] (x_i
-m^*_i -t\tilm_i - x_i + m^*_i)w_i\cr & \phantom{12345}+ \sum_{i\in
Z(\bfx-\bfm^*)} |t\tilm_i| w_i\cr &\;=\; - t \sum_{i=1}^m [\sgn (x_i
-m^*_i)]\tilm_i w_i + |t| \sum_{i\in Z(\bfx-\bfm^*)} |\tilm_i|
w_i.\cr}$$ By assumption,
$$\left| \sum_{i=1}^m [\sgn (x_i-m_i^*)] \tilm_i w_i \right| = \sum_{i\in Z(\bfx-\bfm^*)} |\tilm_i|w_i.$$
Thus
$$ \|\bfx -\bfm^* -t\tilbfm\|_{\bfw} - \|\bfx -\bfm^*\|_{\bfw} =
 - t \sum_{i=1}^m [\sgn (x_i -m^*_i)] \tilm_i w_i + |t| \left| \sum_{i=1}^m [\sgn (x_i-m_i^*)] \tilm_i w_i \right|.$$

If $\sum_{i=1}^m [\sgn (x_i -m^*_i)] \tilm_i w_i=0$, then for all
such $t$, i.e., $|t|<\eps$, we have
$$ \|\bfx -\bfm^* -t\tilbfm\|_{\bfw} = \|\bfx -\bfm^*\|_{\bfw}.$$
If $\del\left( \sum_{i=1}^m [\sgn (x_i -m^*_i)] \tilm_i w_i\right)
>0$, $\del \in \{-1,1\}$, then
$$ \|\bfx -\bfm^* -t\tilbfm\|_{\bfw} = \|\bfx -\bfm^*\|_{\bfw}$$
for all $t$ satisfying $0\le t\del <\eps$. In either case, we have
nonuniqueness of the best approximant. \eop

In the second case, we assume that $\nu$ is a non-atomic positive
measure. In this case, neither strong uniqueness (nor uniqueness) is
always present, but it is nevertheless around. This next result can
be found in Angelos, Schmidt [1983].

\proclaim Theorem 5.4. Let $\nu$ be a non-atomic positive measure.
Let $M$ be a finite-dimensional subspace of $L^1(K,\nu)$. Then the
set of $f\in L^1(K, \nu)$ that have a strongly unique best
approximant from $M$ is dense in $L^1(K, \nu)$.

Before proving this result, we recall another characterization of
the best approximant via linear functionals that, since $\nu$ is a
non-atomic positive measure, has the following form. This result was
first proved in Phelps [1966].

\proclaim Theorem 5.5. Let $\nu$ be a non-atomic positive measure.
Let $M$ be a finite-dimensional subspace of $L^1(K,\nu)$ and $f\in
L^1(K,\nu)$. Then $m^*\in P_M(f)$ if and only if there exists an
$h\in L^\infty(K,\nu)$ satisfying
$$\eqalign{ ({\rm i}) & \, |h(x)|=1,\ {\sl all}\ x\in K;\cr
({\rm ii}) & \int_K hm\dd\nu =0, \ {\sl all}\  m\in M;\cr ({\rm
iii}) & \int_K h(f-m^*)\dd\nu = \|f-m^*\|_1.\cr}$$

\smallskip\noindent {\bf Proof of Theorem 5.4:} Let $f\in L^1(K,\nu)$ and assume that $f$
does not have a strongly unique best approximant from $M$. For
convenience, we translate $f$ by an element of $P_M(f)$ so that we
may assume, without loss of generality, that $0\in P_M(f)$. By
Theorem 5.5, there exists an $h\in L^\infty(K, \nu)$ satisfying
$$\eqalign{ ({\rm i}) & \,|h(x)|=1,\ {\rm all}\ x\in K;\cr
({\rm ii}) & \int_K hm\dd\nu =0, \ {\rm all}\  m\in M;\cr ({\rm
iii}) & \int_K hf\dd\nu = \|f\|_1.\cr}$$

\noindent Note that conditions (i) and (iii) imply that $h=\sgn f$
$\nu$-a.e.\ on $N(f)$.

Given $\eps>0$, it follows from the fact that $\nu$ is a non-atomic
positive measure and the absolute continuity of integrals that there
exists a $\del>0$ such that if $C\subseteq K$ with $\nu(C)<\del$,
then
$$\int_C|f| \, d\nu < \eps.$$
Let $C\subseteq K$, $\nu(C)<\del$, be such that if $hm\le 0$ on $C$,
then $m=0$. (Such a $C$ can be found; see Pinkus [1989, p.~22].) Set
$$f_\eps := \cases{f, & off $C$\cr 0, & on $C$.\cr}$$
Then
$$\|f-f_\eps\|_1 = \int_C |f| \dd\nu < \eps.$$

Since $\sgn f = \sgn f_\eps$ off $C$, it follows that for the above
$h$, we also have
\item{({\rm iii}$'$)} $\int_K hf_\eps\dd\nu = \|f_\eps\|_1.$

\noindent Thus $0\in P_M(f_\eps)$.

As $h=\sgn f =\sgn f_\eps$ on $N(f_\eps)\subseteq N(f)$, we have
from (ii)
$$0 = \int_K hm\dd\nu = \int_{N(f_\eps)} [\sgn f_\eps] m \dd\nu + \int_C hm\dd\nu + \int_{Z(f_\eps)\\ C} hm\dd\nu$$
for all $m\in M$. If
$$\left| \int_{N(f_\eps)} [\sgn f_\eps] m \dd\nu \right| = \int_{Z(f_\eps)} |m| \dd\nu$$
for some $m\in M$, then we must have
$$\left| \int_C hm\dd\nu + \int_{Z(f_\eps)\\ C} hm\dd\nu\right| = \int_C |m| \dd\nu + \int_{Z(f_\eps)\\ C}
|m| \dd\nu$$ implying that $hm$ is of one sign on $C$. However this
in turn implies that $m=0$. Thus
$$\left| \int_{N(f_\eps)} [\sgn f_\eps] m \dd\nu \right| < \int_{Z(f_\eps)} |m| \dd\nu$$
for all $m\in M$, $m\ne 0$, and the zero function is therefore the
strongly unique best approximant to $f_\eps$ from $M$. \eop

In Smarzewski [1988], it is proven, in the general $L^1(K,\nu)$
space without any assumptions on the measure $\nu$, that those
elements with a strongly unique best approximant are dense in the
set of elements with a unique best approximant. His proof in this
case parallels his proof of Theorem 2.5. This proof immediately
implies Theorem 5.3. It does not directly imply Theorem 5.4. For
$\nu$ a non-atomic positive measure, the density in $L^1(K, \nu)$ of
the set of functions that have a unique best approximant from $M$
goes back to Havinson, Romanova [1972] and Rozema [1974].

In Theorem 3.1, we proved that strong uniqueness at $f$ implies that
the best approximation operator from $M$ is locally Lipschitz
continuous at $f$, i.e., if $f$ has a unique best approximant from
$M$ and
$$ \|f-m\| - \|f-P_M(f)\| \ge \gamma \|m-P_M(f)\|$$
for all $m\in M$ and some $\gamma>0$, then for each $g\in X$ and any
element of $P_M(g)$, we have
$$\| P_M(f) - P_M(g)\| \le {2\over {\gamma}} \|f-g\|.$$
This implies that in the inequality
$$\| P_M(f) - P_M(g)\| \le \sig \|f-g\|,$$
we can always take the minimal (optimal) $\sig(f)$ therein to
satisfy
$$\sig(f) \le {2\over {\gamma(f)}}$$
where $\gamma(f)$ is the strong uniqueness constant for $f$. In
Theorem 3.2, we proved that on $C(B)$, assuming $M$ is a
finite-dimensional subspace, local Lipschitz continuity at $f$ of
the best approximation operator and strong uniqueness at $f$ are
equivalent. But other than the above inequality, the relationship
between $\sig(f)$ and $\gamma(f)$ is unclear.

In the case of $L^1(K, \nu)$, we again have that local Lipschitz
continuity at $f$ of the best approximation operator from $M$ and
strong uniqueness at $f$ are equivalent. Full details may be found
in Angelos, Kro\'o [1986]. We consider here the particular case of
where $\nu$ is a non-atomic positive measure. In this case, despite
the positive density statement of Theorem 5.4, for many $f$, we do
not have strong uniqueness. But when we do have strong uniqueness
then we also have
$$ {1\over \gamma(f)} \le \sig(f).$$

\proclaim Theorem 5.6. Let $\nu$ be a non-atomic positive measure.
Let $M$ be a finite-dimensional subspace of $L^1(K,\nu)$. Then $f\in
L^1(K, \nu)$ has a strongly unique best approximant from $M$ if and
only if the best approximation operator from $M$ is locally
Lipschitz continuous at $f$. Furthermore in this case, we have
$${1\over \gamma(f)} \le \sig(f)\le {2\over {\gamma(f)}}$$
where $\gamma(f)$ and $\sig(f)$ are as detailed above.

\pf If $f$ has a strongly unique best approximant from $M$, then we
have from Theorem 3.1 that the best approximation operator from $M$
is locally Lipschitz continuous at $f$, and $\sigma(f) \le
2/\gamma(f)$. Assume the converse, i.e., the best approximation
operator from $M$ is locally Lipschitz continuous at $f$. We will
prove that $f$ has a strongly unique best approximant from $M$ and
$1/\gamma(f) \le \sigma(f)$.

We assume that $f\notin M$ and the best approximant to $f$ is
unique. If either of these assumptions does not hold, then the
result follows easily. For ease of exposition, we assume that $0\in
P_M(f)$. Let $m\in M$, $m\ne 0$, and assume
$$\|f-m\|_1 -\|f\|_1 =\del >0.$$
Set
$$B := \{ x: f(x)(f-m)(x) >0\}$$
and
$$ A := K\\ (B\cup Z(f)).$$
Since $\nu$ is a non-atomic positive measure, Theorem 5.5 is valid.
Let $h\in L^\infty(K,\nu)$ be as therein.

Define $g\in L^1(K, \nu)$ by
$$g(x) :=\cases{ m(x), & $x\in A$\cr
f(x), & $ x\in B$\cr |m(x)|h(x) +m(x), & $x\in Z(f)$.\cr}$$ We shall
prove that $m\in P_M(g)$ and $\|f-g\|_1\le \del$ from which the
results will follow.

We note that as $|h|=1$ on $K$ then
$$\int_K h(g-m)\dd\nu = \|g-m\|_1$$
if we can prove that $h(g-m)\ge 0$. Now on $A$, we have $h(g-m) =
0$. On $B$, we have $h(g-m) = h(f-m)=|f-m|$, since $h=\sgn f$ on
$N(f)$, and $B\incl N(f)$ is where $\sgn f =\sgn (f-m)$. On $Z(f)$,
we have $h(g-m) = |m|$. Thus $h(g-m)\ge 0$ implying that $m\in
P_M(g)$. From Theorem 5.5, it therefore follows that $m\in P_M(g)$.

By the above, we have
$$\eqalign{\del &\;=\; \|f-m\|_1 -\|f\|_1  =  \int_K |f-m| \dd\nu -\int_K
|f|\dd\nu \cr &\;=\; \int_{A\cup B} |f-m|\dd \nu + \int_{Z(f)}
|m|\dd\nu - \int_{A\cup B} |f|\dd \nu \cr  &\;=\; \int_{A\cup B}
\!\!\!\! (f-m) \sgn(f-m)\dd \nu - \int_{A\cup B}\!\!\!\! (f-m) \sgn
f\dd \nu + \int_{Z(f)} |m|\dd\nu - \int_{A\cup B} \!\!\!\! m \sgn
f\dd \nu\cr &\;=\; \int_{A\cup B}(f-m) [\sgn(f-m) - \sgn f]\dd \nu +
\int_{Z(f)} |m|\dd\nu + \int_{Z(f)} h m \dd \nu,\cr}$$ since $\int_K
hm \dd\nu=0$ for all $m\in M$ and $h=\sgn f$ on $A\cup B$. On $B$,
we have $\sgn(f-m) =\sgn f$, while on $A$, we have $\sgn(f-m) \ne
\sgn f$. Thus
$$\del =  2 \int_{A} |f-m| \dd \nu
+ \int_{Z(f)} |m| + h m \dd \nu.$$ Now
$$\eqalign{\|f-g\|_1 &\;=\; \int_A |f-m| \dd\nu + \int_B |f-f| \dd\nu  + \int_{Z(f)} \left||m|h +m\right|
\dd\nu \cr  &\;=\; \int_A |f-m| \dd\nu + \int_{Z(f)} |m| +h m \dd\nu
\cr &\;\le\;   2\int_A |f-m| \dd\nu + \int_{Z(f)} |m| +h m \dd\nu
=\del.\cr}$$

As the best approximation operator from $M$ is locally Lipschitz
continuous at $f$, then
$$\| P_M(f) - P_M(g)\|_1 \le \sig(f) \|f-g\|_1$$
for all $g\in L^1(K,\nu)$. Substituting the above $g$, where we
recall that $0\in P_M(f)$, $m\in P_M(g)$ and $\|f-g\|_1 \le
\del=\|f-m\|_1 -\|f\|_1$, we have
$$\|m\|_1 = \|P_M(f) - P_M(g)\|_1 \le \sig(f) \|f-g\|_1 \le \sig(f) [\|f-m\|_1
-\|f\|_1],$$ i.e.,
$$\|m\|_1 \le \sig(f) [\|f-m\|_1 -\|f\|_1]$$
for every $m\in M$. This implies that $f$ has a strongly unique best
approximant from $M$, and
$$ {1\over {\gamma(f)}}\le \sig(f). \meop$$

\medskip
We now consider strong uniqueness in the one-sided $L^1(K)$ case.
Let $K$ be a compact set in $\RR^d$ satisfying $K=\overline{{\rm
int\,}K}$. We consider $f\in C(K)$ with norm
$$\|f\|_1 =\int_K |f| \dd\mu$$
where $\mu$ is a non-atomic positive finite measure with the
property that every real-valued continuous function is
$\mu$-measurable, and such that if $f\in C(K)$ satisfies $\|f\|_1=0$
then $f=0$, i.e., $\|\cdot\|_1$ is truly a norm on $C(K)$. For each
$f\in C(K)$, set
$$M(f):= \{ m: m\in M, m\le f\},$$
where $M$ is a finite-dimensional subspace of $C(K)$. We consider
the problem
$$\inf_{m\in M(f)} \|f-m\|_1 = \inf_{m\in M(f)} \int_K |f-m|\dd\mu = \inf_{m\in M(f)} \int_K f-m \dd\mu,$$
since $f-m \ge 0$ for all $m\in M(f)$. This is equivalent to
considering
$$\sup_{m\in M(f)} \int_K m \dd\mu.$$
We let $P_{M(f)}(f)$ denote the set of one-sided best approximants
to $f$ from $M$. That is,
$$P_{M(f)}(f) := \{ m^*: m^*\in M(f), \|f-m^*\|_1 \le \|f-m\|_1,\ {\rm all\ }m\in M(f)\}.$$

The following is one characterization of best one-sided
$L^1$-approximations.

\proclaim Theorem 5.7. Let $f\in C(K)$. Assume $M$ is an
$n$-dimensional subspace of $C(K)$ containing a strictly positive
function. Then a one-sided best $L^1$-approximation to $f$ from $M$
exists and $m^*\in P_{M(f)}(f)$ if and only if there exist distinct
points $\{x_i\}_{i=1}^k$ in $K$, $1\le k \le n$, and positive
numbers $\{\lam_i\}_{i=1}^k$ for which
\item {(a)} $(f-m^*)(x_i) = 0, \quad i=1\nek k.$
\item {(b)} For all $m\in M$ $$\int_K m \dd\mu = \sum_{i=1}^n \lam_i m(x_i).$$

\smallskip
We have the following result whose conditions are, unfortunately,
difficult to verify and generally rarely hold.

\proclaim Theorem 5.8. Let $M$ be an $n$-dimensional subspace of
$C(K)$ and $f\in C(K)$. Assume there exist distinct points
$\{x_i\}_{i=1}^k$ in $K$, and an $m^*\in M(f)$ satisfying
\item {(a)} $(f-m^*)(x_i) = 0, \quad i=1\nek n.$
\item {(b)} If $m\in M$ satisfies $m(x_i)=0$, $i=1\nek n$, then $m=0$.
\item {(c)} There exist strictly positive values $\lam_i$, $i=1\nek n$, such that
$$\int_K m \dd\mu = \sum_{i=1}^n \lam_i m(x_i)$$
for all $m\in M$.

\noindent {\sl In this case, $m^*$ is the unique best approximant to
$f$ from $M(f)$ and there exists a $\gamma>0$ such that
$$\|f-m\|_1 -\|f-m^*\|_1 \ge \gamma \|m-m^*\|_1$$
for all $m\in M(f)$.}

\smallskip
\pf Assume $m\in M(f)$. Then from (c) and (a), we have
$$\int_K m\dd\mu = \sum_{i=1}^n \lam_i m(x_i) \le \sum_{i=1}^n \lam_i f(x_i) =
\sum_{i=1}^n \lam_i m^*(x_i) = \int_K m^*\dd\mu.$$ If equality
holds, then we must have $m(x_i) = f(x_i) = m^*(x_i)$, $i=1\nek n$.
Thus from (b), we have $m=m^*$ and therefore $m^*$ is the unique
best approximant to $f$ from $M(f)$.

Let $m_i\in M$ satisfy $m_i(x_j) = \del_{ij}$, $i,j=1\nek n$. From
(b), it follows that we can construct these $m_i$. Thus, we can
write each $m\in M$ is the form $m = \sum_{i=1}^n m(x_i) m_i$. Set
$$\gamma := \min_{i=1\nek n} {{\lam_i}\over {\|m_i\|_1}}>0.$$
Thus $\lam_i \ge \gamma \|m_i\|_1$ for $i=1\nek n$. Assume $m\in M$
satisfies $m(x_i)\ge 0$, $i=1\nek n$. Then, using (c), we obtain
$$\int_K m\dd\mu =  \sum_{i=1}^n \lam_i m(x_i) \ge \gamma \sum_{i=1}^n m(x_i) \|m_i\|_1
 \ge \gamma \|\sum_{i=1}^n m(x_i) m_i\|_1 = \gamma \|m\|_1.$$

Now for each $m\in M(f)$, we have
$$\|f-m\|_1 - \|f-m^*\|_1 = \int_K (f-m) - (f-m^*)\dd\mu = \int_K(m^*-m)\dd\mu.$$
From (a), and since $m\in M(f)$, we have $(m^* - m)(x_i) =
(f-m)(x_i) \ge 0$, $i=1\nek n$. Thus
$$ \int_K(m^*-m)\dd\mu \ge  \gamma \|m-m^*\|_1$$
implying
$$\|f-m\|_1 -\|f-m^*\|_1 \ge \gamma \|m-m^*\|_1. \meop$$

\smallskip
Both Theorems 5.7 and 5.8 are essentially to be found in
N\"urnberger [1985]. In Section 9, we discuss non-classical strong
uniqueness in the $L^1$ setting.

\goodbreak

\subsect{{6.} Strong Uniqueness of Rational Approximation
 in the Uniform Norm}

\medskip\noindent
In this section, we consider strong uniqueness of the best uniform
approximant of functions $f\in C[a,b]$ by rational functions of the
form
$$R_{m,n}:= \{r=p/q: p\in \Pi_m, q\in \Pi_n, q(x)>0, x\in [a,b]\},$$
where $\Pi_n = \span\{1,x\nek x^n\}$. In contrast with the previous
approximation problems, we are here faced with a {\sl nonlinear} set
of approximants that makes even the question of existence of a best
approximant a nontrivial problem. The general technique, based on
bounded compactness of the unit ball in finite-dimensional spaces,
is not applicable here. Nevertheless, it turns out that existence
and uniqueness of the best rational approximant is still valid for
arbitrary $f\in C[a,b]$. However, there is an essential difference
from uniform polynomial approximation: the possibility of
degeneracy. A rational function $r=p/q\in R_{m,n}$ is said to be
{\sl degenerate} if $r=0$ or $r\in R_{m-1,n-1}$. This degeneracy of
rational approximants does not effect the uniqueness of the best
approximant, but it does spoil the continuity of the best
approximation operator and also the property of strong uniqueness.

To explain all this, let us introduce the following quantity called
the {\sl defect} of the irreducible rational function $r=p/q\in
R_{m,n}$:
$$d(r):=\cases{\min \{m-\partial p, n-\partial q\},& $ r\neq 0$;\cr n,& $r=0,$\cr}$$
where $\partial p$ is the degree of the polynomial $p$. Clearly,
$r=p/q\in R_{m,n}$ is degenerate if and only if $d(r)>0$. Moreover,
$d(r)$ is the greatest number such that $r=p/q\in
R_{m-d(r),n-d(r)}.$

First, we are going to address the questions of existence,
characterization and uniqueness of the best uniform rational
approximant. This next result is due to Walsh [1931].

\proclaim Theorem 6.1. {\rm (Existence of Best Rational
Approximation)} Any $f\in C[a,b]$ possesses a best approximant from
$R_{m,n}$.

\pf Let $r_k=p_k/q_k\in R_{m,n}$ be such that
$$\|f-r_k\|\rightarrow E(f):=\inf_{r\in R_{m,n}}||f-r||,\qquad k\rightarrow\infty.$$
Then $\|r_k\|\leq M, k\geq 0$, and we can also assume the
normalization $\|q_k\|=1, k\geq 0 $, i.e., $\|p_k\|\leq M, k\geq 0$.
Thus, passing to a subsequence, we may assume that $p_k\rightarrow
p\in \Pi_m$, $q_k\rightarrow q\in \Pi_n$, $k\rightarrow\infty$,
$q\ge 0$, and $||q||=1.$ Furthermore, the inequality $|p_k|\le M
q_k$, $k\ge 0$, that clearly holds on $[a,b]$ implies that $|p|\le M
q$, on $[a,b]$. Thus $p$ must vanish at every zero of $q$ in $[a,b]$
to at least the same multiplicity, which yields that $r:=p/q\in
R_{m,n}$ and
$$r_k(x)\rightarrow r(x), \quad x\in [a,b]\setminus Z, \eqno (6.1)$$
where $Z$ is the discrete set of zeros of $q$. Choose now an
arbitrary $\eps>0$. Since $[a,b]\setminus Z$ is dense in $[a,b]$,
the continuity of functions involved yields that for some $x\in
[a,b]\setminus Z$, we have
 $$\|f-r\| \le |f(x)-r(x)|+\eps.\eqno(6.2)$$
Since (6.1) holds for this $x$, we have for $k$ large enough
$$|f(x)-r(x)| \le |f(x)-r_k(x)|+\eps\le ||f-r_k||+\eps\le
E(f)+2\eps.$$ Combining the two inequalities, we obtain $\|f-r\| \le
E(f)+3\eps$, i.e., $\|f-r\|= E(f)$. Thus $r$ is a best approximant
to $f$. \eop

The next step will be to derive an analogue of the Alternation
Theorem for rational approximation. We start with two auxiliary
lemmas.

\proclaim Lemma 6.2.  Let $r\in R_{m,n}$ be irreducible. Then
$\Pi_m+r\Pi_n$ is a Chebyshev ($T$-) space of dimension $m+n+1-d(r)$
on $\RR$.

\pf  Let us first show that the dimension of $\Pi_m+r\Pi_n$ equals
$m+n+1-d(r)$. If $r=0$, then $d(r)=n$ and this becomes trivial. So
assume $r\neq 0$. Clearly,
$$\dim(\Pi_m+r\Pi_n)=\dim (\Pi_m)+\dim(r\Pi_n)-\dim(\Pi_m\cap r\Pi_n)=m+n+2-\dim(\Pi_m\cap r\Pi_n).$$
Furthermore, if $p_1\in \Pi_m\cap r\Pi_n$ with $r=p/q$ then
$p_1q=pq_1$ for some $q_1\in \Pi_n$. Since the polynomials $p,q$ are
relatively prime, it follows that $p_1=gp$ and $q_1=gq$, where $g$
is an arbitrary polynomial of degree at most $\min \{m-\partial p,
n-\partial q\}=d(r)$. Thus $\dim(\Pi_m\cap r\Pi_n)=d(r)+1$ and hence
$\dim(\Pi_m+r\Pi_n)=m+n+1-d(r)$.

To prove the $T$-space property, assume some $p_1+rq_1\in
\Pi_m+r\Pi_n$, $p_1\in \Pi_m$, $q_1\in \Pi_n $, $r=p/q$, has
$m+n+1-d(r)$ distinct zeros on $\RR$. Then $p_1q+pq_1$ has these
same zeros. But evidently $\partial (p_1q+pq_1)\le m+n-d(r)$,
implying that $p_1 + r q_1=0$. Thus $\Pi_m+r\Pi_n$ is a $T$-space.
\eop

\proclaim  Lemma 6.3.  Let $f\in C[a,b]$ and assume that $r^*\in
R_{m,n}$ is its best approximant, and is irreducible. Then for every
$r\in \Pi_m+r^*\Pi_n, r\neq 0$, we have
$$\max_{x\in A_{f-r^*}}[\sgn (f-r^*)(x)]\, r(x)>0. \eqno (6.3)$$

\pf Set $r^*=: p^*/q^*$ and consider an arbitrary $r = p/q$, $p\in
\Pi_m$, $q\in \Pi_n$, without the assumption that $ q > 0$ on the
interval $[a,b]$. Then for every $t > 0$, small, we have
$$ {{p^* -tp}\over {q^*+tq}}\in R_{m,n}$$
and
$$\|f-r^*\| \le \| (f-r^*) + \left(r^*- {{p^* -tp}\over {q^*+tq}}\right)\|.$$
Clearly
$$r^*- {{p^* -tp}\over {q^*+tq}} = t r_1 + O(t^2)$$
where $r_1 := {{p^*q + q^* p}\over {(q^*)^2}} \in {1\over {q^*}}
\left(\Pi_m+r^*\Pi_n\right)$. Thus by the last two relations, we
have
$$\|f - r^* + tr_1\| -\|f-r^*\| \ge C t^2$$
for $t>0$, sufficiently small, and some $C\in \RR$. This implies
that $\tau_+(f-r^*, r_1) \ge 0$ for every $r_1 \in {1\over {q^*}}
\left(\Pi_m+r^*\Pi_n\right)$. Applying Theorem 2.1, we obtain that
$$ \max_{x\in A_{f-r^*}}[\sgn (f-r^*)(x)]\, r(x)\ge 0, \qquad r\in \Pi_m+r^*\Pi_n.$$
In particular, the last relation yields by Theorem 2.2 that the zero
function is a best approximant to $f$ from $\Pi_m+r^*\Pi_n$. But, by
Lemma 6.2, this space is a Chebyshev ($T$-) space. Thus, in view of
Theorems 4.5 and 2.2, the strict inequality (6.3) must hold. \eop

We can now verify the Alternation Theorem which is the analogue of
Theorem 4.4 for rational approximation. It may be found in Achieser
[1930].

\proclaim  Theorem 6.4.  {\rm (Characterization of Best Rational
Approximation)} Let $f\in C[a,b]$, and $r^*\in R_{m,n}$ be
irreducible. Then the following statements are equivalent:
\smallskip\noindent
(i) $r^*$ is a best approximant to $f$;
\smallskip\noindent
(ii) the error function $f-r^*$ equioscillates on at least
$m+n+2-d(r^*)$ points in $[a,b]$.

\pf (i) $\Rightarrow$ (ii). By Lemma 6.3 and Theorem 2.2, the $0$
function is a best approximant to $f-r^*$ from  $\Pi_m+r^*\Pi_n$.
Since, by Lemma 6.2, this space is a $T$-space of dimension
$m+n+1-d(r^*)$, Theorem 4.4 implies that (ii) must hold.

\smallskip\noindent
(ii) $\Rightarrow$ (i). Assume $\|f-r\| < \|f-r^*\|$ for some $r\in
R_{m,n}$. Since $f-r^*$ equioscillates on at least $m+n+2-d(r^*)$
points in $[a,b]$, it follows that $r-r^*$ has at least
$m+n+1-d(r^*)$ distinct zeros in $[a,b]$. Setting $r^*=:p^*/q^*$ and
$r=:p/q$, we obtain that the polynomial $p^*q-pq^*$ has at least
$m+n+1-d(r^*)$ distinct zeros in $[a,b]$. But clearly $\partial
(p^*q-pq^*)\leq \max\{\partial p^*+n, \partial q^*+m\}\leq
m+n-d(r^*)$, a contradiction. Thus $r^*$ is a best approximant. \eop

As a byproduct of Theorem 6.4, we obtain that every $f\in C[a,b]$
must have a {\sl unique} best approximant from $R_{m,n}$; see
Achieser [1930] and Achieser [1947].

\proclaim Corollary 6.5. {\rm (Uniqueness of Best Rational
Approximation)} Every $f\in C[a,b]$ possesses a unique best
approximant from $R_{m,n}$.

\pf Assume that $r^*=p^*/q^*$ is a best approximant to $f$ and
$\|f-r\|=\|f-r^*\|$ for some $r=p/q\in R_{m,n}$.  By the same
argument as in the proof of (ii) $\Rightarrow$ (i) in Theorem 6.4,
it follows that the polynomial $g:=p^*q-pq^*$ of degree $\leq
m+n-d(r^*)$ has at least $N:=m+n+2-d(r^*)$ points of weak sign
change. That is, for some points $a\leq x_1< \cdots <x_N\leq b$, we
have $g(x_j)g(x_{j-1})\leq 0$ for every $2\le j\le N$. Since
$\partial g\le N-2$, the Lagrange Interpolation Formula yields that
$$\sum _{j=1}^{N}g(x_j){1\over {\omega^{'}(x_j)}}=0,$$
where $\omega (x):=(x-x_1)\cdots (x-x_N)$. Since $ \omega^{'}(x_j)$
also alternate in sign, we must have $ g(x_j)=0$, $1\leq j\leq N$,
i.e., $g=0$, and hence $r^*=r$. \eop

We now have the necessary prerequisites to verify that strong
uniqueness of best rational approximation from $R_{m,n}$ holds if
and only if the best approximant is non-degenerate. Moreover, our
next theorem shows that the continuity of the best rational
approximant operator also fails to hold exactly in the case of
degeneracy of the best approximant. This theorem combines results of
Maehly, Witzgall [1960], Cheney, Loeb [1964] and Werner [1964] (with
regards to continuity of the best approximation operator) and
Cheney, Loeb [1964] and  Cheney [1966] (with regards to strong
uniqueness).

\proclaim Theorem 6.6.  Let $f\in C[a,b]\\ R_{m,n}$ and let $r^*=
p^*/q^* \in R_{m,n}$ be its best rational approximant. Then the
following statements are equivalent:
\smallskip\noindent
(i) $d(r^*)=0$;
\smallskip\noindent
(ii) $r^*$ is a strongly unique best approximant;
\smallskip\noindent
(iii) the operator of best rational approximation from $R_{m,n}$ is
continuous at $f$.

\pf  (i) $\Rightarrow$ (ii). For an arbitrary $r\in R_{m,n}$, $r\neq
r^*$, set
$$\gamma (r):= {{\|f-r\|-\|f-r^*\|}\over {\|r-r^*\|}}.$$
In order to verify strong uniqueness, we need to show that there
exists $\gamma>0$ such that $\gamma (r)\geq \gamma$ for every $r\in
R_{m,n}$, $r\neq r^*$. Assume, to the contrary, that there exists a
sequence $r_k=p_k/q_k \in R_{m,n}$, $r_k\neq r^*$, such that $\gamma
(r_k)\rightarrow 0.$ We may assume, without loss of generality, that
$\|p_k\|+\|q_k\|=1$, $k\in \NN$. Passing, if necessary, to a
subsequence, it can also  be assumed that for some $p\in \Pi_m$,
$q\in \Pi_n$, we have $p_k\rightarrow p$, $q_k\rightarrow q$, as
$k\rightarrow \infty$, uniformly on $[a,b]$. In particular,
$\|p\|+\|q\|=1$.  Assume, in addition, that $\|p^*\|+\|q^*\|=1$
(where $r^*=p^*/q^*$).

Consider now an arbitrary $x\in A_{f-r^*}$. Then
$$\eqalign{\gamma (r_k)\|r^*-r_k\|&\;=\;\|f-r_k\|-\|f-r^*\|\cr
&\;\geq\; (f-r_k)(x) \sgn (f-r^*)(x)-(f-r^*)(x) \sgn (f-r^*)(x)\cr
&\;=\;(r^*-r_k)(x) \sgn (f-r^*)(x).\cr}$$ Thus using the fact that
$q_k>0$ and $\|q_k\|\leq 1$, we obtain, for any $x\in A_{f-r^*}$,
$$[\sgn(f-r^*)(x)] (q_kr^*-p_k)(x) \le \gamma (r_k)\|r^*-r_k\|q_k(x)\leq
\gamma (r_k)\|r^*-r_k\|. \eqno (6.4)$$ Note that since $\gamma
(r_k)\rightarrow 0$, we must have $\|r_k\|\le M$, $k\in \NN$, for
some $M>0$. (Otherwise $\|r_k\|\rightarrow \infty$ for a subsequence
of $k$'s that in turn would yield that $\gamma (r_k)\rightarrow 1$
for the same subsequence.) This means that the right hand side of
(6.4) tends to 0 as $k\rightarrow \infty$ yielding
$$[\sgn (f-r^*)(x)] (qr^*-p)(x) \le 0,\qquad x\in A_{f-r^*}.$$
In view of Lemma 6.3, this implies that $qr^*=p$, i.e., $p\in
\Pi_m\cap r^*\Pi_n$. On the other hand, using the fact that
$d(r^*)=0$, we have $\dim(\Pi_m\cap r^*\Pi_n)=d(r^*)+1=1$ (see the
proof of Lemma 6.2 for details). Since evidently $p^*$ is also an
element of $\Pi_m\cap r^*\Pi_n$, we must have $p=cp^*$ for some
nonzero constant $c$. But then in view of  $qr^*=p$, the relation
$q=cq^*$ also holds. Recalling the conditions $\|p\|+\|q\|= \|p^*\|
+ \|q^*\| = 1$, and $q, q^*>0$, we clearly have $c=1$, i.e.,
$p=p^*$, $q=q^*$. Since $q^*>0$ on $[a,b]$, it follows that for some
$\delta>0$ the relation $q_k(x)\ge \delta$ holds for all $x\in
[a,b]$, for $k$ sufficiently large.

Furthermore, since inequality (6.3) holds for every $r$ in the
finite-dimensional linear space $\Pi_m+r^*\Pi_n$, and because the
left hand side of (6.3) is a continuous function of $r$, we obtain
that there exists an $\eta>0$ such that
$$\max_{x\in A_{f-r^*}}[\sgn (f-r^*)(x)] \,r(x)>\eta ||r||,\qquad r\in \Pi_m+r^*\Pi_n.$$
Thus, in particular, there exists an $x^*_k\in A_{f-r^*}$ for which
$$ [\sgn(f-r^*)(x^*_k)]\, (q_kr^*-p_k)(x^*_k) \ge \eta \|q_kr^*-p_k \| .$$
Using this last inequality together with (6.4) applied with
$x=x^*_k$, we arrive at
$$\eta \|q_kr^*-p_k \|\le \gamma (r_k)\|r^*-r_k\|=\gamma
(r_k)\|{{1}\over {q_k}}(q_kr^* -p_k)\|\le {{\gamma (r_k)}\over
{\delta}}\|q_kr^*-p_k\|.$$ Since $r_k\neq r^*$, we can divide both
sides of this inequality by $\|q_k r^* -p_k\|$  yielding $\eta\leq
{{\gamma (r_k)}\over {\delta}}$. But this contradicts the assumption
that $\gamma (r_k)\rightarrow 0.$ This completes the proof of strong
uniqueness.

\smallskip\noindent
(ii) $\Rightarrow$ (iii). This implication is trivial since strong
uniqueness implies Lipschitz continuity of the best approximation
operator; see Theorem 3.1.

\smallskip\noindent
(iii) $\Rightarrow$ (i). Assume that to the contrary $d(r^*)>0$. We
separate the proof into two cases.

\smallskip\noindent
{\it Case 1}. $f-r^*$ has less than $m+n+2$ points of
equioscillation. In this case, $f-r^*$ has $m+n+2-d$ points of
equioscillation for some $0<d\leq d(r^*)$. Evidently, for any
$\eps>0$, there exists an $r_1\in R_{m,n}$ such that $\|r^*-r_1\|\le
\eps$ and $d(r_1)=0$. Set $g:=f+(r_1-r^*)$, and denote by $r_2\in
R_{m,n}$ the best approximant to $g$. Clearly, $r_2\neq r_1$ since
otherwise the error function $g-r_2$ would have less than $m+n+2$
points of equioscillation, in contradiction to the fact that
$d(r_1)=0$ and Theorem 6.4. Thus by the uniqueness of best rational
approximation (Corollary 6.5),
$$\|g-r_2\|<\|g-r_1\|=\|f-r^*\|. \eqno (6.5)$$
Now, the function $g-r_1=f-r^*$ has $m+n+2-d$ points of
equioscillation, which by (6.5) yields that $r_1-r_2$ must have at
least $m+n+1-d$ zeros in $[a,b]$. Set $r_1=:p_1/q_1$,
$r_2=:p_2/q_2$. By the above observation that the nonzero polynomial
$p_1q_2-p_2q_1$ has at least $m+n+1-d$ zeros, it follows that
$\partial (p_1q_2-p_2q_1)\geq m+n+1-d$.  On the other hand,
$$\partial (p_1q_2-p_2q_1)\le \max \{m+\partial q_2, n+\partial
p_2\} \le m+n-d(r_2).$$ Thus, we obtain $d(r_2)\leq d-1.$ This
implies that $g-r_2$ must have at least $m+n+2-d(r_2)\geq m+n+3-d$
points of equioscillation. Moreover, $\|f-g\|=\|r^*-r_1\|\leq \eps$
and hence in view of the continuity of the best approximation
operator at $f$, we have that $g-r_2$ tends uniformly to $f-r^*$ as
$\eps\rightarrow 0$. Recall that $f-r^*$ has exactly $m+n+2-d$
points of equioscillation, while $g-r_2$ has at least $ m+n+3-d$
points of equioscillation. This leads to a contradiction as, in the
limit, we cannot lose points of equioscillation.

\smallskip\noindent
{\it Case 2}. We now assume that $f-r^*$ possesses at least $m+n+2$
points of equioscillation. We also, for simplicity, set
$[a,b]=[0,1].$ We first consider the subcase
$|(f-r^*)(0)|<\|f-r^*\|.$ There then exist $0<\alpha<1$ and
$\delta>0$ such that
$$|(f-r^*)(x)|<\alpha\|f-r^*\|,\qquad 0\leq x\leq \delta. \eqno (6.6)$$
Note that all points of equioscillation of the error curve $f-r^*$
must be in $[\delta,1]$. Set $$A:=(1-\alpha){{\|f-r^*\|}\over
{\|1/q^*\|}},\quad r_k:={ {A}\over {q^*(1+kx)}}, \quad
r_k^{**}:=r^*-r_k.$$ Since $r^*$ is degenerate, it is easy to see
that $r_k, r_k^{**}\in R_{m,n}$. Consider now the functions
$$g_k:={A\over {q^*}}\min \{{1 \over{1+kx}},{1\over
{1+k\delta}}\},\quad f_k:=f-g_k.$$ Evidently, $g_k=r_k$ if $x\geq
\delta$. Furthermore, whenever $0\leq x\leq \delta$, we have by the
choice of $A$
$$|(r_k-g_k)(x)|={{Ak(\delta-x)}\over {q^*(1+kx)(1+k\delta)}}\leq
{{Ak\delta}\over {q^*(1+k\delta)}}\leq (1-\alpha){{k\delta
\|f-r^*\|}\over {1+k\delta}}\leq (1-\alpha)\|f-r^*\|.$$ Using the
last inequality together with (6.6), we obtain for $0\leq x\leq
\delta$
$$|(f_k-r_k^{**})(x)|=|(f-r^*)(x)+(r_k-g_k)(x)| <
\alpha \|f-r^*\|+(1-\alpha)\|f-r^*\|=\|f-r^*\|.$$ Recalling that
$f_k-r_k^{**}=f-r^*$ whenever $x\geq \delta$ and all points of
equioscillation of the error curve $f-r^*$ are in $[\delta,1]$, it
follows that the error curve $f_k-r_k^{**}$ equioscillates at least
$m+n+2$ times, i.e., by Theorem 6.4, $r_k^{**}$ is the best
approximant to $f_k$. On the other hand, we clearly have
$$\|f_k-f\|=\|g_k\|\leq {{A\|1/q^*\|}\over {1+k\delta}}\rightarrow
0, \qquad k\rightarrow \infty,$$ while for every $k$
$$\|r_k^{**}-r^*\|=\|r_k\|\geq {A\over {q^*(0)}}.$$
Thus the continuity of the operator of best rational approximation
fails to hold at $f$. This verifies the needed statement in the case
$|(f-r^*)(0)|<\|f-r^*\|.$

Finally, if $|(f-r^*)(0)|=\|f-r^*\|$, we can introduce a small
perturbation to the function $f$ resulting in $f_\eps \in C[0,1]$
such that $\|f-f_\eps\|<\eps$, $r^*$ is still the best approximant
to $f_\eps$ and $|(f_\eps-r^*)(0)|<\|f_\eps-r^*\|.$ By the above
argument, the continuity of best rational approximation fails
uniformly at $f_\eps$ by at least $A/q^*(0)$, and hence it will fail
at $f$ as well. \eop

The ideas in this section have been further generalized to other
non-linear families, cf., Barrar, Loeb [1970], Braess [1973], and
Braess [1986].

\goodbreak

\bigskip
\sectwopn{{Part II.} Non-Classical Strong Uniqueness}

\subsect{{7.} Uniformly Convex Space}

\medskip\noindent
Let $X$ be a uniformly convex space (of dimension at least 2) with
norm $\|\cdot\|$. In such a case, we know that we always have
uniqueness of the best approximant. But we rarely have strong
uniqueness in the classical sense. As such, we look for a different
non-classical form of strong uniqueness. This we do in the present
context using the modulus of convexity. The {\sl modulus of
convexity} on $X$ is defined as follows. For $\eps\in (0,2]$, we set
$$\del(\eps) :=\inf \{ 1 - {{\|f+g\|}\over 2}: \|f\|, \|g\| = 1, \|f-g\|\ge \eps\}.$$
We recall that $X$ is uniformly convex if $\del(\eps)>0$ for every
$\eps>0$. If $X$ is a uniformly convex space, then the following
estimate is valid; see Bj\"{o}rnest{\aa}l [1979].

\proclaim Theorem 7.1. Let $M$ be a linear subspace of $X$. Given
$f\in X$ assume that $m^*$ is the best approximant to $f$ from $M$.
Then for all $m\in M$, $m\ne m^*$, we have
$$\|f-m\| - \|f-m^*\| \ge  \|f-m\|\, \del \left({{\|m-m^*\|}\over {\|f-m\|}}\right).$$

\smallskip\noindent
{\bf Remark.} Note that as $\|m-m^*\| \le \|f-m\| + \|f-m^*\|$, it
follows that
$${{\|m-m^*\|}\over {\|f-m\|}} \le {{\|f-m\| + \|f-m^*\|}\over
{\|f-m\|}}= 1+ {{\|f-m^*\|}\over {\|f-m\|}} < 2$$ since $\|f-m^*\| <
\|f-m\|$. Thus $ \del({{\|m-m^*\|}\over {\|f-m\|}})$ is
well-defined.

\smallskip\noindent
{\bf Remark.} $\del$ is a non-decreasing function. If
$$\|f-m\| -\|f-m^*\| \le \sig ,$$
then, since $\|f-m\|\ge \|f-m^*\|$, we obtain
$$\|f-m\| - \|f-m^*\| \ge  \|f-m^*\|\, \del\left({{\|m-m^*\|}\over {\|f-m^*\|+\sig}}\right),$$
which is a non-classical strong uniqueness inequality.

\smallskip\pf
As $M$ is a linear subspace, the above claim is equivalent to
verifying
$$\|h-m\|-\|h\| \ge \|h-m\|\, \del\left({{\|m\|}\over {\|h-m\|}}\right)$$
for any $h$ such that $0$ is a best approximant to $h$ from $M$.

It easily follows from the definition of $\del$ that for any $f,
g\in X$ satisfying $\|f\|, \|g\|\le r$, $r>0$, we have
$$\left\| {{f+g}\over 2}\right\| \le r \left[ 1 - \del
\left({{\|f-g\|}\over r}\right)\right].$$

Let, in the above, $f=h$ and $g=h-m$, and $r = \|h-m\|\ge \|h\|$.
Substituting, we therefore obtain
$$\left\| {{2h-m}\over 2}\right\| \le \|h-m\| \left[ 1 - \del
\left({{\|m\|}\over {\|h-m\|}}\right)\right].$$ Note that as $0$ is
a best approximant to $h$ from $M$, it follows that
$$\left\| {{2h-m}\over 2}\right\| =\left\|h -{m\over 2}\right\|
\ge \|h\|.$$ Thus
$$ \|h\| \le \|h-m\| \left[ 1 - \del
\left({{\|m\|}\over {\|h-m\|}}\right)\right]$$ that immediately
translates into the desired inequality. \eop

For a similar result, see Wegmann [1975], and also Lin [1989]. From
Theorem 7.1, we may obtain the following result that was proven in
Smarzewski [1986], Prus, Smarzewski [1987], Smarzewski [1987] by
different methods.

\proclaim Corollary 7.2. Let $M$ be linear subspace of $X$. Assume
$$\del(\eps) \ge C \eps^q$$
for $\eps\in (0,2]$ and $q\ge 1$. Then for $m^*$, the best
approximant to $f$ from $M$, we have
$$\| f-m\|^q - \|f-m^*\|^q \ge C \|m-m^*\|^q,$$
for all $m\in M$, $m\ne m^*$. Furthermore, if
$$\|f-m\| -\|f-m^*\| \le \sig$$
then we also have
$$\|f-m\| - \|f-m^*\| \ge C k \|m-m^*\|^q,$$
where $k:= \|f-m^*\|/(\|f-m^*\| +\sig)^q$.

\pf The first inequality follows immediately from Theorem 7.1. If
$q=1$, there is nothing to prove. For $q>1$, we have
$$\|f-m\| - \|f-m^*\| \ge  \|f-m\| \del\left({{\|m-m^*\|}\over {\|f-m\|}}\right)
\ge \|f-m\| C \left({{\|m-m^*\|}\over {\|f-m\|}}\right)^q.$$
Multiply through by $\|f-m\|^{q-1}$ and use the fact that $\|f-m^*\|
\le \|f-m\|$.

The second inequality follows from the second remark after the
statement of Theorem 7.1.\eop

\smallskip\noindent
{\bf Remark.} Moduli of convexity satisfying
$$\del(\eps) \ge C \eps^q$$
are called moduli of convexity of power type $q$. It is known that
for any modulus of convexity, we can only have the above holding for
$q\ge 2$. (The condition $q\ge 1$ in the above corollary is
superfluous.)

\medskip
It is also known, see Hanner [1956], that if $\del_p(\eps)$ denotes
the modulus of convexity of $L^p$ or $\ell_p$, then we have
$$\del_p(\eps) =\cases{{{(p-1)\eps^2}\over 8} +o(\eps^2), & $1<p<2$,\cr
{{\eps^p}\over {p2^p}} + o(\eps^p), & $2\le p <\infty$.\cr}$$ In
addition, in Hilbert spaces, we always have $\del(\eps) = 1-
(1-\eps^2/4)^{1/2} = \eps^2/8 + O(\eps^4)$ and Hilbert spaces are
the most convex in that $\del(\eps) \le  1- (1-\eps^2/4)^{1/2}$
holds for any uniformly convex Banach space. For $L_p$, we in fact
have
$$\del_p(\eps) \ge c_p \eps^q$$
for all $\eps\in [0,2]$ where $q=\max\{2, p\}$ and
$$c_p := \cases{{{p-1}\over 8}, & $1<p<2$,\cr
{1\over {p2^p}}, & $2\le p <\infty$.\cr}$$

The following result can be found in Smarzewski [1987]. Assume
$\phi$ is an increasing convex function defined on $\RR_+$
satisfying $\phi(0)=0$. Assume in addition that
$$\del(\eps) \ge K \phi(\eps)$$
for $\eps\in [0,2]$, and $\phi$ is submultiplicative in that there
exists a constant $L>0$ such that for all $t,s\in \RR_+$, we have
$$\phi(ts) \le L \phi(t)\phi(s).$$ Then the following holds.

\proclaim Proposition 7.3. Let $M$ be a linear subspace of the
uniformly convex space $X$. Assume $\phi$ is as above. If $m^*$ is
the best approximant to $f$ from $M$ then
$$\phi(\|f-m\|) - \phi(\|f-m^*\|) \ge K L^{-1} \phi(\|m-m^*\|).$$
for all $m\in M$, $m\ne m^*$.

\pf Since $\phi$ is convex and $\phi(0)=0$, we have for any $\lam\in
[0,1]$ that
$$\phi(\lam t) = \phi(\lam t + (1-\lam)0)\le \lam \phi(t).$$
Thus
$$\phi(\|f-m^*\|) = \phi \left({{\|f-m^*\|}\over {\|f-m\|}}
\|f-m\|\right)\le {{\|f-m^*\|}\over {\|f-m\|}} \phi(\|f-m\|).$$
Therefore, applying Theorem 7.1,
$$\eqalign{ \phi(\|f-m\|)- \phi(\|f-m^*\|)&\;\ge\;  \phi(\|f-m\|) -  {{\|f-m^*\|}\over {\|f-m\|}}
\phi(\|f-m\|)\cr &\;=\; {{\|f-m\|- \|f-m^*\|}\over {\|f-m\|}}
\phi(\|f-m\|)\cr &\;\ge\; \del\left({{\|m-m^*\|}\over
{\|f-m\|}}\right) \phi(\|f-m\|)\cr &\;\ge\; K \phi
\left({{\|m-m^*\|}\over {\|f-m\|}}\right) \phi(\|f-m\|)\cr &\;\ge\;
KL^{-1} \phi(\|m-m^*\|). \meop\cr}$$

\smallskip
Corollary 7.2 easily follows from Proposition 7.3 by choosing
$\phi(t)=t^q$ (where $L=1$). A different approach is the following
which may also be found in Smarzewski [1987].

\proclaim Theorem 7.4. Assume there exists a positive constant $K$
and an increasing nonnegative function $\phi$ defined on $\RR_+$
such that
$$\phi\left( \left\| {{g+h}\over 2}\right\|\right)\le {1\over 2}
[\phi(\|g\|) + \phi(\|h\|)] - K\phi (\|g-h\|)$$ holds for all $f,g
\in X$. Let $M$ be a linear subspace of $X$. If $m^*$ is the best
approximant to $f$ from $M$, then
$$\phi(\|f-m\|) - \phi(\|f-m^*\|) \ge 2K \phi(\|m-m^*\|)$$
for all $m\in M$, $m\ne m^*$.

\pf As $m^*$ is the best approximant to $f$ from $M$, we have that
$$\left\| {{f-m + f-m^*}\over 2}\right\| = \left\| f-
{{m+m^*}\over 2}\right\| \ge \|f-m^*\|.$$ Set $g:= f-m^*$ and
$h=f-m$. Then, as $\phi$ is increasing,
$$\phi(\|f-m^*\|) \le \phi\left(\left\| {{f-m + f-m^*}\over
2}\right\|\right) \le {1\over 2}[\phi(\|f-m^*\|) + \phi(\|f-m\|)] -
K\phi (\|m-m^*\|).$$ The desired inequality follows. \eop

For $p\ge 2$, it is known that for $g, h\in L^p$, we have
$$\|g+h\|_p^p + \|g-h\|_p^p \le 2^{p-1} (\|g\|_p^p + \|h\|_p^p).$$
This inequality is called Clarkson's Inequality; see Clarkson
[1936]. Thus, it follows that we can apply Theorem 7.4 with
$\phi(t)=t^p$ and $K= 2^{-p}$ giving
$$ \|f-m\|_p^p - \|f-m^*\|_p^p \ge 2^{1-p} \|m-m^*\|_p^p.$$
This improves somewhat the constant $2^{-p}/p$ obtained as a
consequence of Corollary 7.2.  For a better constant, see Smarzewski
[1986]. For $1<p<2$, the inequality
$$ \left\| {{g+h}\over 2}\right\|_p^2\le {1\over 2}
[\|g\|_p^2 + \|h\|_p^2] - {{p(p-1)}\over 8}\|g-h\|_p^2$$ holds for
all $f,g \in X$. This is called Meir's Inequality; see Meir [1984].
Thus, we can take $\phi(t)=t^2$ and $K= p(p-1)/8$ giving
$$ \|f-m\|_p^2 - \|f-m^*\|_p^2 \ge {{p(p-1)}\over 4} \|m-m^*\|_p^2.$$
This improves somewhat the constant $(p-1)/8$ obtained as a
consequence of Corollary 7.2.

For $p=2$, both these improved constants are $1/2$. However this is
not the optimal constant, as in any Hilbert space, we have
$$\|f-m\|^2 - \|f-m^*\|^2 = \|m-m^*\|^2$$
for all $m\in M$ assuming $m^*$ is a best approximant to $f$ from
the linear subspace $M$.

For other results in this direction, see Angelos, Egger [1984] and
Egger, Taylor [1989].

\goodbreak

\subsect{{8.} The Uniform Norm Revisited}

\medskip\noindent
In this section, we return to the uniform norm over an interval
$[a,b]$. It transpires that if we restrict our approximation to only
``smooth'' functions, then a weaker condition than the Haar
condition is both necessary and sufficient to insure uniqueness of
the best approximant to each ``smooth'' function if the
finite-dimensional subspace $M$ is itself ``smooth''. The following
result, from Garkavi [1959], see also Kro\'o [1984], characterizes
unicity spaces with respect to $C^{r}[a,b]$ for $r\ge 1$. Note that
this characterization is independent of $r$. The characterization
critically uses the fact (independent of $r \ge 1$) that if $x\in
A_{f-m^*} \cap (a,b)$ and $f-m^*\in C^{r}[a,b]$ then we must have
$(f-m^*)'(x)=0$.

\proclaim Theorem 8.1. Let $r$ be a fixed positive integer, and let
$M$ be an $n$-dimensional subspace of $C^{r}[a,b]$. A necessary and
sufficient condition for the uniqueness of the best approximant from
$M$ to each $f\in C^{r}[a,b]$ is that the following does {\bf not}
hold:
\smallskip\noindent
There exist $k$ distinct points $\{x_i\}_{i=1}^k$ in $[a,b]$, $1\le
k \le n+1$, nonzero values $\{\sig_i\}_{i=1}^k$, and an $\tilm\in M\\
\{0\}$ such that
\item{(a)} $\sum_{i=1}^k \sig_i m(x_i) = 0$ all $m\in M$.
\item{(b)} $\tilm(x_i)=0$, $i=1\nek k$.
\item{(c)} ${\tilm}'(x_i)=0$ if $x_i\in (a,b)$.

\medskip
In Section 4, we defined a Chebyshev ($T$-) system as a system of
functions $m_1\nek m_n$ defined on $C[a,b]$ such that no nontrivial
$m\in M=\span \{m_1\nek m_n\}$ vanishes at more than $n-1$ distinct
points of $[a,b]$. We often use the term $T$-system interchangeably
for both $m_1\nek m_n$ and for the space $M$. We say that $M=\span\{
m_1\nek m_n\}$ is an $ET_2$-system if $M\subset C^{1}[a,b]$ and no
nontrivial $m\in M$ vanishes at more than $n-1$ distinct points of
$[a,b]$, where we count multiplicities up to order 2, i.e., we count
$x$ as a double zero if $m(x)=m'(x)=0$. From simple zero counting,
we obtain the following result, that is essentially in Garkavi
[1959].

\proclaim Proposition 8.2.  Let $r$ be a fixed positive integer. Let
$M$ be an $n$-dimensional subspace of $C^{r}[a,b]$. Assume $M$
contains a $T$-system of dimension $m$ and is contained in an
$ET_2$-system of dimension $\ell$, with $m\ge \ell/2$. Then to each
$f\in C^{r}[a,b]$, there exists a unique best approximant from $M$.

\medskip
Based on the above Theorem 8.1, it is natural to ask whether strong
uniqueness always holds in this setting. We first give an example to
show that classical strong uniqueness need not hold. We then show
that a non-classical strong uniqueness with $\phi(t)=t^2$  holds in
this case.

\medskip\noindent {\bf Example.} Let $M= \span \{1, x, x^3\}$ in
$C[-1,1]$. Since $M$ contains a $T$-system of dimension 2 ($\span
\{1, x\}$), and is contained in an $ET_2$-system of dimension 4
($\span \{1, x, x^2, x^3\}$), it follows from Proposition 8.2 that
$M$ is necessarily a unicity space with respect to $C^{1}[-1,1]$.
Set $f(x) := 2x^2 - 1$. As is easily checked, its best approximant
from $M$ is the zero function. Since $A_{f}=\{-1, 0, 1\}$, we have
$$\gamma(f) =\min_{{m\in M}\atop {\|m\|=1}} \max \{ m(-1), - m(0),
m(1)\}.$$ Choosing $m(x)= x-x^3\in M$, we see that $\gamma(f) =0$,
and therefore classical strong uniqueness does not hold.

\medskip
Let us consider the above example in further detail. Set $m_a(x) :=
a(x-x^3)$. For $a$ sufficiently small, we have
$$\|f-m_a\| \ge \|f\| + \gamma \|m_a\|^p$$
for some $\gamma> 0$ if and only if $p\ge 2$. To see this, note that
$\|f\|=1$ and $\|m_a\| = 2|a|/3\sqrt{3}$. Thus the right-hand side
of the above is of the form
$$ 1 + \gamma \left( {{2|a|}\over {3 \sqrt{3}}}\right)^p.$$
Set
$$g(a) := \|f-m_a\| = \max_{x\in [-1,1]} |(2x^2-1) - a(x-x^3)|.$$
For $|a|$ sufficiently small, we have
$$g(a) = {1\over 3} + {2\over 9} \sqrt{4+3a^2} + {8\over {27 a^2}}
\left( \sqrt{4+3a^2} -2\right).$$ Expanding $\sqrt{4+3a^2}$ in a
Taylor series about $a=0$, we obtain
$$g(a) = 1 + {1 \over 8} a^2 - c a^4 + \cdots $$
for some $c>0$. Thus there exists a $\gamma>0$ such that
$$g(a) \ge 1 +\gamma \left( {{ 2|a|}\over {3\sqrt{3}}}\right)^p$$
for $a$ sufficiently small if and only if $p\ge 2$.

\medskip
Let $M\subset C^{r}[a,b]$ be an $n$-dimensional unicity space with
respect to $C^{r}[a,b]$. The previous example showed that we cannot
hope to obtain a classical strong uniqueness theorem. We shall,
however, prove the following result due to Kro\'o [1983a].

\proclaim Theorem 8.3. Let $r$ be a positive integer, $r\ge 2$. Let
$M\subset C^{r}[a,b]$ be a finite-dimensional unicity space with
respect to $C^{r}[a,b]$, and let $f\in C^{r}[a,b]$. Given any
positive constant $\sigma$, there exists a positive constant
$\gamma$, dependent on $f$, $M$ and $\sigma$ (but independent of
specific $m\in M$) such that if $m^*$ is the best approximant to $f$
from $M$, and if $m\in M$ satisfies
$$\|f-m\| - \|f-m^*\| \le \sigma,\eqno{(8.1)}$$
then
$$\|f-m\| - \|f-m^*\| \ge \gamma \|m-m^*\|^2.\eqno{(8.2)}$$

\medskip\noindent {\bf Remark.} Theorem 8.3 is a simpler case of a
general result in Kro\'o [1983a]. This more general result also
covers the case where $M\subset C^{1}[a,b]$ is an $n$-dimensional
unicity space with respect to $C^{1}[a,b]$. The proof thereof is
somewhat more cumbersome as the modulus of continuity of $f'$ and
$m'$ enter both the statement of the result and the analysis, i.e.,
a non-classical strong uniqueness result is proven with a more
complicated $\phi$. We choose to deal with the simpler case as
presented here.

\medskip
Prior to proving this theorem, let us prove an ancillary result
concerning local strong uniqueness.

\proclaim Proposition 8.4. Let $\phi$ be any nonnegative increasing
function defined on $[0,\infty)$. Assume there exists a $\sig_1>0$
such that if $m\in M$ satisfies
$$\|f-m\| - \|f-m^*\| \le \sig_1$$
then
$$\|f-m\| - \|f-m^*\| \ge \gamma_1 \phi(\|m-m^*\|)$$
for some $\gamma_1>0$. Then given any $\sig_2>0$, there exists a
$\gamma_2>0$, dependent on $\gamma_1$, $\sig_1$, $\sig_2$, $\phi$
and $\|f-m^*\|$, such that, for every $m\in M$ satisfying
$$\|f-m\| - \|f-m^*\| \le \sig_2,$$
we have
$$\|f-m\| - \|f-m^*\| \ge \gamma_2 \phi(\|m-m^*\|).$$

\pf We need only consider the case where $\sig_2 > \sig_1$. For all
$m\in M$ satisfying
$$\|f-m\| - \|f-m^*\| \le \sig_2,$$
we have
$$\|m-m^*\| \le \|f-m\| + \|f-m^*\| \le 2 \|f-m^*\| + \sig_2.$$
Thus for all such $m$, we have
$$\phi(\|m-m^*\|) \le \phi( 2 \|f-m^*\| + \sig_2).$$

Choose $\gamma^*>0$ for which
$$\gamma^* \phi( 2 \|f-m^*\| + \sig_2) \le \sig_1$$
and set
$$\gamma_2 =\min \{ \gamma_1, \gamma^*\}.$$

Then for $m\in M$ satisfying
$$\|f-m\| - \|f-m^*\| \le \sig_1,$$
we have
$$\|f-m\| - \|f-m^*\| \ge \gamma_1 \phi(\|m-m^*\|) \ge \gamma_2 \phi(\|m-m^*\|).$$
For $m\in M$ satisfying
$$\sig_1 \le \|f-m\| - \|f-m^*\| \le \sig_2,$$
we have
$$\eqalign{\gamma_2 \phi(\|m-m^*\|)&\;\le\; \gamma^* \phi(\|m-m^*\|)\le  \gamma^* \phi( 2 \|f-m^*\| + \sig_2)\cr
&\;\le\; \sig_1 \le \|f-m\| - \|f-m^*\|. \meop\cr} $$

\medskip\noindent
{\bf Proof of Theorem 8.3:} To prove the theorem, we show that if
$$\|f-m\| - \|f-m^*\| = \delta^2\eqno{(8.3)}$$
then
$$\sqrt{\gamma}\,\|m-m^*\| \le \delta.\eqno{(8.4)}$$
We shall choose a sufficiently small $\delta>0$ depending on $f$ and
$M$. This proves the theorem for $\sig=\delta^2$. We then apply
Proposition 8.4 to get the full result (with a different $\gamma$).

Recall from Theorem 2.3 that the function $m^*\in M$ is a best
approximant to $f$ from the $n$-dimensional subspace $M$ if and only
if there exist $k$ distinct points $x_1\nek x_k \in A_{f-m^*}$, for
$1\le k\le n+1$, and strictly positive values $\lam_1\nek \lam _k $
such that
$$\sum_{i=1}^k \lam_i [\sgn (f-m^*)(x_i)] m(x_i)= 0$$
for all $m\in M$.

Furthermore, from Theorem 8.1, the fact that $M$ is a unicity space
with respect to $C^{r}[a,b]$ implies that for the $\{x_i\}_{i=1}^k$
as above, if $m\in M$ satisfies $m(x_i) =0$, $i=1\nek k$, and
$m'(x_i) =0$ if $x_i\in (a,b)$, then $m$ is necessarily identically
zero. Thus the unicity space property says that the functional
$|||\,\cdot\, |||$ defined on $M$ by
$$|||m||| = \max_{i=1\nek k} |m(x_i)| + \max_{{{i=1\nek k}\atop
{x_i\in (a,b)}}} |m'(x_i)|$$ is a norm on $M$. Since $M$ is
finite-dimensional, this norm is equivalent to the uniform norm.
Thus, there exists a constant $c_1$ such that
$$ \|m-m^*\| \le c_1 ||| m-m^*|||$$
for all $m\in M$. Here and in what follows, the $c_j$ will always
denote positive constants that depend only upon $f$ and $M$. It
therefore suffices, in proving (8.4), to prove that, under
conditions (8.3), we have
$$|(m-m^*)(x_i)| \le c_2 \delta,\qquad i=1\nek k,\eqno{(8.5)}$$
and
$$|(m-m^*)'(x_i)|\le  c_3 \delta,\qquad {\rm if\ } x_i\in (a,b),
\quad i=1\nek k. \eqno{(8.6)}$$

We assume that (8.3) holds and will prove (8.5) and (8.6). We prove
(8.5) and (8.6) in two distinct lemmas. In the first lemma, we prove
more than what is formally needed to verify (8.5), but apply this
better estimate in the second lemma to prove (8.6).

\proclaim Lemma 8.5. Under the above assumptions, (8.5) holds. In
fact, $|(m-m^*)(x_i)| \le c_4 \delta^2$, $i=1\nek k$.

\pf We assume $m \in M$. Since $x_i\in A_{f-m^*}$, we have
$$\eqalign{\|f-m\|^2 &\;\ge\; |(f-m)(x_i)|^2 = |(f-m^*)(x_i) - (m-m^*)(x_i)|^2\cr
&\;=\; \|f-m^*\|^2 -2\|f-m^*\|\sgn [(f-m^*)(x_i)] (m-m^*)(x_i) +
|(m-m^*)(x_i)|^2.\cr}$$ From (8.3), we have
$$\|f-m\| - \|f-m^*\| = \delta^2,$$
and therefore, assuming $\delta< 1$,
$$\|f-m\|^2 = (\|f-m^*\| + \delta^2)^2 \le \|f-m^*\|^2 + c_5
\delta^2.$$ Thus
$$-2\|f-m^*\|\sgn [(f-m^*)(x_i)] (m-m^*)(x_i) + |(m-m^*)(x_i)|^2 \le
c_5\delta^2$$ implying that
$$\sgn [(f-m^*)(x_i)] (m-m^*)(x_i)\ge
-c_6\delta^2,\qquad i=1\nek k.$$ Recall that
$$\sum_{i=1}^k \lam_i [\sgn (f-m^*)(x_i)] m(x_i)= 0$$
for all $m\in M$, and normalize the $\lam_i$ so that $\sum_{i=1}^k
\lam_i=1$. The $\lam_i$ are independent of any specific $m$. Thus
$$ \eqalign{-c_6\delta^2&\;\le\; \sgn [(f-m^*)(x_i)] (m-m^*)(x_i) \;=\; -{1\over
{\lam_i}} \sum_{{j=1}\atop {j\ne i}}^k \lam_j [\sgn (f-m^*)(x_j)]
(m-m^*)(x_j)\cr &\;=\; {1\over {\lam_i}} \sum_{{j=1}\atop {j\ne
i}}^k \lam_j\left( -[\sgn (f-m^*)(x_j)] (m-m^*)(x_j)\right)\;\le\;
c_7\delta^2.\cr}$$ We have shown that
$$-c_6\delta^2 \le [\sgn (f-m^*)(x_i)] (m-m^*)(x_i)\le c_7\delta^2,$$
implying that
$$ |(m-m^*)(x_i)|\le c_4\delta^2,\qquad i=1\nek k.\ \ \meop$$

\proclaim Lemma 8.6. Under the above assumptions, (8.6) holds, i.e.,
$|(m-m^*)'(x_i)| \le c_3 \delta$, if $x_i\in (a,b)$, $i=1\nek k$.

\pf Assume $x_i\in (a,b)$. From Lemma 8.5,
$$|(m-m^*)(x_i)| \le c_4\delta^2,$$
where
$$\|f-m\| - \|f-m^*\| =\delta^2.$$

Now for $x_i\in (a,b)$ and $\delta>0$ such that $x_i\pm \delta \in
[a,b]$,
$$\eqalign{[\sgn (f-m^*)(x_i) ] (m-m^*)(x_i  \pm \delta)
&\;=\;  [\sgn (f-m^*)(x_i) ] \left( (f-m^*)(x_i \pm \delta) -
(f-m)(x_i \pm \delta)\right)\cr &\;\ge\;  [\sgn (f-m^*)(x_i)
](f-m^*)(x_i \pm \delta) - \|f-m\|\cr &\;=\; [\sgn (f-m^*)(x_i)
](f-m^*)(x_i \pm \delta) - \|f-m^*\| -\delta^2\cr &\;=\; [\sgn
(f-m^*)(x_i) ] \left( (f-m^*)(x_i \pm \delta) - (f-m^*)(x_i)\right)
- \delta^2.\cr}$$ Set
$$h(y) := (f-m^*)(x_i +y ) - (f-m^*)(x_i).$$
As $x_i\in A_{f-m^*}\cap (a,b)$, we have $(f-m^*)'(x_i) =0$. Thus we
see that $h(0) = h'(0) =0$. Since $h\in C^{r}$, $r\ge 2$, in a
neighborhood of the origin, it therefore follows that $|h(y)| \le
c_8 y^2$ for some constant $c_8$. Thus
$$ [\sgn (f-m^*)(x_i) ] (m-m^*)(x_i \pm \delta)  \ge -c_8\delta^2.\eqno{(8.7)}$$

In addition, from Taylor's Theorem,
$$[\sgn (f-m^*)(x_i) ] (m-m^*)(x_i\pm \delta)\le [\sgn (f-m^*)(x_i) ]
(m-m^*)(x_i)$$  $$\pm [\sgn (f-m^*)(x_i) ] (m-m^*)'(x_i)\,\delta +
c_9\delta^2.$$ Thus, as $|(m-m^*)(x_i)|\le c_4 \delta^2$ and using
(8.7), we have
$$\eqalign{ \pm \delta &[\sgn (f-m^*)(x_i) ] (m-m^*)'(x_i)\cr
&\;\ge\; [\sgn (f-m^*)(x_i) ] (m-m^*)(x_i\pm \delta)- [\sgn
(f-m^*)(x_i) ] (m-m^*)(x_i) - c_9 \delta^2\cr &\;\ge\; [\sgn
(f-m^*)(x_i) ] (m-m^*)(x_i\pm \delta) - c_{10}\delta^2\ge
-c_{11}\delta^2.\cr}$$ We have proven that
$$\pm \delta\, [\sgn (f-m^*)(x_i) ] (m-m^*)'(x_i)\ge -c_{11}\delta^2,$$
implying
$$|(m-m^*)'(x_i)|\le c_{11}\delta.\meop$$

\medskip\noindent
{\bf Proof of Theorem 8.3 (cont'd):} Lemmas 8.5 and 8.6 together
prove Theorem 8.3. \eop

The example prior to Theorem 8.3 shows that, in general, we cannot
expect to obtain a power less than 2 if $M$ is not a Haar space. In
fact, a result is proved in Kro\'o [1983a] which shows that this
estimate is optimal for any such subspace. We state this here
without proof.

\proclaim Theorem 8.7. Let $r$ be a positive integer, $r\ge 2$. Let
$M\subset C^{r}[a,b]$ be a $n$-dimensional unicity space with
respect to $C^{r}[a,b]$, and assume $M$ is not a Haar space, i.e.,
there exists an $\tilm\in M$, $\tilm\ne 0$, with $n$ distinct zeros
in $[a,b]$. Then there exists an $f\in C^\infty [a,b]$ with best
approximant $m^*\in M$, and a sequence of elements $m_k\in M$,
$m_k\to m^*$, such that
$$\|f-m_k\|-\|f-m^*\| \le \gamma \|m_k-m^*\|^2 $$
for some constant $\gamma >0$ that depends only on $f$ and $M$ and
is independent of $k$.

\goodbreak

\subsect{{9.} The $L^1$-Norm Revisited}

\medskip\noindent
In this section, we consider the space $L^1[a,b]$ with measure
$\mu$. We assume in what follows that $\mu$ is a non-atomic positive
finite measure on $[a,b]$ with the property that every real-valued
continuous function is $\mu$-measurable, and such that if $f\in
C[a,b]$ satisfies $\|f\|_1=0$ then $f=0$, i.e., $\|\cdot\|_1$ is
truly a norm on $C[a,b]$. We consider finite-dimensional subspaces
$M\subset C[a,b]$ that are unicity spaces with respect to all $f\in
C[a,b]$ in the $L^1[a,b]$ norm, i.e., for which there exists a
unique best approximant from $M$ to each $f\in C[a,b]$ in the
$L^1[a,b]$ norm. We restrict our consideration to approximation to
continuous functions because if we consider approximation to all
functions in $L^1[a,b]$ then there are no finite-dimensional unicity
spaces. But there are many unicity spaces in this continuous
setting, cf.\ Pinkus [1989]. One characterization of these unicity
spaces is the following which may be found in Cheney, Wulbert [1969,
Theorem 24].

\proclaim Theorem 9.1. Let $\mu$ be a non-atomic positive measure,
as above. Let $M$ be a finite-dimensional subspace of $C[a,b]$. Then
$M$ is a unicity space for all $f\in C[a,b]$ in the $L^1[a,b]$ norm
if and only if there does not exist an $h\in L^\infty[a,b]$ and an
$\tilm\in M$, $\tilm\ne 0$, for which
$$\eqalign{ ({\rm i}) & \, |h(x)|=1,\ {\sl all}\ x\in [a,b],\cr
({\rm ii}) & \int_a^b hm\dd\mu =0, \ {\sl all}\  m\in M,\cr ({\rm
iii}) &  \, h|\tilm|\ {\sl is \ continuous}.\cr}$$

Let $\ome(f;\del)$ denote the standard modulus of continuity of
$f\in C[a,b]$, i.e.,
$$\ome(f;\del) = \max\{ |f(x)-f(y)|: x,y\in [a,b],\, |x-y|\le \del\}.$$
We recall that $\ome(f;\cdot)$ is a continuous nonnegative
non-decreasing function for $\del \ge 0$ with $\ome(f; 0)=0$. For
any such continuous nonnegative non-decreasing function $F$, let
$F^{-1}(x) = \min\{ y: F(y)=x\}$ denote its inverse. Note that
$F^{-1}$ is nonnegative and increasing, but may have jumps, i.e., it
need not be continuous. This next result from Kro\'o [1981a]
provides a non-classical strong uniqueness result in the case where
$M$ is a finite-dimensional unicity space in $C[a,b]$ in the
$L^1[a,b]$ norm. It is well-known, for example, that Haar spaces
satisfy this condition. For Haar spaces, this next result was
obtained by Bj\"{o}rnest{\aa}l [1975].

\proclaim Theorem 9.2. Let $M$ be a finite-dimensional unicity space
in $C[a,b]$ in the $L^1[a,b]$ norm. If $f\in C[a,b]$, $m^*$ is the
best approximant to $f$ from $M$, and $\sig>0$ is an arbitrary
positive constant, then for any $m\in M$ satisfying
$\|f-m\|_1-\|f-m^*\|_1\le \sig$, we have
$$\|f-m\|_1-\|f-m^*\|_1 \ge \gamma \|m-m^*\|_1 \ome^{-1}(f; D\|m-m^*\|_1)$$
where $\gamma, D>0$ are constants that depend only on $f$, $\sig$
and $M$.

\pf For convenience, we assume that $m^* = 0$, i.e., replace $f-m^*$
by $f$. Set $\eps := \|f-m\|_1 - \|f\|_1$. We shall prove that if
$\eps \le 1$ then
$$\eps \ge \gamma \|m\|_1 \ome^{-1}(f; D \|m\|_1)$$
for some constants $\gamma, D>0$ that depend only on $f$ and $M$.
From Proposition 8.4, this suffices to prove the theorem.

As $0\in P_M(f)$, there exists, by Theorem 5.5, an $h\in L^\infty
[a,b]$ satisfying
$$\eqalign{ ({\rm i}) & \,|h(x)|=1,\ {\rm all}\ x\in [a,b],\cr
({\rm ii}) & \int_a^b hm\dd\mu =0, \ {\rm all}\  m\in M,\cr ({\rm
iii}) & \int_a^b hf \dd\mu = \|f\|_1.\cr}$$ From (i) and (iii), it
follows that $h(x) = \sgn f(x)$ $\mu$-a.e.\ on $N(f)$.

For each $x\in [a,b]$, we define the $\del$-neighborhood of $x$ as
$$O_\del := (x-\del, x+\del) \cap [a,b].$$
Given the above $h\in L^\infty[a,b]$, let $S(f)$ denote the set of
$x\in [a,b]$ for which
$$\mu(O_\del(x) \cap \{x: h(x)=1\})>0$$
and
$$\mu(O_\del(x) \cap \{x: h(x)=-1\})>0$$
for all $\del>0$. That is, $S(f)$ is the set of essential sign
changes of $h$. It is easily verified that if $S(f)$ is not empty
then it is closed and hence compact. Thus for each fixed $\del
>0$, we have
$$\min_{\sig=\pm 1} \inf_{x\in S(f)} \mu(O_\del(x) \cap \{x: h(x)=\sig\}) = \beta_f(\del)>0. \eqno{(9.1)} $$

We have yet to show that $S(f)\ne \emptyset$. This follows
immediately from Theorem 9.1. In fact, we have that
$$|||m||| := \max_{x\in S(f)} |m(x)|$$
is a norm on $M$. We prove this as follows. By the equivalence of
norms on finite-dimensional subspaces of $C[a,b]$ it suffices to
prove that only the zero function in $M$ vanishes identically on
$S(f)$. Assume to the contrary that there exists an $\tilm\in M$,
$\tilm \ne 0$, that vanishes on $S(f)$. We claim that we may choose
$h$ as above so that $h|\tilm|\in C[a,b]$. To see this, assume
$\tilm$ does not vanish on an interval $(c,d)$. As $(c,d) \cap
S(f)=\emptyset$ for each $x\in (c,d)$, there exists a $\del_x>0$
such that $h=1$ or $h=-1$ $\mu$-a.e.\ on $O_{\del_x}(x)$. From the
connectedness of $(c,d)$, this implies that we may take $h=1$ or
$h=-1$ on all of $(c,d)$. Thus $h|\tilm|$ is continuous on $(c,d)$,
and hence it easily follows that $h|\tilm|\in C[a,b]$.  But, by
Theorem 9.1, this contradicts the fact that $M$ is a unicity space
for all $f\in C[a,b]$ in the $L^1[a,b]$ norm. No such $\tilm$ exists
and therefore $|||\cdot|||$ is a norm on $M$.

As $|||\cdot|||$ is a norm on $M$, it therefore follows that there
exists a constant $c_1>0$, that depends only on $M$, satisfying
$$|||m||| \ge c_1 \|m\|_\infty$$
for all $m\in M$. Here and in what follows, the $c_j$ will always
denote positive constants that depend only upon $f$ and $M$.

Returning to the proof of the theorem, we have
$$\eqalign{\eps &\;=\; \|f-m\|_1 - \|f\|_1 \cr
&\;=\; \int_a^b (f-m) \sgn (f-m)\dd\mu - \int_a^b f \sgn f\dd\mu\cr
&\;=\; \int_{N(f)} (f-m) \sgn (f-m) \dd\mu + \int_{Z(f)} (-m) \sgn
(-m) \dd\mu - \int_{N(f)} f \sgn f \dd\mu\cr &\;=\; \int_{N(f)}
(f-m) \left[\sgn (f-m) -\sgn f\right] \dd\mu - \int_{N(f)} m \sgn
f\dd\mu + \int_{Z(f)} |m| \dd\mu .\cr}$$ From (ii), we have
$$\int_a^b h m  \dd\mu =0$$
for all $m\in M$, while $h=\sgn f$ on $N(f)$. Thus
$$\eqalign{\eps &\;=\; \int_{N(f)} (f-m) \left[\sgn (f-m) -\sgn f\right]
\dd\mu + \int_{Z(f)} m h\dd\mu + \int_{Z(f)} |m|\dd\mu\cr &\;=\; 2
\int_{\{x: f(f-m)<0\}} |f-m|\dd\mu + \int_{Z(f)} |m| + m
h\dd\mu.\cr}\eqno{(9.2)}$$

For $m\in M$, as above, let $x^*\in S(f)$ satisfy
$$|m(x^*)|= |||m||| \ge c_1 \|m\|_\infty$$
and for convenience set $\sig^* :=\sgn m(x^*)$. Thus, there exists a
$\del>0$ such that for all $x\in O_\del(x^*)$,
$$\sig^* m(x) = |m(x)| \ge {{c_1}\over 2} \|m\|_\infty \eqno{(9.3)}$$
and thus by (9.1)
$$\mu(O_\del(x^*) \cap \{x: h(x)=\sig^*\}) \ge
\beta_f(\del)=c_2>0.\eqno{(9.4)}$$ Let
$$B := O_\del(x^*) \cap \{x: x\in N(f),\, \sgn f(x)=\sig^*\}$$
and
$$ Q := O_\del(x^*) \cap \{x: x\in Z(f),\, h(x)=\sig^*\}.$$
By (9.4), we have
$$\mu(B) + \mu(Q)\ge c_2.$$
We now consider three options.

\smallskip\noindent
{\bf I.} Assume $\mu(Q) \ge c_2/2$. Since $Q\subseteq Z(f)$, we have
from (9.2) that
$$\eqalign{\eps &\;\ge\; \int_{Z(f)} |m| + mh\dd\mu \ge \int_Q |m|+mh\dd\mu
\;=\; 2\int_Q |m|\dd\mu\cr &\;\ge\; c_1\|m\|_\infty \mu(Q) \ge c_3
\|m\|_1\cr}\eqno{(9.5)}$$ due to the equivalence of the
$\|\cdot\|_1$ and $\|\cdot\|_\infty$ norms on $M$. Now if $m\in M$
satisfies $\|f-m\|_1-\|f\|_1 \le 1$ then
$$\|m\|_1 \le \|f-m\|_1 +\|f\|_1 \le 2\|f\|_1 +1,$$
i.e., all such $m$ are uniformly bounded in norm, and therefore
$\ome^{-1}(f; \|m\|_1)$ is also uniformly bounded above for all such
$m$. Thus (9.5) implies the desired result.

\smallskip\noindent
{\bf II.} Assume $\mu(B) \ge c_2/2$ and $\sig^* f(x) < (c_1/4)
\|m\|_\infty$ for all $x\in O_\del(x)$. From (9.3) and our
assumption, we have that
$$B \subseteq \{x: f(x)(f-m)(x)<0\}$$
and therefore from (9.2) and (9.3)
$$\eqalign{\eps &\;\ge\; 2\int_{\{x: f(f-m)<0\}} |f-m| \dd\mu \;\ge\; 2\int_B
|f-m|\dd\mu\cr &\;\ge\; 2\int_B {{c_1}\over 2} \|m\|_\infty -
{{c_1}\over 4} \|m\|_\infty \dd\mu \;=\; {{c_1}\over 2} \|m\|_\infty
\mu(B) \;\ge\; c_4 \|m\|_1\cr}$$ which, as above, implies the
result.

\smallskip\noindent {\bf III.} Assume $\mu(B) \ge c_2/2$ and there exists an
$\tilx\in O_\del(x^*)$ for which  $\sig^* f(\tilx) \ge (c_1/4)
\|m\|_\infty$. Since $x^*\in Z(f)$, we have $f(x^*)=0$. Assume,
without loss of generality, that $\tilx > x^*$. Choose $x^* \le
x_1<x_2<\tilx$ so that $f(x_1)=0$, $\sig^* f(x_2)= (c_1/4)
\|m\|_\infty$ and $0<\sig^* f(x) < (c_1/4) \|m\|_\infty$ on $(x_1,
x_2)$. As $f(x_2) - f(x_1) = \sig^* (c_1/4) \|m\|_\infty$, we have
$$x_2-x_1 \ge \ome_f^{-1}\left({{c_1}\over 4}
\|m\|_\infty\right).\eqno{(9.6)}$$ Since $(x_1, x_2)\subseteq O_\del
(x^*)$ from (9.3), we see that
$$f(x) (f-m)(x)<0$$
for $x\in (x_1, x_2)$. Thus from (9.2) and (9.6)
$$\eqalign{\eps &\;\ge\; 2\int_{x_1}^{x_2} |f-m| \dd\mu \;\ge\; 2\int_{x_1}^{x_2} {{c_1}\over 2} \|m\|_\infty -
{{c_1}\over 4} \|m\|_\infty \dd\mu\cr&\;=\;  {{c_1}\over 2}
\|m\|_\infty (x_2-x_1) \;\ge\; {{c_1}\over 2} \|m\|_\infty
\ome_f^{-1}\left({{c_1}\over 4} \|m\|_\infty\right) \;\ge\; \gamma
\|m\|_1 \ome_f^{-1}\left(D \|m\|_1\right).\cr}$$ This proves the
theorem. \eop

\medskip
We now turn to the one-sided $L^1$-approximation problem. For this
problem, unicity spaces are rare when considering approximating all
$f\in C[a,b]$ from some finite-dimensional subspace $M$ in $C[a,b]$;
see Pinkus [1989]. The situation is different if we consider the
case where both the approximating subspace $M$ and the functions to
be approximated are in $C^1[a,b]$, i.e., are continuously
differentiable. In this case, there are many unicity spaces.

For convenience, we again consider the interval $[a,b]$. The more
general case can be found in Pinkus [1989].  For each $f\in
C^1[a,b]$, we define $Z_1(f)$ to be the set of zeros of $f$ in
$[a,b]$ that only includes interior zeros if the function and its
derivative vanish, i.e., $x\in (a,b)$ is in $Z_1(f)$ only if
$f(x)=f'(x)=0$. Note that if $f\in C^1[a,b]$ satisfies $f\ge 0$,
then $Z_1(f)=Z(f)$. In what follows, we also always assume that $M$
contains a strictly positive function. With these assumptions, we
have the following characterization of unicity subspaces that is
from Pinkus, Strauss [1987]; see also Strauss [1982] and Pinkus
[1989].

\proclaim Theorem 9.3. An $n$-dimensional subspace $M$ of $C^1[a,b]$
is a unicity space for $C^1[a,b]$ in the one-sided $L^1$-norm
problem if and only if there does not exist an $m^*\in M$, $m^*\ne
0$, points $\{x_i\}_{i=1}^k$ in $Z_1(m^*)$, $1\le k\le n$, and
positive values $\{\lam_i\}_{i=1}^k$ satisfying
$$\int_a^b m  \dd\mu = \sum_{i=1}^k \lam_i m(x_i)$$
for all $m\in M$.

Let $m_1\nek m_n$ be any basis for the $n$-dimensional subspace $M$.
In what follows, we assume $M\subset C^1[a,b]$. Given an $f\in
C^1[a,b]$, we set
$$H_f(\del) := \ome(f';\del) + \sum_{i=1}^n \ome(m_i'; \del).$$
With these assumptions and definitions, we can now state the next
result which is in Kro\'o, Sommer, Strauss [1989]. Recall that
$M(f):= \{ m: m\in M, m\le f\}.$

\proclaim Theorem 9.4. Assume $M\subset C^1[a,b]$ is a
finite-dimensional unicity subspace for $C^1[a,b]$ in the one-sided
$L^1$-norm problem. Assume $f\in C^1[a,b]$, $m^*$ is the best
approximant to $f$ from $M(f)$, and $\sig>0$ is an arbitrary
positive constant. If $m\in M(f)$ is such that
$\|f-m\|_1-\|f-m^*\|_1 \le \sig$, then we have
$$\|f-m\|_1-\|f-m^*\|_1 \;\ge\; \gamma \|m-m^*\|_1 H_f^{-1}(D\|m-m^*\|_1)$$
for some constants $\gamma, D>0$ that depend only on $f$, $\sig$ and
$M$.

\smallskip
If, for example, $f$ and the elements of $M$ are such that their
derivatives are all in Lip $\alp$, $0<\alp \le 1$, then $H_f(\del) =
O(\del^\alp)$, and from the above, we obtain
$$\|f-m\|_1-\|f-m^*\|_1 \;\ge\; \gamma \|m-m^*\|^{(\alp+1)/\alp}_1$$
for some other constant $\gamma >0$ that depends only on $f$ and
$M$.

\smallskip
\pf For convenience, we assume that $m^*=0$, i.e., replace $f-m^*$
by $f$. Set
$$\eps := \|f-m\|_1 - \|f\|_1,$$
and assume $0<\eps \le 1$, as in the proof of Theorem 9.2. Since the
zero function is a best one-sided $L^1$-approximation from the
$n$-dimensional subspace $M$, we have that $f\ge 0$, and from
Theorem 5.7, there exist distinct points $\{x_i\}_{i=1}^k$ in
$Z_1(f)$, $1\le k \le n$, and positive numbers $\{\lam_i\}_{i=1}^k$
for which
$$\int_a^b m \dd\mu = \sum_{i=1}^k \lam_i m(x_i)$$
for all $m\in M$.

Now
$$\eps = \|f-m\|_1 - \|f\|_1 = \int_a^b (f-m)\, d\mu - \int_a^b f\dd\mu =
-\int_a^b m\,\dd \mu = - \sum_{i=1}^n \lam_i m(x_i).$$ As $f(x_i)=0$
and $f-m\ge 0$, we have $m(x_i)\le 0$ and thus
$$\eps = \sum_{i=1}^k \lam_i |m(x_i)|.$$
Recall that our aim is to prove that there exist constants $\gamma,
D>0$ that depend only on $f$ and $M$ such that
$$\eps = \sum_{i=1}^k \lam_i |m(x_i)| \ge \gamma \|m\|_1 H_f^{-1}(D\|m\|_1).$$

As $M$ is a unicity space for $C^1[a,b]$ in the one-sided $L^1$-norm
problem, we have from Theorem 9.3 that there does not exist an
$\tilm\in M$, $\tilm\ne 0$, points $\{y_i\}_{i=1}^k$ in
$Z_1(\tilm)$, $1\le k\le n$, and positive values $\{\mu_i\}_{i=1}^k$
satisfying
$$\int_a^b m \dd \mu = \sum_{i=1}^k \mu_i
m(y_i)$$ for all $m\in M$. Thus, for the above $\{x_i\}_{i=1}^k$, if
$m\in M$ satisfies $m(x_i)=0$, $i=1\nek k$, and $m'(x_i) =0$,
$x_i\in (a,b)$, then necessarily $m=0$. This implies that
$$|||m||| := \sum_{i=1}^k |m(x_i)| + \sum_{x_i\in (a,b)}|m'(x_i)|$$
is a norm on $M$ and therefore
$$\|m\|_1 \le c_1 |||m|||\eqno(9.7)$$
for all $m\in M$. As usual, the $c_j$ will always denote positive
constants that depend only upon $f$ and $M$.

As $\eps = \sum_{i=1}^n \lam_i |m(x_i)|$ where $\lam_i>0$ for all
$i$, we have
$$|m(x_i)| \le c_2 \eps,\qquad i=1\nek k.\eqno(9.8)$$
To use the norm $|||m|||$, we must also bound the $|m'(x_i)|$ for
$x_i\in (a,b)$. To this end, assume $m=\sum_{i=1}^n b_i m_i$. Then
$$\ome(m';\del) \le \sum_{i=1}^n |b_i|\, \ome(m'_i;\del) \le (\max_{j=1\nek n} |b_j|)
\sum_{i=1}^n \ome(m'_i;\del).$$ As $m_1\nek m_n$ is a basis for $M$,
it follows that $\max_{j=1\nek n} |b_j|$ is also a norm on $M$.
Therefore $\max_{j=1\nek n} |b_j|\le c_3 \|m\|_1$. Furthermore, $
\sum_{i=1}^n \ome(m'_i;\del) \le H_f(\del)$. Thus
$$\ome(m';\del) \le c_3\|m\|_1 H_f(\del).\eqno(9.9)$$
For $m$ satisfying  $\|f-m\|_1-\|f\|_1 \le 1$, it follows that
$$\|m\|_1 \le \|f-m\|_1 + \|f\|_1 \le 2\|f\|_1 +1.$$
Thus $\|m\|_1$ is uniformly bounded above and from (9.9), we have
$$\ome(m';\del) \le c_4 H_f(\del).\eqno(9.10)$$

Consider $x_i\in (a,b)$. Choose $1 \ge \eta >0$ such that the
$(x_i-\eta, x_i+\eta)$ are disjoint intervals of $(a,b)$. For each
$m\in M$, $0<h<\eta$ and $b=\sgn m'(x_i)$, we have
$$m(x_i +bh) = m(x_i) + bh\, m'(x_i) + E = m(x_i) + h |m'(x_i)| + E.\eqno(9.11)$$
An easy estimate and (9.10) shows that
$$|E| \le h \ome(m';h) \le c_4 h H_f(h).\eqno(9.12)$$
Since $f(x_i)=f'(x_i)=0$, we therefore have
$$|f(x_i + bh)| =  |f(x_i + bh)- f(x_i)| \le h\ome(f';h) \le h H_f(h).\eqno(9.13)$$
From $m(x_i + bh) \le f(x_i+bh)$, we have by (9.7), (9.11), (9.12)
and (9.13)
$$\eqalign{ |m'(x_i)| = & \,{{m(x_i+bh) - m(x_i) -E}\over h} \cr
\le & \, {1 \over h} \left(f(x_i+bh) + c_2 \eps - E\right)\cr \le &
\, {1 \over h} \left( h H_f(h) + c_2\eps + c_4 h H_f(h)\right)\cr =
& \, c_5 \left(H_f(h) + {\eps \over h}\right). \cr}\eqno(9.14)$$

Using the equivalence of norms, from (9.7), and (9.14), we finally
obtain
$$\eqalign{\|m\|_1 \le & \, c_1 \left( \sum_{i=1}^k |m(x_i)| + \sum_{x_i\in (a,b)}|m'(x_i)|\right)\cr
\le & \, c_1 \left(k c_2\eps + k c_5 \left(H_f(h) + {\eps \over
h}\right)\right) \cr \le & \, c_6 \left( H_f(h) + {\eps \over h}
\right)\cr}\eqno(9.15)$$ for all $h\in (0,\eta)$.

Let us recall that we wish to bound $\eps = \|f-m\|_1 - \|f\|_1$
from below. To this end, we consider two cases.

\smallskip \noindent
{\bf Case 1:} $\|m\|_1 \le 2 c_6 \eps/\eta$. In this case, it
follows from the monotonicity of $H_f$ that
$$H_f^{-1}(\|m\|_1) \le H_f^{-1}(2 c_6 \eps/\eta) \le c_7$$
and therefore
$$\|m\|_1 H_f^{-1}(\|m\|_1) \le 2 c_6 c_7\eps/\eta,$$
i.e.,
$$\gamma \|m\|_1 H_f^{-1}(\|m\|_1) \le \eps.$$

\smallskip \noindent
{\bf Case 2:} $\|m\|_1 > 2 c_6 \eps/\eta$. Set $h:= 2 c_6
\eps/\|m\|_1 < \eta$ in (9.15) to obtain
$$\|m\|_1 \le  c_6 \left( H_f\left({{2 c_6 \eps}\over {\|m\|_1}}\right)
+ {{ \|m\|_1}\over {2 c_6 }} \right).$$ Thus
$$\|m\|_1 \le {{ \|m\|_1}\over {2 }} + c_6 H_f\left({{2 c_6 \eps}\over {\|m\|_1}}\right)$$
implying
$$ c_8 \|m\|_1 \le  H_f\left({{2 c_6 \eps}\over {\|m\|_1}}\right)$$
and therefore
$$H_f^{-1}\left( c_8 \|m\|_1 \right) \le {{2 c_6 \eps}\over {\|m\|_1}}$$
whence
$$ {{\|m\|_1}\over {2 c_6}} H_f^{-1}\left( c_8 \|m\|_1 \right) \le  \eps.$$
This proves the theorem. \eop

A converse result is proved in Kro\'o, Sommer, Strauss [1989] that
shows that this estimate is optimal. We state it here without proof.

\proclaim Theorem 9.5. Assume $M\subset C^1[a,b]$ is an
$n$-dimensional unicity subspace for $C^1[a,b]$ in the one-sided
$L^1$-norm problem. Then there exists an $f\in C^1[a,b]$ with $m^*$
the best approximant to $f$ from $M(f)$, and a sequence of elements
$m_k\in M(f)$, $m_k\to m^*$, such that
$$\|f-m_k\|_1-\|f-m^*\|_1 \le \gamma \|m_k-m^*\|_1 H_f^{-1}(D\|m_k-m^*\|_1)$$
for some constants $\gamma , D>0$ that depend only on $f$ and $M$
and are independent of $k$.

\goodbreak

\subsect{{10.} Strong Uniqueness in Complex Approximation in the
Uniform Norm}

\medskip\noindent
In this section, we consider the question of strong uniqueness in
the space of continuous complex-valued functions endowed with the
uniform norm. The principal result here will indicate that classical
strong uniqueness fails, in general, in the complex setting in a
Haar space. Instead, we shall derive a non-classical strong
uniqueness type result of order 2, i.e., prove that for $m^*$ the
best approximant to $f$ from $M$, and any $m\in M$ sufficiently
close to $m^*$, we have
$$\|f-m\|-\|f-m^*\|\geq \gamma \|m-m^*\|^2,$$
where $\gamma >0$ depends only on $f$ and $M$.  To this end, we
shall need to extend some results from Section 2 to the complex
setting. Throughout this section, $C(K)$ will denote the space of
complex-valued functions, continuous on the compact Hausdorff space
$K$. First, we shall require a formula for the complex
$\tau$-functional that is similar to the one provided by Theorem
2.1. We recall that $A_f = \{ x : |f(x)|=\|f\|\}$, and if $z=x+\ii
y$ where $x,y\in \RR$, then we set $\Re z=x$ and $\Im z = y$.

\proclaim Theorem 10.1. For any functions $f,g\in C(K)$, $f\ne 0$,
we have
$$ \tau_+(f,g)={1 \over {\|f\|}}\max_{x\in A_f}\Re
(\overline{f(x)}g(x)).$$

\pf Note that for any $z_1, z_2 \in \CC$ and $t\in \RR$, we have
$$ |z_1+tz_2|^2= |z_1|^2  +2t\, \Re\, \overline{z_1}z_2+t^2|z_2|^2.  \eqno (10.1)$$

Using this relation, we obtain for every $x\in A_f$, $f\ne 0$,
$$\tau_+(f,g)\geq \lim_{t\rightarrow
0^+}{{|f(x)+tg(x)|^2-|f(x)|^2}\over {t(|f(x)+tg(x)|+|f(x)|)}}={{\Re
(\overline{f(x)}g(x))}\over {|f(x)|}},$$ i.e.,
$$ \tau_+(f,g)\geq {1\over {\|f\|}}\max_{x\in A_f}\Re (\overline{f(x)}g(x)).$$
This proves the desired lower bound for $\tau_+(f,g)$.

In order to verify the upper bound, choose $t_n\rightarrow 0^+$ and
let $x_n\in K$ be such that $\|f+t_ng\|$ is attained at $x_n$.
Similar to the argument used in the proof of Theorem 2.1, we can
assume, without loss of generality, that $x_n\rightarrow x^*\in
A_f$. Then using the fact that $|f(x_n)|^2\leq |f(x^*)|^2$, we
obtain
$$\eqalign{\tau_+(f,g)\;=\;\lim_
{n\rightarrow\infty} & {{|f(x_n)+t_ng(x_n)|^2-|f(x^*)|^2}\over
{t_n(|f(x_n)+t_ng(x_n)|+|f(x^*)|)}}\;=\;\lim_
{n\rightarrow\infty}{{|f(x_n)|^2-|f(x^*)|^2}\over
{t_n(|f(x_n)+t_ng(x_n)|+|f(x^*)|)}}\cr & + \lim_
{n\rightarrow\infty}{{2\Re (\overline{f(x_n)}g(x_n))}\over
{|f(x_n)+t_ng(x_n)|+|f(x^*)|}}\;\leq\; {{\Re
(\overline{f(x^*)}g(x^*))}\over {|f(x^*)|}}.\cr}$$ (Note that since
the limit on the left hand side and the second limit on the right
hand side exist, the first limit on the right hand side must exist,
as well.) Clearly this provides the needed upper bound. \eop

Applying the above and the argument of Theorem 1.4, we obtain

\proclaim Corollary 10.2. Let $f\in C(K)$ and $M$ be a linear
subspace of $C(K)$. Then $m^*\in M$ is a best approximant to $f$
from $M$ if and only if
$$\max_{x\in A_{f-m^*}}\Re (\overline{(f-m^*)(x)}m(x))\geq 0,\quad m\in M. $$
Moreover, $m^*$ is the strongly unique best approximant to $f$ from
$M$ if and only if $$\inf _{m\in M, \|m\|=1}\max_{x\in A_{f-m^*}}\Re
(\overline{(f-m^*)(x)}m(x))>0.$$

The next theorem provides an analogue of Theorem 2.3 in the complex
setting. Since this is a key result in this section, we include a
proof.

\proclaim Theorem 10.3. Let $f\in C(K)$ and $M$ be an
$n$-dimensional subspace of $C(K)$. Then $m^*\in M$ is a best
approximant to $f$ from $M$ if and only if there exist points
$x_1\nek x_k\in A_{f-m^*},$ and positive numbers $\lambda_1 \nek
\lambda_k, 1\leq k\leq 2n+1$, such that for every $m\in M$, we have
$$\sum_{j=1}^{k}\lambda_j\overline{(f-m^*)(x_j)}m(x_j)=0. \eqno (10.2)$$

\pf ($\Leftarrow$) We may assume that $\sum_{j=1}^{k}\lambda_j=1$
and $m^*=0$. Then using (10.2), we have for every $m\in M$
$$\|f\|^2=\sum_{j=1}^{k}\lambda_j\overline{f(x_j)}f(x_j)=\sum_{j=1}^{k}\lambda_j\overline{f(x_j)}(f-m)(x_j)\leq
\|f\| \cdot \|f-m\|.$$ Dividing both sides of the above inequality
by $\|f\|$ yields that $m^*=0$ is a best approximant to $f$ from
$M$.

\smallskip\noindent
($\Rightarrow$) Assume again that $m^*=0$. Let $M=\span\{m_1\nek
m_n\}$. For any $x\in A_{f}$, consider the vectors $$\bfu_x:=(\Re
(\overline{f(x)}m_1(x)), -\Im (\overline{f(x)}m_1(x))\nek
\Re(\overline{f(x)}m_n(x)), -\Im (\overline{f(x)}m_n(x)))
\in\RR^{2n}.$$ Furthermore, let $D:=\{\bfu_x, x\in
A_f\}\subset\RR^{2n}$. We claim that ${\bf 0}$ lies in the convex
hull of $D$. Indeed, if this is not the case, then by the separating
hyperplane theorem, there exists a $\bfc=(c_1\nek
c_{2n})\in\RR^{2n}$ such that $(\bfc, \bfu_x)<0$ for every $x\in
A_f$. This clearly yields that for some $\tilm \in M$, we have $\Re
(\overline{f(x)}\tilm(x)) <0$ for every $x\in A_f$. But in view of
Corollary 10.2, this contradicts the assumption that $0$ is a best
approximant to $f$ from $M$. Thus ${\bf 0}$ belongs to the convex
hull of $D\subset\RR^{2n}$ and, by the Carath\'eodory Theorem, ${\bf
0}$ is a convex linear combination of at most $2n+1$ points of $D$.
Evidently, this verifies relations (10.2). \eop

We are now in a position to prove a strong uniqueness type result
for complex Chebyshev approximation. This is the main result of this
section and can be found in Newman, Shapiro [1963]. The proof as
presented here is from Smarzewski [1989]. As previously, we shall
call the $n$-dimensional subspace $M$ of $C(K)$ a Haar space if
every nontrivial element of $M$ has at  most $n$ distinct zeros in
$K$.

\proclaim Theorem 10.4. Let $M$ be a finite-dimensional Haar space
in $C(K)$. Let $f\in C(K)$ and $m^*$ be the best approximant to $f$
from $M$. Then there exists a $c>0$ that depends only upon $f$ and
$M$ such that
$$\|f-m\|^2-\|f-m^*\|^2\ge c\|m-m^*\|^2, \eqno (10.3)$$
for all $m\in M$. In addition, if $m\in M$ satisfies
$$\|f-m\| - \|f-m^*\| \le \sig$$
then
$$\|f-m\|-\|f-m^*\|\ge \gamma \|m-m^*\|^2, \eqno (10.4)$$
where $\gamma  :=c/(2\|f-m^*\| +\sig)>0$ and thus also depends only
upon $f$, $M$ and $\sig$.

\pf Inequality (10.4) is an immediate consequence of (10.3) since if
$\|f-m\| -\|f-m^*\| \le \sig$, then
$$\eqalign{ c \|m-m^*\|^2 &\;\le\; \|f-m\|^2 - \|f-m^*\|^2 \;=\; (\|f-m\| + \|f-m^*\|) (\|f-m\| - \|f-m^*\|)\cr
&\;\le\;  (2\|f-m^*\| +\sig) (\|f-m\| - \|f-m^*\|).\cr}$$

It remains to prove (10.3). Assume $M$ is a Haar space of dimension
$n$. We recall from Theorem 10.3 that
$$\sum_{j=1}^{k}\lambda_j\overline{(f-m^*)(x_j)}m(x_j)=0$$
for all $m\in M$, where $x_1\nek x_k \in A_{f-m^*}$, $\lam_1\nek
\lam_k >0$ and, since $M$ is a Haar space, we have $n+1\le k \le
2n+1$ and, more importantly,
$$\left( \sum_{j=1}^k \lam_j |m(x_j)|^2\right)^{1/2} =: |||m|||$$
is a norm on $M$. Hence there exists a $c>0$ such that
$$|||m|||^2 \ge c \|m\|^2\eqno (10.5)$$
for all $m\in M$.

For each $m\in M$ and $j\in \{1\nek k\}$, we have
$$\eqalign{\|f-m\|^2 &\;\ge\; |(f-m)(x_j)|^2 \;=\; |(f-m^*)(x_j) +
(m-m^*)(x_j)|^2\cr &\;=\; \|f-m^*\|^2 + 2 \Re
\left(\overline{(f-m^*)(x_j)}(m-m^*)(x_j)\right) +
|(m-m^*)(x_j)|^2.\cr}$$ Multiply the above by $\lam_j>0$ and sum
over $j$. Assuming, without loss of  generality, that $\sum_{j=1}^k
\lam_j=1$, we obtain
$$\|f-m\|^2 \ge  \|f-m^*\|^2 + 2 \Re \left(\sum_{j=1}^k \lam_j \overline{(f-m^*)(x_j)}(m-m^*)(x_j)\right)$$
$$+ \sum_{j=1}^k \lam_j |(m-m^*)(x_j)|^2.$$
Applying (10.2) and (10.5), this gives
$$\|f-m\|^2-\|f-m^*\|^2\ge c\|m-m^*\|^2. \meop $$

Clearly, estimations (10.3) and (10.4) provide strong uniqueness
type results that are weaker than the classical strong uniqueness
result that holds in the real case when approximating by elements of
Haar spaces. This raises the natural question of whether classical
strong uniqueness can also hold for every complex function. Our next
proposition, see Gutknecht [1978] for a particular example thereof,
shows that this is not the case. In fact, it turns out that, in
general, (10.4) provides the best possible estimate.

\proclaim Proposition 10.5. Let $M$ be any subspace in $C(K)$
containing the constant function. Assume that there exists an $f\in
C(K)$ such that $\Im f\equiv 0$ on $K$ and $m^*\equiv 0$ is the best
approximant to $f$ from $M$. (Such an $f$ will clearly exist, for
instance, when $M$ possesses a basis consisting of real functions.)
Then for $m:= \ii b \in M$, $b\in \RR$, we have the following
converse to (10.4):
$$\|f-m\|-\|f-m^*\|\leq {{\|m-m^*\|^2}\over {2\|f\|}}.\eqno (10.6)$$

\pf Indeed, by the above assumptions, $\|f-m\|^2=\|f\|^2+\|m\|^2$
and since $m^*\equiv 0$, the above inequality easily follows. \eop

From (10.6), it follows that there cannot exist a $\gamma>0$ such
that
$$\|f-m\|-\|f-m^*\|\ge \gamma \|m-m^*\|$$
for all $m\in M$.

Despite the fact that classical strong uniqueness fails, in general,
in the complex case, it is possible to give a relatively simple
sufficient condition for it to hold. Let $f\in C(K)$ and assume that
$m^*$ is its best approximant. We shall say that the set of points
$A\subset A_{f-m^*}$ is {\sl extremal} if $m^*$ is a best
approximant to $f$ from $M$ also on the set $A$. Moreover, $A$ is a
minimal extremal set if no proper subset of $A$ is extremal.
Clearly, any set of points $x_1 \nek x_k\in A_{f-m^*}$, $k\in
\{1\nek 2n+1\}$, satisfying conditions (10.2) of Theorem 10.3 must
be extremal. Hence the cardinality of any minimal extremal set can
be at most $2n+1$. This also follows easily from Corollary 10.2.

We now give a sufficient condition for classical strong uniqueness
to hold in the complex case. This sufficient condition is to be
found in Theorem 2 in Gutknecht [1978].

\proclaim Theorem 10.6. Let $f\in C(K)$ and $M$ be an
$n$-dimensional subspace of $C(K)$. Assume that $m^*$ is a best
approximant to $f$ from $M$ possessing a minimal extremal set of
cardinality $2n+1$. Then $m^*$ is a strongly unique best approximant
of $f$, i.e., for all $m\in M$, we have
$$\|f-m\|-\|f-m^*\|\geq \gamma \|m-m^*\|,$$
for some $\gamma>0$.

\pf Let $\{x_1 \nek x_{2n+1}\}\subset A_{f-m^*}$ be a minimal
extremal set. Then by Theorem 10.3, there exist corresponding
positive numbers $\lambda_1 \nek \lambda_{2n+1}$ such that (10.2)
holds with $k=2n+1$. Assume that to the contrary $m^*$ is not a
strongly unique best approximant to $f$ from $M$. Then by Corollary
10.2, there exists an $m_0\in M$, $m_0 \ne 0$, such that
$$\max_{x\in A_{f-m^*}}\Re (\overline{(f-m^*)(x)}m_0(x))=0.$$
Since relations (10.2) hold for this $m_0$, we obtain that
$$\Re (\overline{(f-m^*)(x_j)}m_0(x_j))=0, \quad j=1 \nek 2n+1.  \eqno
(10.7)$$

Let $\{m_1 \nek m_n\}$ be a basis for $M$. For each $j\in \{1\nek
2n+1\}$, consider the vectors
$$\bfu_j:= (\Re (\overline{(f-m^*)(x_j)}m_i(x_j)), -\Im
(\overline{(f-m^*)(x_j)}m_i(x_j)), \ i=1\nek n)\in\RR^{2n},$$ and
set $B:=\{\bfu_j\}_{j=1}^{2n+1} \subset \RR^{2n}$. Then relations
(10.7) yield that the set $B$ lies in a ($2n-1$)-dimensional
hyperplane $H$ of $\RR^{2n}$. Moreover, by the standard argument
repeatedly used above, ${\bf 0}$ belongs to the convex hull of $B$
(otherwise relations (10.2) would fail for some $m\in M$). Now,
using the fact that the convex hull of $B$ is of dimension $2n-1$,
it follows from the Carath\'eodory Theorem that ${\bf 0}$ is a
convex linear combination of at most $2n$ points from $B$. But this
in turn means that (10.2) holds for a proper subset of $\{x_1 \nek
x_{2n+1}\}$ (with some $\lambda$'s). Clearly this set of at most
$2n$ points is extremal too, contradicting the condition of
minimality of the extremal set $\{x_1 \nek x_{2n+1}\}$. \ \eop

Thus if $m^*$ is a best approximant with minimal extremal set
consisting of exactly $2n+1$ points then the best approximation is
strongly unique in the classical sense. This result is not vacuous.
It is easy to construct examples that satisfy the criteria of
Theorem 10.6. Our next example, see Gutknecht [1978], shows that we
do not necessarily have classical strong uniqueness for a minimal
extremal set consisting of fewer than $2n+1$ points.

\medskip\noindent
{\bf Example.} Let $K:=\{|z|=1 \,:\, z\in\CC\}$ be the unit circle,
$M:=\Pi_{n-1}$ the space of polynomials of degree at most $n-1$, and
$$f(z):={{z^n+z^{3n}}\over {2}}.$$ Then it is easy to see that
$\|f\|=1$ and $A_f=\{z_k:=\ee^{\ii\pi k/n}: k=1 \nek 2n\},$ i.e.,
$A_f$ is the set of all roots of unity of order $2n$. In particular,
$f(z_k)=(-1)^k$, $k=1\nek 2n$. We claim that $m^*=0$ is the best
approximant of $f$. Indeed, if for some $m\in M$, we had
$\|f-m\|<\|f\|=1$ then evidently
$$(-1)^k\Re m(z_k)>0, \quad k=1\nek 2n, \eqno (10.8)$$ where
$m=\sum_{j=0}^{n-1}d_j z^j$, $d_j=a_j+\ii b_j$. Now $\Re m(\ee^{\ii
t})=\sum_{j=0}^{n-1}(a_j\cos jt-b_j\sin jt)$, i.e., $T_n:=\Re
m(\ee^{\ii t})$ is a trigonometric polynomial of degree at most
$n-1$. By (10.8), $T_n$ has at least $2n-1$ distinct zeros on the
unit circle. Thus, we must have $T_n=0$, and that contradicts
(10.8). Hence $m^*=0$ is a best approximant to $f$ from $M$. It is
also the unique best approximant since $M$ is a Haar space. Note
that by Corollary 10.2 this best approximant is not strongly unique
in the classical sense since choosing $m$ to be an imaginary
constant, we have
$$\max_{x\in A_{f-m^*}}\Re (\overline{(f-m^*)(x)}m(x))= \max_{x\in
A_{f}}\Im (\overline{f(x)})= 0.$$ It remains to show that the
extremal set $A_f=\{z_k:=\ee^{\ii\pi k/n}\,: \,k=1 \nek 2n\}$ is
minimal. Let us delete any point of $A_f$, say $z_1$, and assume
that the remaining set $\{z_k:=\ee^{\ii\pi k/n}: k=2 \nek 2n\}$ is
still extremal. (From symmetry considerations, deleting $z_1$ does
not restrict the generality.) Clearly, we can choose a trigonometric
polynomial $T^*=\Re m(\ee^{\ii t})$, $m\in  \Pi_{n-1}$, of degree
$n-1$ such that $(-1)^kT^*(k\pi/n)<0$, $k=2 \nek 2n $. But in view
of Corollary 10.2 this means that $m^*=0$ is no longer a best
approximant on this set of $2n-1$ points. Thus the extremal set
$A_f$ consisting of $2n$ points is minimal, but the best approximant
is not strongly unique in the classical sense.

\smallskip
Hence the condition of Theorem 10.6 is sharp, in a certain sense.
Nevertheless, as the next example verifies, this condition is not
necessary for classical strong uniqueness to hold.

\medskip\noindent
{\bf Example.} Let $K:=\{|z|=1: z\in\CC\}$ be the unit circle,
$M:=\Pi_{n-1}$ be the space of polynomials of degree at most $n-1$,
and $f(z):=z^n$. Then $A_f=K$ and as in the previous example the
function $m^*=0$ is the unique best approximant to $z^n$ from $M$
with the set $\{z_k:=\ee^{\ii\pi k/n}: k=1 \nek 2n\}$ being a
minimal extremal set of $2n$ points. Let us show that this best
approximant is strongly unique in the classical sense. Indeed, if
classical strong uniqueness is not valid then by Corollary 10.2 for
some nontrivial $m=\sum_{j=0}^{n-1}d_j z^j$, $d_j=a_j+ib_j$, we
would have
$$\max_{|z|=1}\Re (\overline{z^n} m(z))= 0.$$ On the other hand,
evidently $T_n:=\Re(\overline{z^n} m(z)) = \sum_{j=0}^{n-1} a_j \cos
(n-j)t + b_j \sin (n-j)t $ is a nonpositive trigonometric polynomial
of degree $n$ such that
$$\int_{0}^{2\pi}T_n(t)\dd t=0.$$
Thus $T_n$ must be identically zero implying that $m$ must be
identically zero. By this contradiction, it follows that the best
approximation is strongly unique in the classical sense. In fact,
the optimal strong uniqueness constant in this case is $1/n$; see
Rivlin [1984a].

\smallskip
An exact necessary and sufficient condition for when one has strong
uniqueness in the classical sense is to found in Blatt [1984], who
credits the result to Brosowski [1983]. These conditions are
somewhat technical and are not detailed here.

\sectwopn{{Part III.} Applications of Strong Uniqueness}

\noindent We shall briefly consider various applications of strong
uniqueness type results. The main idea behind these applications is
the following: instead of solving the best approximation problem in
a given norm, we replace it by considering another norm, close to
the original norm, one that leads to a simpler approximation
problem. Strong uniqueness is then applied in order to show that the
best approximant in this new norm is sufficiently close to the
original best approximant. Typically, the original norm is modified
by replacing it by a similar discrete norm, or by introducing a
weight function into the norm. We first start with some general
remarks concerning approximation in nearby norms, and then proceed
to a discussion of specific examples.

\subsect{{11.} Strong Uniqueness and Approximation in Nearby Norms}

\medskip\noindent
Let $X$ be a normed linear space with norm $\|\cdot \|$, and let $M$
a finite-dimensional subspace of $X$. Assume that the sequence of
seminorms $\|\cdot \|_k$, $k\in \NN$, approximate the given norm,
i.e., $\lim_{k\to\infty} \|f\|_k = \|f\|$, for all $f \in X$. We
also assume that for a given $f \in X$, its best approximant
$P_M(f)\in M$ with respect to the norm $\|\cdot \|$ is unique, and
denote by $P_M^{(k)}(f)\in M$ the set of best approximants to $f$ in
the seminorms $\|\cdot \|_k$, $k\in \NN$. We are interested in
estimating the deviation of $P_M^{(k)}(f)$ from $P_M(f)$. This
approach was first considered by Kripke [1964] and Peetre [1970].

To estimate the deviation of $P_M^{(k)}(f)$ from $P_M(f)$, we
introduce a quantity measuring the deviation of $\|\cdot \|_k$ from
$\|\cdot \|$ uniformly on the set $S(M):=\{m\in M: \|m\|=1\}$, the
unit sphere in $M$, namely
$$\eta_k(M):= \sup_{m\in S(M)}|\;\|m\|_k-1|.$$

\proclaim Lemma 11.1. For any finite-dimensional subspace $M\subset
X$, we have $$\lim_{k\to\infty} \eta_k(M)= 0.$$

\pf  Assume to the contrary that for a subsequence $k_i$, there
exist elements $u_i\in S(M)$ such that
$$|\;\|u_i\|_{k_i}-1| \geq \delta >0, \qquad i\in \NN. \eqno (11.1)$$
Let $u_i=\Sigma_{j=1}^na_j^i m_j$, where $m_1\nek m_n$ is a basis
for $M$. Since $S(M)\subset X$ is compact, we may assume, without
loss of generality, that $u_i$ converges to $u^*$ as $i\to\infty$
for some $u^*=\Sigma_{j=1}^na_j^*m_j\in S(M)$, i.e., $\|u^*\|=1$.
Moreover, by the equivalence of norms in finite-dimensional spaces,
$$\lim_{i\to\infty} \max_{1\leq j\leq n}|a_j^i-a_j^*|= 0.$$
In addition, $$\lim_{i\to\infty} \sum _{j=1}^n\|m_j\|_{k_i}=
\sum_{j=1}^n\|m_j\|.$$ Thus
$$ \|u^*-u_i\|_{k_i}\leq \max_{1\leq j\leq n}|a_j^i-a_j^*|\sum _{j=1}^n\|m_j\|_{k_i},
\eqno (11.2)$$ and the right-hand-side converges to $0$ as
$i\to\infty$.

Hence, we obtain by (11.2) that
$$ \lim_{i\to \infty}|1-\|u_i\|_{k_i}| \leq \lim_{i\to\infty}
(|1-\|u^*\|_{k_i}|+\|u^*-u_i\|_{k_i})=|1-\|u^*\|\;|=0.$$ But this
clearly contradicts (11.1). \eop

\proclaim Corollary 11.2. If $k\in \NN$ is sufficiently large so
that $\eta_k(M)<1/2$, then for all $m\in M$, we have
$${2\over 3}\|m\|_k\leq \|m\|\leq 2\|m\|_k. \eqno (11.3)$$

\smallskip
Assume now that for a given $f\in X$, with unique best approximant
$P_M(f)\in M$, non-classical strong uniqueness holds. That is, we
assume there exists a nonnegative strictly increasing function
$\phi$ (depending only on $f$ and $M$), defined on $[0,\sig]$, such
that for all $m\in M$ satisfying $\|f-m\|-\|f-P_M(f)\|\leq \sigma$,
we have
$$ \|f-m\|-\|f-P_M(f)\|\geq \phi (\|m-P_M(f)\|). \eqno (11.4)$$
Note that the above relation yields that $\phi(t)\leq t$ on
$[0,\sig]$. Moreover, if $\phi(t)\geq ct$ thereon then classical
strong uniqueness holds at $f$. Using Lemma 11.1 and (11.4), we will
estimate the deviation of $P_M^{(k)}(f)$ from $P_M(f)$. In what
follows, $\eta_k(f+M)$ will be an abuse of notation for
$\eta_k(\span\{f, M\})$.

\proclaim Proposition 11.3. Assume that $M$ is a finite dimensional
subspace in $X$ and a strong uniqueness type estimate of the form
(11.4) holds for a given $f\in X$. Then for any best approximant
$P_M^{(k)}(f)$ and any $k$ sufficiently large so that
$\eta_k(f+M)<1/2$, we have
$$\phi (\|P_M^{(k)}(f)-P_M(f)\|)\leq 8\|f\| \eta_k(f+M).$$
In particular,
$$\lim_{k\to\infty} \|P_M^{(k)}(f)-P_M(f)\|= 0.$$

\pf From (11.3) applied to both $P_M^{(k)}(f)$ and $f$, and since
$\|P_M^{(k)}(f)\|_k \le 2 \|f\|_k$, we have
$$\|f-P_M^{(k)}(f)\|\leq \|f\|+\|P_M^{(k)}(f)\|\leq
\|f\|+2\|P_M^{(k)}(f)\|_k\leq \|f\|+4\|f\|_k\leq 7\|f\|.$$ From the
definition of $\eta_k(f+M)$ and the above estimate, we have
$$|\;\|f-P_M^{(k)}(f)\|-\|f-P_M^{(k)}(f)\|_k|\leq \|f-P_M^{(k)}(f)\|\eta_k(f+M)
\le 7\|f\|\eta_k(f+M).$$ Thus
$$\eqalign{ \|f-P_M^{(k)}(f)\|&\;\leq\; \|f-P_M^{(k)}(f)\|_k+7\|f\|\eta_k(f+M)\;\leq\; \|f-P_M(f)\|_k+7\|f\|\eta_k(f+M)\cr&\;\leq\;
\|f-P_M(f)\|+\|f-P_M(f)\|\eta_k(f+M)+7\|f\|\eta_k(f+M)\cr &\;\leq\;
\|f-P_M(f)\|+8\|f\|\eta_k(f+M).\cr}$$ Combining now the last
estimate with (11.4) immediately yields the desired statement. The
limit follows from Lemma 11.1 since $$\lim_{k\to\infty}
\eta_k(f+M)=0$$ and $\phi(0)=0$ while $\phi(t)>0$ for $t>0$. \eop

\medskip\noindent {\bf Remark.} The above proposition yields the general
estimate
$$\|P_M^{(k)}(f)-P_M(f)\|\leq \phi^{-1}(c\eta_k(f+M)), \eqno (11.5)$$
where, by Lemma 11.1, the quantity on the right hand side tends to 0
as $k\rightarrow \infty$. In order to use this estimate to obtain
rates of convergence, one needs to find sharp bounds for both
$\eta_k$ and $\phi$. In many instances, this leads to sharp
estimates for the deviation of $P_M^{(k)}(f)$ from $P_M(f)$ (for
example, discrete Chebyshev approximation, P\'olya algorithm).
However, in some cases, this approach does not yield sharp bounds,
even when both $\eta_k$ and $\phi$ are precisely determined (this
happens, for example, in the case of discrete $L_1$-approximation).
In what follows, we provide a brief summary of these results.

\subsect{{12.} Discretization of Norms}

\medskip\noindent
It is considerably easier to solve an approximation problem when the
$L_p$-norm is replaced by a discrete $L_p$-norm. In this section, we
shall investigate the size of the error of this discretization
technique.

Let us denote by
$$\|f\|_k:= \max_{0\leq j\leq k}|f(j/k)|$$
and
$$\|f\|_{k,p}:={1\over {k^{1/p}}}(\sum_{j=0}^{k-1}|f(j/k)|^p)^{1/p},
\qquad 1\leq p<\infty$$ the discrete uniform and $L_p$-norms on
$[0,1]$, respectively.

In what follows, we deal with the example $M=\Pi_n$. We first
consider the case of the discrete uniform norm.

\proclaim Theorem 12.1. Let $f\in C^2[0,1]$ and $M=\Pi_n$. Set
$X:=C[0,1]$, with the usual uniform norm thereon, and
$$\|f\|_k:= \max_{0\leq j\leq k}|f(j/k)|.$$
Then there exists a constant $C$ such that
$$\eta_k(f+\Pi_n)\leq (C \|f''\|+4n^4(1+C\|f\|))/(8k^2)$$
for all $k$, and therefore
$$\|P_M^{(k)}(f)-P_M(f)\|=O(k^{-2}).$$

\pf We know from preceding results that classical strong uniqueness
holds in this case, i.e., we can use estimate (11.5) with $\phi(t) =
\gamma_n(f)t$, where $\gamma_n(f)$ is the strong unicity constant
for the function $f$. We now need to estimate $\eta_k(f+\Pi_n)$ for
$f\in C^2[0,1]\\ \Pi_n$. Let $g\in S(f+\Pi_n)$ and assume that
$\|g\|=g(x^*)=1$, $x^*\in [0,1]$. If $x^*=j/k$ for some $j\in
\{0,1\nek k\}$, then $\|g\|=\|g\|_k=1$. Thus, we may assume that
$x^*\neq j/k$, $j\in \{0,1\nek k\}$, i.e., in particular
$x^*\in(0,1)$ and $g'(x^*)=0$. Using Taylor's formula, we have that
for every $x\in [0,1]$,
$$g(x)=g(x^*)+{{g''(\xi)}\over 2}(x-x^*)^2 \eqno (12.1)$$
with some $\xi$ between $x$ and $x^*$. Obviously, we can choose a
$j\in \{0,1\nek k\}$ so that $|x^*-j/k|\le 1/2k.$ Setting $x=j/k$ in
(12.1) yields
$$1=\|g\|\leq \|g\|_k+{1 \over {8k^2}}\|g''\|.$$
Finally, $g=\alp f+p_n$ for some $p_n\in \Pi_n$. As $\|g\|=1$ and
$f\notin \Pi_n$, it easily follows that $|\alp| \le C = [\min_{p\in
\Pi_n} \|f-p\|]^{-1}$. Thus $\|p_n\|\leq 1+ C\|f\|$. Therefore by
the Markov inequality,
$$\|g''\|\leq C\|f''\|+4n^4(1+C\|f\|).$$
From these estimates, we obtain
$$\eta_k(f+\Pi_n)\leq (C \|f''\|+4n^4(1+C\|f\|))/(8k^2).$$
In particular, $\eta_k(f+\Pi_n)=O(1/k^2),$ and hence by (11.5), we
have $$\|P_M^{(k)}(f)-P_M(f)\|=O(k^{-2}). \meop$$

The above estimate turns out to be sharp, in general.

Let us now turn our attention to the discrete $L_p$-norms, $1<
p<\infty$.

\proclaim Theorem 12.2. Let $f\in C^1[0,1]$ and $M=\Pi_n$. Set
$X:=L_p[0,1]$, $1<p<\infty$, and
$$\|f\|_{k,p}:={1\over
{k^{1/p}}}(\sum_{j=0}^{k-1}|f(j/k)|^p)^{1/p}.$$ Then there exist
constants $C_p$ and $c_p$ such that
$$\eta_k(f+\Pi_n)\leq pc_p{{n^2}\over k}(1+C_p\|f\|_p+C_p\|f'\|_p)$$
for all $k$, and therefore
$$\|P_M^{(k)}(f)-P_M(f)\|_p=O(k^{-\theta_p}),\qquad 1< p<\infty,$$
where $\theta_p:=\min\{1/2,1/p\}$.

\pf Consider $f\in C^1[0,1]$ and let us estimate $\eta_k(f+\Pi_n)$
from above. Setting $\Delta_k:=[j/k,(j+1)/k]$, we obtain, by
repeated application of H\"{o}lder's inequality, for any $g\in
C^1[0,1]$,
$$\eqalign{\Bigl|\, \|g\|_p^p-\|g\|_{k,p}^p\Bigr| &\;=\;\Bigl|\int_0^1|g(x)|^p\dd x-{1 \over
k}\sum_{j=0}^{k-1}|g(j/k)|^p\Bigr|\cr &\;=\;\Bigl|
\,\sum_{j=0}^{k-1}\int_{\Delta_k}(|g(x)|^p-|g(j/k)|^p)\dd
x\Bigr|\;\leq\;
\sum_{j=0}^{k-1}\int_{\Delta_k}\left(\int_{j/k}^xp|g(\xi)|^{p-1}|g^{'}(\xi)|\dd
\xi\right)\dd x\cr &\;\leq\; {p \over
k}\sum_{j=0}^{k-1}\left(\int_{\Delta_k}|g(\xi)|^p\dd
\xi\right)^{1-1/p} \left(\int_{\Delta_k}|g'(\xi)|^p\dd
\xi\right)^{1/p}\cr &\;\leq\; {p\over k}
\left(\sum_{j=0}^{k-1}\int_{\Delta_k}|g(\xi)|^p\dd
\xi\right)^{1-1/p}
\left(\sum_{j=0}^{k-1}\int_{\Delta_k}|g'(\xi)|^p\dd \xi\right)^{1/p}
\;=\;{p\over k} \|g\|_p^{p-1}\|g'\|_p.\cr}$$ Thus
$$\Bigl|\, \|g\|_p^p-\|g\|_{k,p}^p\Bigr| \le {p\over
k} \|g\|_p^{p-1}\|g'\|_p$$ which in turn implies
$$\Bigl|\;\|g\|_p-\|g\|_{k,p}\Bigr|\leq {p \over k}\|g'\|_p\ . \eqno (12.2)$$

Let now $g:=\alp f+p_n\in S(f+\Pi_n)$. As $\|g||_p=1$, it follows,
as above, that $|\alp| \le C_p = [\min_{p\in \Pi_n}
\|f-p\|_p]^{-1}$. Thus $\|p_n\|_p\leq 1+C_p \|f\|_p$ and using the
$L_p$-Markov inequality, we obtain $\|p'_n\|_p\le
c_pn^2(1+C_p\|f\|_p).$ Combining this estimate with (12.2) yields
$$\eta_k(f+\Pi_n)\leq pc_p{{n^2}\over k}(1+C_p\|f\|_p+C_p\|f'\|_p).$$
Thus, we can use the estimate
$$\eta_k(f+\Pi_n)=O(k^{-1}),\qquad 1\leq p<\infty, \eqno (12.3)$$
in (11.5). Recall that by Corollary 7.2 and the subsequent Remark in
the case of the $L_p$-norm, $1<p<\infty$, we have the estimate
$\phi^{-1}(t)=O(t^{\theta_p})$ where $\theta_p:=\min\{1/2,1/p\}$.
Thus, we obtain by (11.5),
$$\|P_M^{(k)}(f)-P_M(f)\|_p=O(k^{-\theta_p}),\qquad 1< p<\infty. \meop$$

\medskip\noindent
{\bf Remark.} When $p=1$, Theorem 9.2 implies that for $f\in
C^1[0,1]$, non-classical strong uniqueness holds with
$\phi(t)=ct^2$. Thus repeated application of (12.3) and (11.5)
yields
$$\|P_M^{(k)}(f)-P_M(f)\|_1=O(k^{-1/2}).$$
When $p=1$, it is shown in Kro\'{o} [1981b] that if $f\in C^2[0,1]$
and $f-P_M(f)$ has a finite number of zeros then the above estimate
can be replaced by the sharp upper bound
$$\|P_M^{(k)}(f)-P_M(f)\|_1=O(k^{-1}).$$

\subsect{{13.} Asymptotic Representation of Weighted Chebyshev
Polynomials}

\medskip\noindent In this section, we shall use strong uniqueness
results in order to solve approximation problems in the case where
the norm is altered by a weight function. This approach will be used
to derive asymptotic representations for weighted Chebyshev
polynomials. Let $w>0$ be a positive weight on $[-1,1]$ and denote
by
$$T_{n,p}(x,w):=x^n+q_{n-1}(x), \qquad q_{n-1}\in \Pi_{n-1},$$
the monic polynomial that deviates least from $0$ in the weighted
$L_p$-norm on $[-1,1]$. That is,
$$\|T_{n,p}(x,w)w\|_p=\inf\{\|(x^n+g_{n-1}(x))w\|_p:g_{n-1}\in
\Pi_{n-1}\},\qquad 1\leq p\leq \infty.$$ In addition, let
$T_{n,p}^*(x,w):=T_{n,p}(x,w)/\|T_{n,p}(x,w)w\|_p$ be the normalized
Chebyshev polynomial. The problem of finding asymptotic
representations for weighted Chebyshev polynomials has been studied
by many authors, but satisfactory solutions were given only in the
case $p=2$; see Bernstein [1930]. In addition, in the case
$p=\infty$,  Bernstein [1930] gave an asymptotic formula for
$\|T_{n,\infty}(x,w)w\|_\infty$, Fekete, Walsh [1954/55] found the
$n$-th root asymptotics of these Chebyshev polynomials, while in a
more recent paper, Lubinsky, Saff [1987] gave their asymptotics
outside of $[-1,1]$. In this section, based on results of Kro\'o,
Peherstorfer [2007], [2008], we shall outline a complete solution to
this classical problem. This solution will be based on the strong
uniqueness estimates derived in the previous sections. Another basic
tool consists of the fact that for the specific weight $w=1/\rho_m$,
$\rho_m>0$, $\rho_m\in \Pi_m$, $m<n$, an explicit formula for the
minimal polynomials was already found by Chebyshev (see Akhiezer
[1947]).

\medskip\noindent{\bf The case  $p=\infty$.}
It is known (see Akhiezer [1947] ) that
$$T_{n,\infty}^*(\cos\phi, 1/\rho_m)=\Re\{z^{-n}g_m^2(z)\},\qquad z=\ee^{\ii\phi}, \eqno (13.1)$$
where $\rho_m\in \Pi_m$ is a polynomial positive on $[-1,1]$ and
$g_m\in \Pi_m$ is the certain polynomial (unique up to a
multiplicative constant of modulus $1$) all of whose zeros lie in
$|z|>1$ and such that
$$|g_m(\ee^{\ii\phi})|^2=\rho_m(\cos\phi),\qquad \phi\in[0,\pi].$$

Moreover, the corresponding $L_1$-Chebyshev polynomial is given by
$$T_{n-1,1}^*(\cos\phi, 1/\rho_m)=-\Im\{z^{-n}g_m^2(z)\}/\sin\phi,\qquad z=\ee^{\ii\phi}.$$
In addition, these polynomials satisfy the relation
$$\Re\{z^{-n}g_m^2(z)\}^2+\Im\{z^{-n}g_m^2(z)\}^2=\rho_m^2. \eqno (13.2)$$

It follows immediately from (13.2) that the zeros of
$\Im\{z^{-n}g_m^2(z)\}$ are the $n+1$ equioscillation points of
$T_{n,\infty}^*(\cos\phi, 1/\rho_m)$. Moreover, in view of the
alternation theorem (Theorem 4.4), (13.2) also implies that
$$T_{n-1,\infty}^*(\cos\phi,
\sqrt{1-x^2}/\rho_m)=\Im\{z^{-n}g_m^2(z)\}/\sin\phi. \eqno (13.3)$$
So a natural idea is to approximate a positive weight $w$ by
reciprocals of positive polynomials $1/\rho_m$ (the degree of
approximation can be estimated by the Jackson Theorem), and then
apply strong uniqueness type results in order to obtain the required
asymptotics. Clearly, we shall need rather precise strong uniqueness
type estimates that also take into account the dependence upon $n$.

We first recall that the positive continuous weight function
$w(\cos\phi)$ can be represented in the form
$$w(\cos\phi)={1\over {|\pi(\ee^{\ii\phi})|^2}},$$
where $\pi(z)$, the so called Szeg\H{o} function of $w$, is the
function nonzero and analytic in $|z|<1$ given by the formula
$$\pi(z):=\exp\left\{-{1\over {4\pi}}\int_0^{2\pi}{{\ee^{\ii\phi}+z}\over {\ee^{\ii\phi}-z}}\log
w(\cos\phi)\dd\phi\right\}.$$

Let us denote by $C^{k+\alpha}[-\pi,\pi]$ the class of
$2\pi$-periodic functions whose $k$-th derivative, $k\in \NN$,
satisfies the Lip $\alpha$ property. This next result provides an
asymptotic formula for the weighted $L_{\infty}$-Chebyshev
polynomials for positive weights $w(\cos\phi)\in
C^{2+\alpha}[-\pi,\pi]$.

\proclaim Theorem 13.1. Let $w(\cos\phi)\in C^{2+\alpha}[-\pi,\pi]$
with $0<\alpha<1$, $w(x)>0$, $x\in [-1,1]$. Then
$$T_{n,\infty}^*(\cos\phi,w)=\Re
\{\ee^{-\ii n\phi}(\pi(\ee^{\ii\phi}))^2\}+O( n^{-\alpha})$$
uniformly for $\phi\in[0,\pi]$.

\medskip
First, we shall need a lemma estimating the deviation between
Chebyshev polynomials corresponding to different weights via the
strong unicity constant.

\proclaim Lemma 13.2.  Let $w_1, w_2\in C[-1,1]$ be positive weight
functions. Then
$$\|T_{n,\infty}^*(\cdot,w_1)w_1-T_{n,\infty}^*(\cdot,w_2)w_2\|
\leq {c \over {\gamma_{n-1}(T_{n,\infty}^*(w_1))}}\|w_1-w_2\|,$$
where $\gamma_{n-1}(T_{n,\infty}^*(w_1))$ is the strong unicity
constant of $T_{n,\infty}^*(w_1)w_1$ with respect to the Haar space
$w_1\Pi_{n-1}$, and $c>0$ depends only on $w_1, w_2$.

\pf  We shall denote below by $c_j$ constants depending only on
$w_1, w_2$. Clearly,
$$T_{n,\infty}^*(\cdot,w_1)=a_nx^n-p_1=:T_1^*,\quad
T_{n,\infty}^*(\cdot,w_2)=b_nx^n-p_2=:T_2^*,\qquad p_1, p_2\in
\Pi_{n-1}.$$ Denote by $T_1:=T_1^*/a_n, T_2:=T_2^*/b_n$ the
corresponding monic polynomials. Then using the extremal property of
$T_1$, we obtain
$${1\over {a_n}}=\|T_1w_1\|\leq \|T_2w_1\|\leq
\left\|T_2w_2 {{w_1}\over {w_2}}\right\|\leq {1 \over
{b_n}}\left\|{{w_1}\over {w_2}}\right\|.$$ Replacing $T_1$ by $T_2$,
we can obtain an analogous inequality
$${1 \over {b_n}}\leq {1\over {a_n}}\left\|{{w_2}\over {w_1}}\right\|.$$  Thus
combining these estimates, we have
$$1-c_1\|w_1-w_2\|\leq {1\over {\left\|{{w_1}\over {w_2}}\right\|}}\leq
{{a_n}\over {b_n}}\leq \left\|{{w_2}\over {w_1}}\right\|\leq
c_2\|w_1-w_2\|+1.$$ It therefore follows that
$$\left|{{a_n}\over {b_n}}-1\right|\leq c_3\|w_1-w_2\|. \eqno (13.4)$$
Set now
$$q:=p_1-{{a_n}\over {b_n}}p_2=T_1^*- {{a_n}\over
{b_n}}T_2^* \in \Pi_{n-1}. \eqno (13.5)$$ Since $0$ is the best
approximant to $w_1T_1^*$ in $w_1\Pi_{n-1}$ in the uniform norm, it
follows by the classical strong uniqueness inequality (applied to
the Haar space $w_1\Pi_{n-1}$ ) together with (13.4) and (13.5) that
$$\eqalign{\gamma_{n-1}(T_1^*)\|qw_1\|&\;\leq\;
\|(T_1^*-q)w_1\|-\|T_1^*w_1\|\;=\;{{a_n}\over
{b_n}}\|T_2^*w_1\|-\|T_1^*w_1\|\cr &\;=\;{{a_n}\over
{b_n}}\|T_2^*w_1\|-1\;\leq\; {{a_n}\over {b_n}}\left\|{{w_1}\over
{w_2}}\right\|-1\;\leq\; {{a_n}\over
{b_n}}(1+c_4\|w_1-w_2\|)-1\cr&\;\leq\;
(1+c_3\|w_1-w_2\|)(1+c_4\|w_1-w_2\|)-1\;\leq\; c_5\|w_1-w_2\|,\cr}
\eqno (13.6)$$ where $\gamma_{n-1}(T_1^*)$ is the strong unicity
constant of $T_1^*w_1$ with respect to the Haar space
$w_1\Pi_{n-1}$.

Finally, using (13.4) and (13.6), we arrive at
$$\eqalign{\|T_1^*w_1-T_2^*w_2\|&\;\leq\;
\|qw_1\|+\|T_2^*(w_2- {{a_n}\over {b_n}}w_1)\|\cr&\;\leq\;
\|qw_1\|+\|T_2^*(w_2-w_1)\|+\|T_2^*w_1\left(1- {{a_n}\over
{b_n}}\right)\|\leq {{c_6}\over {\gamma_{n-1}(T_1^*)}}
\|w_1-w_2\|,\cr}$$ which is the required estimate. \eop

We now need to estimate the strong unicity constant
$\gamma_{n-1}(T_{n,\infty}^*(w_1))$ (which is the strong unicity
constant of $T_{n,\infty}^*(w_1)w_1$ with respect to the Haar space
$w_1\Pi_{n-1}$) in the special case when $w_1=1/\rho_m$. Moreover,
it is important in deriving asymptotic relations that the dependence
of this quantity on $n$ is revealed. We shall use Theorem 4.10 of
Section 4 in order to provide a precise bound.

\proclaim Lemma 13.3.  Let $\rho_m\in \Pi_{m}$, $m<n$, be such that
$0<A\leq \rho_m\leq B$ on $[-1,1]$. Then
$$\gamma_{n-1}(T_{n,\infty}^*(\cdot, 1/\rho_m))\geq {c \over {n^2}}, \eqno (13.7)$$
where the constant $c$ depends only on $A$ and $B$.

\pf Let $-1=x_1< \cdots <x_{n+1}=1$ be the equioscillation points of
$T_{n,\infty}^*(\cdot, 1/\rho_m)/\rho_m $ and $m_k\in \Pi_{n-1}$,
$1\leq k\leq n+1$, be the polynomials of Theorem 4.10 with respect
to the Haar space $\Pi_{n-1}/\rho_m$, i.e.,
$$ m_k(x_j)/\rho_m(x_j)= \sgn T_{n,\infty}^*(x_j, 1/\rho_m),\qquad 1\leq j\leq n+1,\quad j\neq
k.$$ By the second statement of Theorem 4.10,
$$\gamma_{n-1}(T_{n,\infty}^*(\cdot, 1/\rho_m))=\min_{1\leq k\leq
n+1}{1 \over {\|m_k/\rho_m\|}}. \eqno (13.8)$$

Set now for any $1\leq k\leq n+1$,
$$Q_k(x):=(x-x_k)T_{n,\infty}^*(x, 1/\rho_m)+(1-x^2)T_{n-1,\infty}^*(x,
\sqrt{1-x^2}/\rho_m).$$ It follows by (13.2) and (13.3) that the
equioscillation points of $T_{n,\infty}^*(x, 1/\rho_m)/\rho_m$ are
the zeros of $(1-x^2)T_{n-1,\infty}^*(x, \sqrt{1-x^2}/\rho_m)$.
Moreover, since $m<n$, the relation (13.2) also implies that the
leading coefficients of  $T_{n,\infty}^*(x, 1/\rho_m)$ and
$T_{n-1,\infty}^*(x, \sqrt{1-x^2}/\rho_m)$ coincide. It therefore
follows that $Q_k\in \Pi_{n}$, $Q_k(x_k)=0$, $1\leq k\leq n+1$, and
$$m_k(x)={{Q_k(x)}\over {x-x_k}},\qquad 1\leq k\leq n+1.$$
Obviously, $\|Q_k\|\leq c$, $1\leq k\leq n+1$, with a constant
depending only on $A,B$. This easily implies,  by the Markov
inequality, that $ \|m_k\|\leq c_1n^2$, $1\leq k\leq n+1$. Applying
this inequality together with (13.8) yields (13.7). \eop

\noindent{\bf Proof of Theorem 13.1:} When $0<\alpha<1$, the
smoothness condition imposed on the weight $w$ implies that its
Szeg\H{o} function $\pi(z)$ is also $C^{2+\alpha}$ on $|z|=1$ (see
Kro\'o, Peherstorfer [2008] for details). Therefore there exists a
sequence of polynomials $g_m\in \Pi_m$ such that uniformly on
$|z|\leq 1$
$$|\pi(z)-g_m(z)|\leq {c\over {m^{2+\alpha}}}. \eqno (13.9)$$
Since $\pi(z)\neq 0$ on $|z|\leq 1$, it follows that for $m\geq m_0$
the function $g_m$ also does not vanish in $|z|\leq 1$ and
$$\rho_m(\cos\phi):=|g_m(\ee^{\ii\phi})|^2\geq c>0. \eqno (13.10)$$
Furthermore, by Chebyshev's result (see (13.1)),
$$T_{n,\infty}^*(\cos\phi,1/\rho_m)=\Re\{z^{-n}g_m^2(z)\},\qquad z=\ee^{\ii\phi}. \eqno (13.11)$$

We are now ready to verify the statement of the theorem. Indeed,
since
$$w(\cos\phi)={1\over {|\pi(\ee^{\ii\phi})|^2}},$$
we have by (13.9) and (13.10) that
$$\left|w(x)-{1\over {\rho_m(x)}}\right|\leq {{c_1}\over {m^{2+\alpha}}},\qquad x\in [-1,1]. \eqno
(13.12)$$ Thus using Lemmas 13.2 and 13.3 together with (13.9),
(13.10) and (13.12), see also (13.11), we obtain
$$\|T_{n,\infty}^*(\cos\phi,w)-\Re\{\ee^{-\ii n\phi}(\pi(\ee^{\ii\phi}))^2\}\|\leq
\|T_{n,\infty}^*(\cos\phi,w)-T_{n,\infty}^*(\cos\phi,1/\rho_m)\|+$$
$$\|T_{n,\infty}^*(\cos\phi,1/\rho_m)-\Re\{\ee^{-\ii n\phi}(\pi(\ee^{\ii\phi}))^2\}\|\leq
{{c_2n^2}\over {m^{2+\alpha}}},\qquad x=\cos\phi.$$ Putting in the
last estimate $m:=\lfloor n/2\rfloor$ yields the statement of the
theorem. \eop

\medskip\noindent{\bf The case  $1\leq p<\infty$.}
Similarly to the case $p=\infty$, the following explicit
representation for weighted Chebyshev polynomials is given in
Akhiezer [1947] for the weight function $\rho_{m,p}:={1\over
{\sqrt{1-x^2}\rho_m^{p/2}}}$, where $\rho_m\in \Pi_m$ is positive on
$[-1,1]$,
$$T_{n,p}^*(\cos\phi,\rho_{m,p})=\lambda_p\Re\{z^{-n}g_m(z)\},\qquad z=\ee^{\ii\phi},$$
where $\lambda_p$ is a constant depending only on $p$ and, as above,
$g_m$ is related to $\rho_m$ by
$$|g_m(\ee^{\ii\phi})|^2=\rho_m(\cos\phi),\qquad \phi\in[0,\pi].$$ Let now $w$
be a positive weight and set $w_{p}:={{w^{p/2}}\over
{\sqrt{1-x^2}}}$, $1\leq p<\infty$.

The next theorem provides an asymptotic representation for
$L_p$-Chebyshev polynomials when $1< p<\infty$.

\proclaim Theorem 13.4. Let $1<p<\infty$. Then for any positive
weight $w$ such that $w\in C^{\alpha}[-\pi,\pi]$, $0<\alpha<1$,
$$T_{n,p}^*(\cos\phi,w_{p})=\lambda_p\Re\{z^{-n}\pi(z)\}+O(n^{-\alpha\theta_p}),\qquad z=\ee^{\ii\phi},$$
where $\pi(z)$ is the Szeg\H{o} function of $w$,
$\theta_p:=\min\{1/2,1/p\}$, and the $O(\cdot)$ is taken with
respect to the $L_p$-norm.

Since the main idea in the proof of Theorem 13.4 is similar to the
proof of Theorem 13.1, we shall just briefly outline the proof.
First, it is shown that, similarly to Lemma 13.2, under suitable
conditions on the weight functions $w_1, w_2>A>0$, we have
$$\|T_{n,p}^*(\cos\phi,w_{1})-T_{n,p}^*(\cos\phi,w_{2})\|_p\leq
c_{A,p}\|w_1-w_2\|_p^{\theta_p}.$$ It should be noted that instead
of classical strong uniqueness, we use here the fact that in
$L_p$-approximation with $1<p<\infty$, non-classical strong
uniqueness with $\phi^{-1}(t)=O(t^{\theta_p})$ is satisfied, and the
constants involved depend only on $p$.

Finally, similar asymptotic relations can also be given when $p=1$.
In this case, non-classical strong uniqueness with
$\phi^{-1}(t)=O(t^{1/2})$ takes place (this follows from Theorem 9.2
applied in the case of Lip 1 functions). In addition, similarly to
Lemma 13.3, one has to study how the constants involved depend on
$n$. See Kro\'o, Peherstorfer [2008] for details.




\sect{ References}

N.~I.~Achieser [1930], On extremal properties of certain rational
functions, {\sl Doklady Akad. Nauk}, 495--499.

N.~I.~Achieser [1947], {\sl Lectures on the Theory of
Approximation}, OGIZ, Moscow-Leningrad, (in Russian) 1947, reprinted
in English as {\sl Theory of Approximation}, Frederick Ungar
Publishing Co., New York, 1956, and reissued by Dover Publications
in 1992.

J.~Angelos and A.~Egger [1984], Strong uniqueness in $L^p$ spaces,
{\sl J. Approx. Theory} {\bf 42}, 14--26.

J.~R.~Angelos, M.~S.~Henry, E.~H.~Kaufman, Jr., and T.~D.~Lenker
[1985], Local Lipschitz constants, {\sl J. Approx. Theory} {\bf 43},
53--63.

J.~R.~Angelos, M.~S.~Henry, E.~H.~Kaufman, A.~Kro\'o and
T.~D.~Lenker [1986a], Local Lipschitz and strong unicity constants
for certain nonlinear families, in {\sl Approximation Theory V},
C.~K.~Chui, L.~L.~Schumaker, J.~D.~Ward, Eds., 239--242, Academic
Press, New York.

J.~R.~Angelos, M.~S.~Henry, E.~H.~Kaufman, A.~Kro\'o and
T.~D.~Lenker [1986b], Local and global Lipschitz constants, {\sl J.
Approx. Theory} {\bf 46}, 137--156.

J.~R.~Angelos, M.~S.~Henry, E.~H.~Kaufman and T.~D.~Lenker  [1986],
Extended Lipschitz constants, in {\sl Approximation Theory V},
C.~K.~Chui, L.~L.~Schumaker, J.~D.~Ward, Eds., 243--246, Academic
Press, New York.

J.~R.~Angelos, M.~S.~Henry, E.~H.~Kaufman, Jr., and T.~D.~Lenker
[1988], Bounds for extended local Lipschitz constants, in {\sl
Constructive Theory of Functions} (Varna, 1987),  17--26, Publ.
House Bulgar. Acad. Sci., Sofia.

J.~R.~Angelos, M.~S.~Henry, E.~H.~Kaufman, Jr., T.~D.~Lenker and
A.~Kro\'o  [1989], Local Lipschitz and strong unicity constants for
certain nonlinear families, {\sl J. Approx. Theory} {\bf 58},
164--183.

J.~Angelos, E.~Kaufman, Jr., T.~Lenker and M.~S.~Henry  [1991],
Bounds for extended Lipschitz constants, {\sl Acta Math. Hungar.}
{\bf 58}, 81--93.

J.~R.~Angelos and A.~Kro\'o  [1986], The equivalence of the moduli
of continuity of the best approximation operator and of strong
unicity in $L^1$, {\sl J. Approx. Theory} {\bf 46}, 129--136.

J.~Angelos and D.~Schmidt  [1983], Strong uniqueness in $L^1(X,
\Sigma, \mu)$, in {\sl Approximation Theory IV}, C.~K.~Chui,
L.~L.~Schumaker, J.~D.~Ward, Eds., 297--302, Academic Press, New
York.

J.~R.~Angelos and  D.~Schmidt  [1988], The prevalence of strong
uniqueness in $L^1$, {\sl Acta Math. Hungar.} {\bf 52}, 83--90.

R.~B.~Barrar and H.~L.~Loeb [1970], On the continuity of the
nonlinear Tschebyscheff operator, {\sl Pacific J. Math.} {\bf 32},
593--601.

R.~B.~Barrar and H.~L.~Loeb [1986], The strong uniqueness theorem
for monosplines, {\sl J. Approx. Theory} {\bf 46}, 157--169.

M.~W.~Bartelt  [1974], Strongly unique best approximates to a
function on a set, and a finite subset thereof, {\sl Pacific J.
Math.} {\bf 53}, 1--9.

M.~Bartelt  [1975], On Lipschitz conditions, strong unicity and a
theorem of A.~K.~Cline, {\sl J. Approx. Theory} {\bf 14}, 245--250.

M.~Bartelt [2001], Hausdorff strong unicity in vector-valued
Chebyshev approximation on finite sets, in {\sl Trends in
Approximation Theory} (Nashville, TN, 2000), K.~Kopotun, T.~Lyche,
M.~Neamtu (eds.), 31--38, Innov. Appl. Math., Vanderbilt Univ.
Press, Nashville, TN.

M.~W.~Bartlet and M.~S.~Henry  [1980], Continuity of the strong
uniqueness constant on $C[X]$ for changing $X$,  {\sl J. Approx.
Theory} {\bf 28}, 87--95.

M.~Bartelt, E.~H.~Kaufman, Jr., and J.~Swetits  [1990], Uniform
Lipschitz constants in Chebyshev polynomial approximation, {\sl J.
Approx. Theory} {\bf 62}, 23--38.

M.~W.~Bartlet, A.~Kro\'o and J.~J.~Swetits  [1989], Local Lipschitz
constants: characterization and uniformity, in {\sl Approximation
Theory VI, Volume I}, C.~K.~Chui, L.~L.~Schumaker, J.~D.~Ward, Eds.,
65--68, Academic Press, New York.

M.~Bartelt and W.~Li [1995a], Error estimates and Lipschitz
constants for best approximation in continuous function spaces, {\sl
Comput. Math. Appl.} {\bf 30}, 255--268.

M.~Bartelt and W.~Li [1995b], Haar theory in vector-valued
continuous function spaces, in {\sl Approximation Theory VIII,
Volume I}, (College Station, TX, 1995), C.~K.~Chui, L.~L.~Schumaker
(eds.), 39--46, Ser. Approx. Decompos., 6, World Sci. Publ., NJ,
1995.

M.~W.~Bartlet and H.~W.~McLaughlin  [1973], Characterizations of
strong unicity in approximation theory, {\sl J. Approx. Theory} {\bf
9}, 255--266.

M.~W.~Bartelt and D.~Schmidt  [1980],  On strong unicity and a
conjecture of Henry and Roulier, in {\sl Approximation Theory III},
E.~W.~Cheney, Ed., 187--191, Academic Press, New York.

M.~W.~Bartelt and D.~Schmidt  [1981], On Poreda's problem on the
strong unicity constants, {\sl J. Approx. Theory} {\bf 33}, 69--79.

M.~W.~Bartelt and D.~Schmidt  [1984], Lipschitz conditions, strong
uniqueness, and almost Chebyshev subspaces of $C(X)$, {\sl J.
Approx. Theory} {\bf 40}, 202--215.

M.~Bartelt and J.~Swetits  [1983], Uniform strong unicity constants
for subsets of $C[a,b]$, in {\sl  Approximation Theory IV}, (College
Station, Tex., 1983), C.~K.~Chui, L.~L.~Schumaker, J.~D.~Ward, Eds.,
329--334, Academic Press, New York.

M.~Bartelt and J.~Swetits  [1986], Uniform Lipschitz constants and
almost alternation sets, in {\sl Approximation Theory V},
C.~K.~Chui, L.~L.~Schumaker, J.~D.~Ward, Eds., 247--250, Academic
Press, New York.

M.~Bartelt and J.~Swetits  [1988], Uniform strong unicity constants
for subsets of $C(X)$, {\sl J. Approx. Theory}  {\bf 55}, 304--317.

M.~W.~Bartelt and J.~J.~Swetits  [1991a], Local Lipschitz constants
and Kolushov polynomials, {\sl Acta Math. Hungar.} {\bf 57},
259--263.

M.~W.~Bartelt and J.~J.~Swetits  [1991b], New classes of local
Lipschitz constants for the best approximation operator on finite
sets, in {\sl Approximation Theory} (Kecskem\'et, 1990),  69--83,
Colloq. Math. Soc. J\'anos Bolyai, {\bf 58}, North-Holland,
Amsterdam.

M.~W.~Bartelt and J.~J.~Swetits  [1993], The strong derivative of
the best approximation operator, {\sl Numer. Funct. Anal. Optim.}
{\bf 14}, 229--248.

M.~W.~Bartelt and J.~J.~Swetits  [1995], Continuity properties of
Lipschitz constants for the best approximation operator, in {\sl
Approximation and Optimization in the Caribbean, II} (Havana, 1993),
51--62, Approx. Optim., {\bf 8}, Lang, Frankfurt am Main.

M.~W.~Bartelt and J.~J.~Swetits  [2002], Uniform strong unicity of
order 2 for generalized Haar sets, in {\sl Approximation Theory X,
Abstract and Classical Analysis} (St. Louis, MO, 2001), C.~K.~Chui,
L.~L.~Schumaker, J.~St\"ockler, Eds., 23--30, Vanderbilt Univ.
Press, Nashville, TN.

M.~Bartelt and  J.~Swetits  [2007], Lipschitz continuity and Gateaux
differentiability of the best approximation operator in
vector-valued Chebyshev approximation, {\sl  J. Approx. Theory} {\bf
148}, 177--193.

M.~Bartelt and J.~Swetits  [2008], Lipschitz continuity of the best
approximation operator in vector-valued Chebyshev approximation {\sl
J. Approx. Theory}  {\bf 152}, 161--166.

S.~N.~Bernstein [1930], Sur les polynomes orthogonaux relatifs \'a
un segment fini, {\sl J. Math.} {\bf 9}, 127--177.

B.~O.~Bj\"{o}rnest{\aa}l [1975], Continuity of the metric projection
operator II, TRITA-MAT-1975-12, Preprint Series of the Royal
Institute of Technology, Stockholm, Sweden.

B.~O.~Bj\"{o}rnest{\aa}l [1979], Local Lipschitz continuity of the
metric projection operator, in  {\sl Approximation Theory} (Papers,
VIth Semester, Stefan Banach Internat. Math. Center, Warsaw, 1975),
pp. 43--53, Banach Center Publ., {\bf 4}, PWN, Warsaw.

H.~P.~Blatt  [1984a], On strong uniqueness in linear complex
Chebyshev approximation, {\sl J. Approx. Theory} {\bf 41}, 159--169.

H.~P.~Blatt  [1984b], Exchange algorithms, error estimations and
strong unicity in convex programming and Chebyshev approximation, in
{\sl Approximation Theory and Spline Functions}, St.~John's, Nfld,
1983, NATO Adv. Sci. Inst. Sec. C Math. Phys. Sci., {\bf 136},
23--63, Reidel, Dordrecht.

H.~P.~Blatt  [1985], Characterization of strong unicity in
semi-infinite optimization by chain of references, in {\sl
Parametric Optimization and Approximation} (Oberwolfach, 1983) ISNM
{\bf 72}, 36--46, Birkh\"auser.

H.~P.~Blatt  [1986], Lipschitz continuity and strong unicity in G.
Freud's work, {\sl J. Approx. Theory} {\bf 46}, 25--31.

H.~P.~Blatt  [1990], On the distribution of the points of Chebyshev
alternation with applications to strong unicity constants, {\sl Acta
Math. Hungar.} {\bf 55}, 75--82.

E.~Borel [1905], {\sl Le{\c c}ons sur les Fonctions de Variables
R\'eelles}, Gauthier-Villars, Paris.

A.~P.~Bosznay  [1988], A remark concerning strong uniqueness of
approximation, {\sl Studia Sci. Math. Hungar.} {\bf 23}, 85--87.

D. Braess [1973], Kritische Punkte bei der nichtlinearen
Tschebyscheff-Approximation, {\sl Math. Z.} {\bf 132}, 327--341.

D.~Braess  [1986], {\sl Nonlinear Approximation Theory},
Springer-Verlag, Berlin Heidelberg.

B.~Brosowski [1983], A refinement of the Kolmogorov-criterion, in
{\sl Constructive Function Theory '81} (Varna, 1981), 241--247,
Bulgar. Acad. Sci., Sofia.

B.~Brosowski and C.~Guerreiro  [1986], An extension of strong
uniqueness to rational approximation, {\sl J. Approx. Theory} {\bf
46}, 345--373.

B.~Brosowski and C.~Guerreiro  [1987], On the uniqueness and strong
uniqueness of best rational Chebyshev-approximation, {\sl Approx.
Theory Appl.} {\bf 3}, 49--70.

B.~L.~Chalmers, F.~T.~Metcalf and G.~D.~Taylor  [1983], Strong
unicity of arbitrary rate, {\sl J. Approx. Theory} {\bf 37},
238--245.

B.~L.~Chalmers and G.~D.~Taylor  [1980], On the existence of strong
unicity of arbitrarily small order, in {\sl Approximation Theory
III}, E.~W.~Cheney, Ed., 293--298, Academic Press, New York.

B.~L.~Chalmers and G.~D.~Taylor [1983], A unified theory of strong
uniqueness in uniform approximation with constraints, {\sl J.
Approx. Theory} {\bf 37}, 29--43.

K.~Y.~Chan, Y.~M.~Chen, M.~C.~Liu and  S.~M.~Ng  [1982], An example
on strong unicity constants in trigonometric approximation, {\sl
Proc. Amer. Math. Soc.} {\bf 84}, 79--84.

E.~W.~Cheney [1965], Approximation by generalized rational
functions, in {\sl Approximation of Functions}, H. L. Garabedian,
Ed., 101--110, Elsevier, Amsterdam.

E.~W.~Cheney  [1966], {\sl Introduction to Approximation Theory},
McGraw-Hill Book Co., New York-Toronto.

E.~W.~Cheney and H. L. Loeb [1964], Generalized rational
approximation, {\sl J. SIAM Numer. Anal. Ser. B} {\bf 1}, 11--25.

E.~W.~Cheney and D.~E.~Wulbert  [1969], The existence and unicity of
best approximations, {\sl Math.~Scand.} {\bf 24}, 113--140.

J.~A.~Clarkson  [1936], Uniformly convex spaces, {\sl Trans. Amer.
Math. Soc.} {\bf 40}, 396--414.

A.~K.~Cline  [1973], Lipschitz conditions on uniform approximation
operators, {\sl J. Approx. Theory} {\bf 8}, 160--172.

L.~Cromme  [1977/78], Strong uniqueness: a far reaching criterion
for the convergence analysis of iterative procedures,  {\sl
Numer.~Math.} {\bf 29}, 179--193.

L.~Danzer, B.~Gr\"unbaum and V.~Klee [1963], Helly's theorem and its
relatives, in {\sl Convexity}, {\sl Proc.~Symposia in Pure
Mathematics, Volume VII}, 101--180, AMS.

F.~Deutsch and W.~Li [1991], Strong uniqueness, Lipschitz
continuity, and continuous selections for metric projections in
$L_1$,  {\sl J. Approx. Theory} {\bf 66}, 198--224.

F.~Deutsch and S.~Mabizela [1996], Best interpolatory approximation
in normed linear spaces, {\sl J. Approx. Theory} {\bf 85}, 250--268.

C.~B.~Dunham [1980], A uniform constant of strong uniqueness on an
interval, {\sl J. Approx. Theory} {\bf 28}, 207--211.

C.~B.~Dunham [1989], Local strong uniqueness for nonlinear
approximation, {\sl Approx. Theory Appl.} {\bf 5}, 43--45.

C.~B.~Dunham [1995], Uniform local strong uniqueness on finite
subsets, {\sl Aequationes Math.} {\bf 49}, 295--299.

A.~G.~Egger and G.~D.~Taylor [1983], Strong uniqueness in convex
$L\sp{p}$ approximation, in {\sl  Approximation Theory IV},
C.~K.~Chui, L.~L.~Schumaker, J.~D.~Ward, Eds., 451--456, Academic
Press, New York.

A.~Egger and G.~D.~Taylor [1989a], A survey of local and directional
local strong uniqueness, in {\sl Approximation Theory VI}, Vol. I,
 C.~K.~Chui, L.~L.~Schumaker, J.~D.~Ward, Eds.,
239--242, Academic Press, Boston, MA.

A.~G.~Egger and G.~D.~Taylor [1989b], Local strong uniqueness, {\sl
J. Approx. Theory} {\bf 58}, 267--280.

N.~H.~Eggert and J.~R.~Lund [1984], Examples of functions whose
sequence of strong unicity constants is unbounded, {\sl J. Approx.
Theory} {\bf 41}, 244--252.

M.~Fang [1990], The Chebyshev theory of a variation of $L_p$
$(1<p<\infty)$ approximation,  {\sl J. Approx. Theory} {\bf 62},
94--109.

M.~Fekete and J.~L.~Walsh [1954/55], On the asymptotic behaviour of
polynomials with extremal properties, and of their zeros, {\sl J.
d'Analyse Math.} {\bf 4}, 49--87.

T.~Fischer [1989], Strong unicity and alternation in linear
approximation and a continuous alternator, in {\sl Approximation
Theory VI}, Vol. I, C.~K.~Chui, L.~L.~Schumaker, J.~D.~Ward, Eds.,
255--258, Academic Press, New York.

T.~Fischer [1990], Strong unicity in normed linear spaces, {\sl
Numer. Funct. Anal. Optim.} {\bf 11}, 255--266.

Y.~Fletcher and J.~A.~Roulier [1979], A counterexample to strong
unicity in monotone approximation, {\sl J. Approx. Theory} {\bf
27}, 19--33.

G.~Freud [1958], Eine Ungleichung f\"ur Tschebyscheffsche
Approximationspolynome, {\sl  Acta Sci. Math. Szeged}  {\bf 19},
162--164.

A.~L.~Garkavi [1959], Dimensionality of polyhedra of best
approximation for differentiable functions, {\sl Izv. Akad. Nauk
SSSR Ser. Mat.} {\bf 23}, 93--114. (Russian)

A.~L.~Garkavi [1964], On \v{C}eby\v{s}ev and almost \v{C}eby\v{s}ev
subspaces, {\sl Izv.~Akad.~Nauk SSSR Ser.~Mat.} {\bf 28}, 799--818;
in English translation in {\sl Amer.~Math. Soc. Transl.} {\bf 96},
(1970) 153--175.

A.~L.~Garkavi [1965], Almost \v{C}eby\v{s}ev subspaces of continuous
functions, {\sl Izv.~Vys\v{s}.~U\v{c}ebn. Zaved.~Matematika} {\bf
1965}, 36--44; in English translation in  {\sl Amer.~Math. Soc.
Transl.} {\bf 96}, (1970) 177--187.

W.~Gehlen [1999], On a conjecture concerning strong unicity
constants, {\sl J. Approx. Theory} {\bf 101}, 221--239.

W.~Gehlen [2000], Unboundedness of the Lipschitz constants of best
polynomial approximation, {\sl J. Approx. Theory} {\bf 106},
110--142.

R.~Grothmann [1988a], On the real $CF$-method for polynomial
approximation and strong unicity constants, {\sl J. Approx. Theory}
{\bf 55}, 86--103.

R.~Grothmann [1988b], Local uniqueness in nonuniqueness spaces, {\sl
Approx. Theory Appl.} {\bf 4}, 35--39.

R.~Grothmann [1989], A note on strong uniqueness constants, {\sl  J.
Approx. Theory} {\bf 58}, 358--360.

R.~Grothmann [1999], On the problem of Poreda, in {\sl Computational
Methods and Function Theory 1997} (Nicosia), N.~Papamichael,
St.~Ruscheweyh and E.~B.~Saff, Eds., 267--273, Ser. Approx.
Decompos., {\bf 11}, World Sci. Publ., NJ.

M.~Gutknecht [1978], Non-strong uniqueness in real and complex
Chebyshev approximation, {\sl J. Approx. Theory} {\bf 23}, 204--213.

A.~Haar [1918], Die Minkowskische Geometrie und die Ann\"aherung an
stetige Funktionen, {\sl Math. Ann.} {\bf 78}, 294--311.

O.~Hanner [1956], On the uniform convexity of $L^p$ and $\ell^p$,
{\sl Ark. Mat.} {\bf 3}, 239--244.

S.~Ja.~Havinson and Z.~S.~Romanova [1972], Approximation properties
of finite-dimensional subspaces in $L_1$, {\sl Mat. Sb.} {\bf 89},
3--15, in English translation in {\sl Math. USSR Sb.} {\bf 18},
1--14.

M.~S.~Henry [1987], Lipschitz and strong unicity constants, in {\sl
A.~Haar Memorial Conference}, Vol. I, (Budapest, 1985), 423--444,
Colloq. Math. Soc. J\'anos Bolyai, {\bf 49}, North-Holland,
Amsterdam.

M.~S.~Henry and L.~R.~Huff [1979], On the behavior of the strong
unicity constant for changing dimension, {\sl J. Approx. Theory}
{\bf 27}, 278--290.

M.~S.~Henry, E.~H.~Kaufman, Jr., and T.~D.~Lenker [1983], Lipschitz
constants for small perturbations, in {\sl Approximation Theory IV},
C.~K.~Chui, L.~L.~Schumaker, J.~D.~Ward, Eds., 515--520, Academic
Press, New York.

M.~S.~Henry, E.~H.~Kaufman, Jr., and T.~D.~Lenker [1983], Lipschitz
constants on sets with small cardinality,  in {\sl Approximation
Theory IV}, C.~K.~Chui, L.~L.~Schumaker, J.~D.~Ward, Eds., 521--526,
Academic Press, New York.

M.~S.~Henry and J.~A.~Roulier [1977], Uniform Lipschitz constants on
small intervals, {\sl J. Approx. Theory} {\bf 21}, 224--235.

M.~S.~Henry and J.~A.~Roulier [1978], Lipschitz and strong unicity
constants for changing dimension, {\sl J. Approx. Theory} {\bf 22},
85--94.

M.~S.~Henry, D.~P.~Schmidt and J.~J.~Swetits [1981], Uniform strong
unicity for rational approximation, {\sl J. Approx. Theory} {\bf
33}, 131--146.

M.~S.~Henry and J.~J.~Swetits [1980a], Growth rates for strong
unicity constants, in {\sl Approximation Theory III},  E.~W.~Cheney,
Ed., 501--505, Academic Press, New York.

M.~S.~Henry and J.~J.~Swetits [1980b], Lebesgue and strong unicity
constants, in {\sl Approximation Theory III}, E.~W.~Cheney, Ed.,
507--512, Academic Press, New York.

M.~S.~Henry and J.~J.~Swetits [1981], Precise orders of strong
unicity constants for a class of rational functions, {\sl J. Approx.
Theory} {\bf 32}, 292--305.

M.~S.~Henry and J.~J.~Swetits [1982], Lebesgue and strong unicity
constants for Zolotareff polynomials, {\sl Rocky Mountain J. Math.}
{\bf 12}, 547--556.

M.~S.~Henry and J.~J.~Swetits [1984], Limits of strong unicity
constants for certain $C^\infty$ functions, {\sl Acta. Math.
Hungar.} {\bf 43}, 309--323.

M.~S.~Henry, J.~J.~Swetits, and S.~E.~Weinstein [1980], Lebesgue and
strong unicity constants, in {\sl Approximation Theory III},
E.~W.~Cheney, Ed., 507--512, Academic Press, New York.

M.~S.~Henry, J.~J.~Swetits, and S.~E.~Weinstein [1981], Orders of
strong unicity constants,  {\sl J. Approx. Theory} {\bf 31},
175--187.

M.~S.~Henry, J.~J.~Swetits, and S.~E.~Weinstein [1983], On extremal
sets and strong unicity constants for certain $C\sp{\infty }$
functions, {\sl J. Approx. Theory} {\bf 37}, 155--174.

R.~Holmes and B.~Kripke [1968], Smoothness of approximation, {\sl
Mich. Math. J.} {\bf 15}, 225--248.

R.~Huotari and S.~Sahab [1994], Strong unicity versus modulus of
convexity, {\sl Bull. Austral. Math. Soc.} {\bf 49}, 305--310.

R.~C.~James [1947], Orthogonality and linear functionals in normed
linear spaces, {\sl Trans. Amer. Math.~Soc.} {\bf 61}, 265--292.

K.~Jittorntrum and M.~R.~Osborne [1980], Strong uniqueness and
second order convergence in nonlinear discrete approximation, {\sl
Numer. Math.} {\bf 34}, 439--455.

P.~Kirchberger, [1902] {\sl \"Uber Tchebychefsche
Ann\"aherungsmethoden}, Dissertation. Univ. G\"ott\-ingen.

A.~V.~Kolushov [1981], Differentiability of the operator of best
approximation, {\sl Math. Notes Acad. Sci. USSR} {\bf 29}, 295--306.

V.~V.~Kovtunec [1984], The Lipschitz property of the operator of
best approximation of complex-valued functions, in {\sl The Theory
of Approximation of Functions and its Applications}, 80--86, {\bf
137}, Akad.~Nauk Ukrain.~SSR, Inst.~Mat., Kiev.

B.~Kripke [1964], Best approximation with respect to nearby norms,
{\sl Numer. Math.} {\bf 6}, 103--105.

A.~Kro\'o [1977a], The continuity of best approximations, {\sl Acta.
Math. Hungar.} {\bf 30}, 175--188.

A.~Kro\'o [1977b], Differential properties of the operator of best
approximation, {\sl Acta. Math. Hungar.} {\bf 30}, 319--331.

A.~Kro\'o [1978], On the continuity of best approximation in the
space of integrable functions, {\sl Acta. Math. Hungar.} {\bf 32},
331--348.

A.~Kro\'o [1980], The Lipschitz constant of the operator of best
approximation, {\sl Acta. Math. Hungar.} {\bf 35}, 279--292.

A.~Kro\'o [1981a], On strong unicity of $L\sb{1}$-approximation,
{\sl Proc. Amer. Math. Soc.} {\bf 83}, 725--729.

A. Kro\'o [1981b], Best $L_1$-approximation on finite point sets:
rate of convergence, {\sl J. Approx. Theory} {\bf 33}, 340--352.

A.~Kro\'o [1981/82], On strong unicity of best approximation in
$C(R)$, {\sl Numer. Funct. Anal. Optim.} {\bf 4}, 437--443.

A.~Kro\'o [1983a], On the strong unicity of best Chebyshev
approximation of differentiable functions, {\sl Proc. Amer. Math.
Soc.} {\bf 89}, 611--617.

A.~Kro\'o [1983b], On unicity and strong unicity of best
approximation in the $L\sb{1}$-norm, {\sl Constructive Function
Theory '81}, (Varna, 1981), 396--399, Bulgar. Acad. Sci., Sofia.

A.~Kro\'o [1984], On the unicity of best Chebyshev approximation of
differentiable functions, {\sl Acta Sci. Math. (Szeged)} {\bf 47},
377--389.

A.~Kro\'o and F.~Peherstorfer [2007], Asymptotic representation of
$L_p$-minimal polynomials, $1<p<\infty$, {\sl Constr. Approx.} {\bf
25}, 29--39.

A.~Kro\'o and F.~Peherstorfer [2008], Asymptotic representation of
weighted $L_{\infty}$- and $L_1$-minimal polynomials, {\sl Math.
Proc. Camb. Phil. Soc.} {\bf 144}, 241--254.

A.~Kro\'o and D.~Schmidt [1991], A Haar-type theory of best uniform
approximation with constraints, {\sl Acta Math.~Hung.} {\bf 58},
351--374.

A.~Kro\'o, M.~Sommer, and  H.~Strauss [1989], On strong uniqueness
in one-sided $L^1$-approx\-imation of differentiable functions, {\sl
Proc. Amer. Math. Soc.} {\bf 106}, 1011--1016.

K.~Kuratowski [1966], {\sl Topology}, Academic Press, New York.

P.-J.~Laurent and D.~Pai, On simultaneous approximation, {\sl Numer.
Funct. Anal. Optim.} {\bf 19}, 1045--1064.

C.~Li [2003], On best uniform restricted range approximation in
complex-valued continuous function spaces, {\sl J. Approx. Theory}
{\bf 120}, 71--84.

C.~Li and G.~A.~Watson [1994], On approximation using a peak norm,
{\sl J. Approx. Theory} {\bf 77}, 266--275.

C.~Li and G.~A.~Watson [1997], Strong uniqueness in restricted
rational approximation, {\sl J. Approx. Theory} {\bf 89}, 96--113.

C.~Li and G.~A.~Watson [1999], On approximation using a generalized
peak norm, {\sl Commun. Appl. Anal.} {\bf 3}, 357--371.

W.~Li [1989], Strong uniqueness and Lipschitz continuity of metric
projections: a generalization of the classical Haar theory, {\sl J.
Approx. Theory} {\bf 56}, 164--184.

P.~K.~Lin [1989], Strongly unique best approximation in uniformly
convex Banach spaces, {\sl J. Approx. Theory} {\bf 56}, 101--107.

H.~L.~Loeb [1966], Approximation by generalized rationals, {\sl J.
SIAM Numer. Anal.} {\bf 3}, 34--55.

D.~S.~Lubinsky and E.~B.~Saff [1987], Strong asymptotics for
$L_p$-minimal polynomials, $1<p<\infty$, in {\sl Approximation
Theory, Tampa}, E.~B.~Saff, Ed., LNM {\bf 1287}, 83--104.

Z.~Ma [1991], Some problems on a variation of $L_1$ approximation,
in {\sl Progress in Approximation Theory}, P.~Nevai, A.~Pinkus,
Eds., 667--692, Academic Press, Boston.

Z.~W.~Ma and Y.~G.~Shi [1990], A variation of rational $L_1$
approximation, {\sl J. Approx. Theory} {\bf 62}, 262--273.

H.~Maehly and Ch.~Witzgall [1960],  Tschebyscheff-Approximationen in
kleinen Intervallen. I. Approximation durch Polynome, {\sl Numer.
Math.} {\bf 2}, 142--150.

P.~F.~Mah [1984], Strong uniqueness in nonlinear approximation, {\sl
J. Approx. Theory} {\bf  41}, 91--99.

L.~K.~Maloz\"emova [1990], On the strong uniqueness constant in a
discrete problem of best approximation in the mean, {\sl Vestnik
Leningrad Univ. Mat. Mekh. Astronom.}  {\bf 125}, 26--30; in English
translation in {\sl Vestnik Leningrad Univ. Math.} {\bf 23}, 33--37.

A.~V.~Marinov [1983], Strong uniqueness constants for best uniform
approximations on compacta, {\sl Mat. Zametki} {\bf 34}, 31--46; in
English translation in {\sl Math. Notes} {\bf 34}, 499--507.

O.~M.~Martinov [2002], Constants of strong unicity of minimal
projections onto some two-dimensional subspaces of $l^{(4)}_\infty$,
{\sl J. Approx. Theory} {\bf  118}, 175--187.

H.~W.~McLaughlin and  K.~B.~Somers [1975], Another characterization
of Haar subspaces, {\sl J. Approx. Theory} {\bf 14}, 93--102.

A.~Meir [1984], On the uniform convexity of $L^p$ spaces, $1<p\le
2$, {\sl Illinois J. Math.} {\bf 28}, 420--424.

V.~N.~Nikolskii [1982], Boundary sets in the strong sense and strong
uniqueness of elements of best approximation, {\sl Application of
Functional Analysis in Approximation Theory}, 126--130, Kalinin.
Gos. Univ., Kalinin, (in Russian).

D.~J.~Newman and H. S. Shapiro [1963], Some theorems on Cebysev
approximation, {\sl Duke Math. J.} {\bf 30}, 673--682. Abstract
appeared in Notices in 1962 {\bf 9}, 143.

G.~N\"urnberger [1979], Unicity and strong unicity in approximation
theory, {\sl J. Approx. Theory} {\bf 26}, 54--70.

G.~N\"urnberger [1980], Strong uniqueness of best approximations and
weak Chebyshev systems, in  {\sl Quantitative Approximation} (Proc.
Internat. Sympos., Bonn, 1979),  pp. 255--266, Academic Press, New
York-London.

G.~N\"urnberger [1982], A local version of Haar's theorem in
approximation theory, {\sl Numer. Funct. Anal. Optim.} {\bf 5},
21--46.

G.~N\"urnberger [1982/83], Strong unicity constants for spline
functions, {\sl Numer. Funct. Anal. Optim.} {\bf 5}, 319--347.

G.~N\"urnberger [1983], Strong unicity constants for
finite-dimensional subspaces, in {\sl Approximation Theory IV},
C.~K.Chui, L.~L.~Schumaker, J.~D.~Ward, Eds., 643--648, Academic
Press, New York.

G.~N\"urnberger [1984],  Strong unicity of best approximations: a
numerical aspect, {\sl Numer. Funct. Anal. Optim.} {\bf 6} (1983),
399--421.

G.~N\"urnberger [1985a], Unicity in one-sided $L_1$-approximation
and quadrature formulae, {\sl J. Approx. Theory} {\bf 45}, 271--279.

G.~N\"urnberger [1985b],  Global unicity in semi-infinite
optimization, {\sl Numer. Funct. Anal. Optim.} {\bf 8}, 173--191.

G.~N\"urnberger [1987], Strong unicity constants in Chebyshev
approximation, in {\sl Numerical Methods of Approximation Theory},
Vol. 8 (Oberwolfach, 1986),  144--154, ISNM, {\bf 81}, Birkh\"auser,
Basel.

G.~N\"urnberger [1994], Strong unicity in nonlinear approximation
and free knot splines, {\sl Constr. Approx.} {\bf 10}, 285--299.

G.~N\"urnberger and I.~Singer [1982], Uniqueness and strong
uniqueness of best approximations by spline subspaces and other
subspaces, {\sl J. Math. Anal. Appl.} {\bf 90}, 171--184.

P.~L.~Papini [1978], Approximation and strong approximation in
normed spaces via tangent functionals, {\sl J. Approx. Theory} {\bf
22}, 111--118.

S.~H.~Park [1989], Uniform Hausdorff strong uniqueness, {\sl J.
Approx. Theory} {\bf 58}, 78--89.

S.~O.~Paur and J.~A.~Roulier [1980], Uniform Lipschitz and strong
unicity constants on subintervals, in {\sl Approximation Theory
III}, E.~W.~Cheney, Ed., 715--720, Academic Press, New York.

S.~O.~Paur and J.~A.~Roulier [1981], Continuity and strong unicity
of the best approximation operator on subintervals, {\sl J. Approx.
Theory} {\bf 32}, 247--255.

J.~Peetre [1970], Approximation of norms, {\sl J. Approx. Theory}
{\bf 3}, 243--260.

R.~R.~Phelps [1966], Cebysev subspaces of finite dimension in $L_1$,
{\sl Proc.~Amer.~Math.~Soc.} {\bf 17}, 646--652.

A.~M.~Pinkus [1989], {\sl On $L^1$-Approximation}, Cambridge
University Press, Cambridge.

A.~Pinkus and H.~Strauss [1987], One-sided $L^1$-approximation to
differentiable functions, {\sl Approx.~Theory Appl.} {\bf 3},
81--96.

S.~J.~Poreda [1976], Counterexamples in best approximation, {\sl
Proc. Amer. Math. Soc.} {\bf 56}, 167--171.

B.~Prus and R.~Smarzewski [1987], Strongly unique best approximation
and centers in uniformly convex spaces, {\sl J.~Math.~Anal.~Appl.}
{\bf 121}, 10--21.

T.~J.~Rivlin [1969], {\sl An Introduction to the Approximation of
Functions}, Blaisdell, Waltham, Mass.

T.~J.~Rivlin [1984a], The strong uniqueness constant in complex
approximation, in {\sl Rational Approximation and Interpolation},
Tampa, 1983, P.~R.~Graves-Morris, E.~B.~Saff, R.~S.~Varga, Eds.,
145--149, LNM {\bf 1105}, Springer.

T.~J.~Rivlin [1984b], The best strong uniqueness constant for a
multivariate Chebyshev polynomial, {\sl J. Approx. Theory} {\bf 41},
56--63.

E.~Rozema [1974], Almost Chebyshev subspaces of $L^1(\mu; E)$, {\sl
Pacific J. Math.} {\bf 53}, 585--604.

R.~Schaback [1978], On alternation numbers in nonlinear Chebyshev
approximation, {\sl J. Approx. Theory} {\bf 23}, 379--391.

D.~P.~Schmidt [1978], On an unboundedness conjecture for strong
unicity constants, {\sl J. Approx. Theory} {\bf 24}, 216--223.

D.~Schmidt [1979], Strong unicity and Lipschitz conditions of order
${1\over 2}$ for monotone approximation, {\sl J. Approx. Theory}
{\bf 27}, 346--354.

D.~Schmidt [1980a], A characterization of strong unicity constants,
in {\sl Approximation Theory III },  E.~W.~Cheney, Ed., 805--810,
Academic Press, New York.

D.~Schmidt [1980b], Strong uniqueness for Chebyshev approximation by
reciprocals of polynomials on $[0,\,\infty ]$, {\sl J. Approx.
Theory} {\bf 30}, 277--283.

N.~Ph.~Seif and G.~D.~Taylor [1982], Copositive rational
approximation, {\sl J. Approx. Theory} {\bf 35}, 225--242.

Y.~G.~Shi [1981], Weighted simultaneous Chebyshev approximation,
{\sl J. Approx. Theory} {\bf 32}, 306--315.

R.~Smarzewski [1986a], Strongly unique best approximation in Banach
spaces, {\sl J. Approx. Theory} {\bf 46}, 184--194.

R.~Smarzewski [1986b], Strongly unique minimization of functionals
in Banach spaces with applications to theory of approximation and
fixed points, {\sl J. Math. Anal. Appl.} {\bf 115}, 155--172.

R.~Smarzewski [1987a], Strong unicity in nonlinear approximation, in
{\sl Rational Approximation and Applications in Mathematics and
Physics}, (Lancut, 1985),  331--350, LNM {\bf 1237}, Springer,
Berlin.

R.~Smarzewski [1987b], Strongly unique best approximation in Banach
spaces II, {\sl J. Approx. Theory} {\bf 51}, 202--217.

R.~Smarzewski [1988a], Strong uniqueness of best approximations in
an abstract $L^1$ space, {\sl J. Math. Anal. Appl.} {\bf 136},
347--351.

R.~Smarzewski [1988b],  Strong unicity of best approximations in
$L\sb \infty(S,\Sigma,µ)$,  {\sl Proc. Amer. Math. Soc.} {\bf 103},
113--116.

R.~Smarzewski [1989], Strong unicity of order $2$ in $C(T)$, {\sl J.
Approx. Theory} {\bf 56}, 306--315.

R.~Smarzewski [1990], Finite extremal characterization of strong
uniqueness in normed spaces, {\sl J. Approx. Theory} {\bf 62},
213--222.

M.~Sommer and H.~Strauss [1993], Order of strong uniqueness in best
$L_\infty$-approximation by spline spaces, {\sl Acta Math. Hungar.}
{\bf 61}, 259--280.

H.~Strauss [1982], Unicity of best one-sided $L_1$-approximations,
{\sl Numer. Math.} {\bf 40}, 229--243.

H.~Strauss [1992], Uniform reciprocal approximation subject to
coefficient constraints, {\sl Approx. Theory Appl.} {\bf 8},
89--102.

J.~Sudolski and A.~W\'ojcik [1987], Another generalization of strong
unicity, {\sl Univ. Iagel. Acta Math.} {\bf 26}, 43--51.

J.~Sudolski and A.~P.~W\'ojcik [1990], Some remarks on strong
uniqueness of best approximation, {\sl Approx. Theory Appl.} {\bf
6}, 44--78.

S.~Tanimoto [1998], On best simultaneous approximation, {\sl Math.
Japon.} {\bf 48}, 275--279.

J.~L.~Walsh [1931], The existence of rational functions of best
approximation, {\sl Trans. Amer. Math. Soc.} {\bf 33}, 668--689.

G.~A.~Watson [1980], {\sl Approximation Theory and Numerical
Methods}, J.~Wiley and Sons, Chichester.

R.~Wegmann [1975], Bounds for nearly best approximations, {\sl
Proc.~Amer.~Math.~Soc.} {\bf 52}, 252--256.

H.~Werner [1964], On the rational Tschebyscheff operator, {\sl Math.
Z.} {\bf 86}, 317--326.

A.~W\'ojcik [1981], Characterization of strong unicity by tangent
cones, in {\sl Approximation and Function Spaces} (Gdansk, 1979),
854--866, North-Holland, Amsterdam-New York.

D.~Wulbert [1971], Uniqueness and differential characterizations of
approximation from manifolds of functions, {\sl Amer.~J.~Math.} {\bf
93}, 350--366.

S.~S.~Xu [1995], A note on a strong uniqueness theorem of Strauss,
{\sl Approx. Theory Appl. (N.S.)} {\bf 11}, 1--5.

C.~Yang [1993], Uniqueness of best approximation with coefficient
constraints, {\sl Rocky Mountain J. Math.} {\bf 23}, 1123--1132.

W.~S.~Yang and C.~Li [1994],  Strong unicity for monotone
approximation by reciprocals of polynomials, {\sl J. Approx. Theory}
{\bf 78}, 19--29.

J.~W.~Young [1907],  General theory of approximation by functions
involving a given number of arbitrary parameters, {\sl Trans. Amer.
Math. Soc.} {\bf 8}, 331--344.

F.~Zeilfelder [1999], Strong unicity of best uniform approximations
from periodic spline spaces, {\sl J. Approx. Theory} {\bf 99},
1--29.

D.~Zwick [1987], Strong uniqueness of best spline approximation for
a class of piecewise $n$-convex functions, {\sl Numer. Funct. Anal.
Optim.} {\bf 9}, 371--379.


{

\bigskip\obeylines
Andr\'as Kro\'o Alfred R\'enyi Institute of Mathematics Hungarian
Academy of Sciences Budapest, Hungary and Budapest University of
Technology and Economics Department of Analysis Budapest, Hungary
{\tt kroo@renyi.hu}

\bigskip
Allan Pinkus Department of Mathematics Technion Haifa, Israel {\tt
pinkus@tx.technion.ac.il}

}

\bye